\documentclass[a4paper,11pt]{article}
\small 
\normalsize

\usepackage[utf8]{inputenc}
\usepackage[T1]{fontenc}
\usepackage[english]{babel}
\usepackage{times}
\usepackage{ulem}

\usepackage[pdftex]{hyperref}

\usepackage{theorem}

\usepackage[all]{xy}

\usepackage{stmaryrd}
\usepackage{graphicx}
\usepackage{paralist}
\usepackage{amsfonts}
\usepackage[leqno]{amsmath}
\usepackage{amssymb,mathrsfs}
\usepackage{color}
\usepackage{times}
\usepackage{enumerate,vmargin}
\usepackage{tocvsec2}
\usepackage{verbatim}
\usepackage[labelfont=bf,font=it]{caption}
\usepackage[numbers,sort&compress]{natbib}

\setmarginsrb{2.9cm}{2.6cm}{2.9cm}{1.7cm}{0cm}{0mm}{0cm}{9mm}

\hypersetup{
    unicode      = false,     
    pdftoolbar   = true,      
    pdfmenubar   = true,      
    pdffitwindow = true,      
    pdfnewwindow = true,      
    colorlinks   = true,      
    linkcolor    = blue,      
    citecolor    = red,      
    filecolor    = blue,      
    urlcolor     = green       
}

\definecolor{vio}{rgb}{.5,.1,.5}
\definecolor{grey}{rgb}{.5,.5,.5}
\definecolor{cnic}{rgb}{.9,0.5,0}
\definecolor{mycolor}{rgb}{.6,0,.8}
\definecolor{mm-color}{rgb}{.6,0,.8}
\definecolor{darkred}{RGB}{165,42,42}

\def \mm#1{{\color{mycolor}#1\color{black}}}

\def \ts#1{{\color{darkred}\sout{#1}\color{black}}}

\newtheorem{lem}{Lemma}[section]
\newtheorem{thm}[lem]{Theorem}
\newtheorem{prop}[lem]{Proposition}

{\theorembodyfont{\rmfamily}}
{\theorembodyfont{\rmfamily}}
{\theorembodyfont{\rmfamily}}
{\theorembodyfont{\rmfamily}\newtheorem{rem}{Remark}[section]}
{\theorembodyfont{\rmfamily}}

\newcommand{\noi}{\noindent}

\newcommand{\un}{{\bf 1}}

\newcommand{\cA}{\mathcal A}

\newcommand{\cE}{\mathcal E}

\newcommand{\cI}{\mathcal I}

\newcommand{\cM}{\mathcal M}
\newcommand{\cN}{\mathcal N}

\newcommand{\cP}{\mathcal P}
\newcommand{\cQ}{\mathcal Q}

\newcommand{\cS}{\mathcal S}

\newcommand{\cV}{\mathcal V}

\newcommand{\bE}{\mathbf E}

\newcommand{\bG}{\mathbf G}
\newcommand{\bP}{\mathbf P}

\newcommand{\bN}{\mathbf N}

\newcommand{\bm}{\mathbf m}
\newcommand{\bp}{\mathbf p}

\newcommand{\bzeta}{\boldsymbol \zeta}

\newcommand{\bbQ}{\mathbb Q}

\newcommand{\bbN}{\mathbb N}

\newcommand{\bbR}{\mathbb R}

\newcommand{\m}{\mathrm m}
\newcommand{\ccF}{\mathscr F}
\newcommand{\ccG}{\mathscr G}
\newcommand{\ccE}{\mathscr E}

\newcommand{\elldo}{{\ell}^{_{\, \downarrow}}}

\def\cq{$\hfill \square$}
\def\cqfd{$\hfill \blacksquare$}
\def\ino{ \! \in \! }

\def\bC{\mathbf{C}}
\def\bD{\mathbf{D}}

\def\bT{\mathbf{T}}

\newcommand{\btheta}{{\boldsymbol{\theta}}}
 
\def\cJ{\mathcal{J}}
\def\cV{\mathscr{V}}
\def\cH{\mathcal{H}}
\def\cK{\mathcal{K}}
\def\ccX{\mathscr{X}}
\def\ccY{\mathscr{Y}}
\def\ccH{\mathscr{H}}
\def\ccL{\mathscr{L}}
\def\ccZ{\mathscr{Z}}

\def\subo{\gamma}

\newcommand{\Ptt}{{\boldsymbol{\Pi}}} 
\newcommand{\Htt}{{\boldsymbol{\mathtt{H}}}} 
\newcommand{\Ytt}{{\boldsymbol{\mathtt{Y}}}}

\newcommand{\cG}{\boldsymbol{\mathcal G}}
\newcommand{\cT}{\boldsymbol{\mathcal T}}
\newcommand{\bw}{\mathtt w}
\newcommand{\JJ}{J}
\newcommand{\Jtt}{{\mathtt J}}
\newcommand{\bq}{\mathbf{q}}

\def\bs{\mathbf{s}}
\def\rbm{\mathtt{m}}

\def\epp{\varepsilon}
\usepackage[refpage,noprefix]{nomencl}                          
\makeglossary
\makenomenclature                                           

\graphicspath{{.}{figures/}}

\title{ \textsc{Limits of multiplicative inhomogeneous random graphs and L\'evy trees: The continuum graphs 
}}
\date{}
\author{Nicolas \textsc{Broutin}
\thanks{Sorbonne Universit\'e, Campus Pierre et Marie Curie, 
Case courrier 158, 4 place Jussieu,  
75252 Paris Cedex 05,  
France. Email: nicolas.broutin@upmc.fr} 
\and Thomas \textsc{Duquesne}
\thanks{Sorbonne Universit\'e, Campus Pierre et Marie Curie, 
Case courrier 158, 4 place Jussieu,  
75252 Paris Cedex 05,  
France.
Email: thomas.duquesne@upmc.fr}
\and Minmin \textsc{Wang}
\thanks{University of Sussex, Department of Mathematics, Falmer, 
Brighton, BN1 9QH,  United Kingdom.  
Email: minmin.wang@sussex.ac.uk }
}

\begin{document}

\maketitle

\begin{abstract}
Motivated by limits of critical inhomogeneous random graphs, we construct a family of sequences of measured metric spaces that we call \emph{continuous multiplicative graphs}, that are expected to be the universal limit of graphs related to the multiplicative coalescent (the Erd\H{o}s--Rényi random graph, more generally the so-called rank-one inhomogeneous random graphs of various types, and the configuration model). At the discrete level, the construction relies on a new point of view on (discrete) inhomogeneous random graphs that involves an embedding into a Galton--Watson forest. 
The new representation allows us to demonstrate that a processus that was already present in the pionnering work of Aldous [\emph{Ann.\ Probab.}, vol.~25, pp.~812--854, 1997] and Aldous and Limic [\emph{Electron. J. Probab.}, vol.~3, pp.~1--59, 1998] about the multiplicative coalescent actually also (essentially) encodes the limiting metric: The discrete embedding of random graphs into a Galton--Watson forest is paralleled by an embedding of the encoding process into a Lévy process which is crucial in proving the very existence of the local time functionals on which the metric is based; it also yields a transparent approach to compactness and fractal dimensions of the continuous objects. In a companion paper, we show that the continuous Lévy graphs are indeed the scaling limit of inhomogeneous random graphs.
\end{abstract}

\section{Introduction}


\noi
\textbf{Motivation and results}. In this paper, we construct a family of sequences of measured metric spaces that we call \emph{continuous multiplicative graphs}. A companion paper \cite{BDW2} is dedicated to showing that these objects are indeed the scaling limits of critical (rank-one) inhomogeneous random graphs. More generally, continuous multiplicative graphs are expected to be the scaling limits of all of models of graphs that are related to the multiplicative coalescent \cite{Al97,AlLi98} which includes the classical Erd\H{o}s--Rényi random graphs (as a special case of rank-one inhomogeneous random graphs) and the configuration model.

The (metric space) scaling limits for critical random graphs have been constructed and limit theorems have been proved for a number of models. At level of discrete objects, the general strategy consists in designing a suitable exploration of the graphs that yields a decomposition into a random function describing a spanning forest, together with relatively few additional ``decorations'' indicating where to add additional edges in order to obtain the graph. By ``suitable'', we mean that the random function should both encode the metric relatively transparently, and be tractable when the time comes to taking limits. So, the random function should take advantage of the independence underlying the random graph construction, and for example should whenever possible remain close to some kind of random walk. Except in very specific cases, this constraint implies that the metric of the spanning forest cannot be a very simple functional of the random function. The natural representations, which go back to the treatment of scaling limits, encode the metric as a (discrete) local-time functional of the random function. 

For a general class of random graphs related to the multiplicative coalescent, Aldous \cite{Al97} and Aldous \& Limic \cite{AlLi98} have devised natural explorations that are suitable for taking limits, and that allowed them to study the scaling limits of the sequence of sizes of the connected components\footnote{to be precise, \cite{AlLi98} deals with ``masses'' for a suitable measure on the vertices that is not the counting measure.}. However, it remains so far unclear how to deal with the metric structure from that approach for two reasons. First, the exploration does not directly yield a nice connection with the metric. Sometimes, a slight modification of the exploration gives access to the metric without changing the law of the random function; this is for instance the case for the Erd\H{o}s--Rényi random graphs: the breadth-first exploration used by Aldous is not convenient for the analysis of the metric, but the depth-first version used in \cite{AdBrGo12a} has the same law. One could envision that such a slight modification of the general exploration of Aldous \& Limic, say from a breadth-first to a depth-first approach, would both give a handle on the metric and leave the processes unchanged in distribution (and thus suitable for taking limits). However, in general, the processes in \cite{AlLi98} are only semi-martingales for which the relevant continuum local-time functionals are now known to exist.

The present paper addresses the issues of the previous paragraph. At the discrete level of finite graphs, we choose a very specific model that is most convenient for the construction of the limit objects. The graphs $G= (\cV(G), \ccE (G))$ that we consider are \textit{not oriented}, without either loops or multiple edges: $\ccE(G)$ is therefore a set consisting of 
unordered pairs of distinct vertices. Let $n\! \geq \! 2$ and let $\bw \! = \! (w_1, \ldots , w_n)$ be a set of vertex \textit{weights}: namely, it is a set of positive real numbers such that 
$ w_1 \! \geq \! w_2 \! \geq \! \ldots \! \geq \! w_n \! >\! 0$. For all $r\ino (0,\infty)$, we let $\sigma_r (\bw)= w^r_1+ \ldots + w^r_n$. 
Then, the random graph $\cG_\bw$ is said to be \textit{$\bw$-multiplicative} if $\cV (\cG_\bw) \! =\!  \{ 1, \ldots, n\}$ and if the random variables $(\un_{\{ \{ i,j \}\in \ccE (\cG_\bw) \}})_{1\leq  i <  j  \leq n}$, are independent and for each $1\le i<j\le n$ we have
\[
\bP \big( \{ i,j\}\ino \ccE (\cG_\bw)\big)=1\! -\! 
e^{-w_iw_j/\sigma_1 (\bw)} . 
\]
We call this graph \textit{multiplicative} because as observed by Aldous, its connected components provide a natural representation of the multiplicative coalescent: see Aldous \cite{Al97} and Aldous \& Limic \cite{AlLi98} for more details.

The first result of the paper (Theorem \ref{loiGG}) shows that $\cG_{\bw}$ is coded by a natural depth-first exploration process that gives access to the metric and that is explained in terms of the following queueing system: \textit{a single server serves at most one client at a time applying the 
Last In First Out policy (LIFO, 
for short); exactly $n$ clients will enter the queue and each client is labelled with a distinct integer of $\{ 1, \ldots , n\}$; Client $i$  enters the queue at a time $E_i$ and requires a time of service $w_i$; we assume that the 
$E_j$ are independent and exponentially distributed r.v.~such that $\bE [E_j] \! =\! \sigma_1 (\bw)/ w_j$}. The LIFO-queue yields the following tree $\cT_{\!\!\bw}$ whose vertices are the clients: \textit{the server is the root (Client $0$) and Client $j$ is a child of Client $i$ in $\cT_{\! \! \bw}$ if and only if Client $j$ interrupts the service of Client $i$ (or arrives when the server is idle if $i\! =\!  0$)}. We claim that the subtrees of $\cT_{\!\!  \bw}$ that stem from the root are spanning trees of the connected components of $\cG_{\!\bw}$ (and thus, $\cT_{\!\!  \bw}$ captures a large part of the metric of $\cG_{\bw}$). If one introduces  
\begin{equation*}
 Y^{\bw}_t= \! - t + \!\! \sum_{1\leq i\leq n} \!\! w_i \un_{\{ E_i \leq t\}}, \quad J^\bw_t\! = \! \! \inf_{s\in [0, t]} \!\!  Y^\bw_s \quad \textrm{and} \quad \cH^\bw_t = \# \Big\{ s\in [0, t] : \inf_{r\in [s, t]} Y^\bw_r > Y^\bw_{s-} \Big\} \; , 
 \end{equation*}
then, $Y^\bw_t\! -\! J^\bw_t$ is the load of the server (i.e.~the amount of service due at time $t$) and $\cH^\bw_t$ is the number of clients waiting in the queue at time $t$ (see Figure \ref{fig:LIFO}). The process $\cH^\bw$ is a contour (or a depth-first exploration) of $\cT_{\!\! \bw}$ and $\cH^\bw$ codes its graph-metric: 
namely, the distance between the vertices/clients served at times $s$ and $t$ in $\cT_{\!\!\bw}$ is $\cH^\bw_t+\cH^\bw_s \! - \!  2 \!  \min_{r \in [s  \wedge t , s\vee t] }   \cH^\bw_r $. 
We obtain $\cG_{\! \bw}$ by adding to $\cT_{\!\! \bw} \backslash \{ 0\}$ a set of surplus edges $\mathcal{S}_{\bw}$ that is derived from
a Poisson point process on $[0,\infty)^2$ of intensity $\frac{{_1}}{{^{\sigma_1 (\bw)}}}  \un_{\{ 0< y <Y^\bw_t -\JJ^\bw_t \}} \, dt\, dy$ as follows: an atom $(t,y)$ of this Poisson point process corresponds to the surplus edge $\{ i,j\}$, where $j$ is the client served at time $t$ and $i$ is the client served at time $\inf \{ s\ino [0, t]: \inf_{r\in [s,t]}Y_r^\bw\! -\! J_r^\bw \! >\! y\}$. Then, Theorem \ref{loiGG} asserts that the resulting graph $\cG_{\bw}\! =\! ( \cT_{\!\! \bw} \backslash \{ 0\}) \cup \mathcal{S}_\bw$ is a $\bw$-multiplicative graph. 
%
%
%
%
%

In order to define the scaling limits of the previous processes and of $\bw$-multiplicative graphs, 
the second main idea of the paper consists in embedding $\cG_\bw$ into a Galton-Watson tree by actually embedding the queue governed by $Y^\bw$ into a Markovian queue that is defined as  follows: \textit{a single server receives in total an infinite number of clients; it applies the LIFO policy; 
clients arrive at unit rate; each client has a type that is an integer ranging in $\{ 1, \ldots, n\}$; the amount of service required by a client of type $j$ is $w_j$; types are i.i.d.~with law $\nu_\bw \! = \! \frac{1}{\sigma_1 (\bw)}\sum_{1 \leq j\leq n} w_j \delta_{j} $.}
If $\tau_k$ and $\Jtt_k$ stand for resp.~the arrival-time and the type of the 
$k$-th client, then, the Markovian LIFO queueing system is entirely characterised by $\sum_{k\geq 1} \delta_{(\tau_k , \Jtt_k)}$ that is a Poisson point measure on $[0, \infty) \! \times \! \{ 1, \ldots, n\}$ with intensity $\ell \! \otimes \! \nu_\bw$, where $\ell $ stands for the Lebesgue measure on $[0, \infty)$. The Markovian queue yields a tree $\bT_{\! \bw}$ that is defined as follows: 
\textit{the server is the root of $\bT_{\! \bw}$ and the $k$-th client to enter the queue is a child of the $l$-th one if the $k$-th client enters when the $l$-th client is being served}. To simplify our explanations, let us focus here on the (sub)critical cases where 
$\sigma_2 (\bw) \! \leq \! \sigma_1 (\bw)$. Then, $\bT_{\! \bw} $ is a sequence of i.i.d.~Galton--Watson trees glued at their root and whose offspring distribution is 
$\mu_\bw(k)\! =\!   \sum_{1\leq j\leq n} w_j^{k+1} \exp (-w_j) /(\sigma_1 (\bw) k!)$, $k\ino \bbN$, that is (sub)critical. 
The tree $\bT_{\! \bw}$ is then coded by its contour process $(H^\bw_t)_{t\in [0, \infty)}$: namely, $H^\bw_t$ stands for the number of clients waiting in the Markovian queue at time $t$ and it is given by 
 \begin{equation}
 \label{HXdisdef}
 H^\bw_t = \# \Big\{ 
 s\in [0, t] : \inf_{r\in [s, t]} X^\bw_r > X^\bw_{s-} \Big\}, \quad \textrm{where} \quad  X^\bw_{t}  =  -t + \sum_{k\geq 1} w_{\Jtt_k}\un_{[0, t]} (\tau_k) 
 \end{equation}
is the (algebraic) load of the Markovian server: namely $X^\bw_t \! -\! \inf_{s\in [0, t]}X_s^\bw$ is the amount of service due at time $t$.  
Note that $X^\bw$ is a compound Poisson process and let us mention that the possible scaling limits 
of $(X^\bw, H^\bw)$ are well-understood. 

The queue governed by $Y^\bw$ is then obtained by pruning clients from the Markovian queue in the following way. We colour each client of the queue governed by $X^\bw$ in blue or in red according to the following rules: \textit{if the type $\Jtt_k$ of the $k$-th client already 
appeared among the types of the blue clients who previously entered the Markovian queue, then the $k$-th client is red; otherwise the $k$-th client inherits her/his colour from the colour of the client who is 
currently served when she/he arrives (and this colour is blue if there is no client served when she/he arrives: namely, we consider that the server is blue)}: see Figure \ref{fig:tree_pruning}. We claim (see Proposition \ref{Xwfrombr} and Lemma \ref{Hthetalem}) that if we skip the periods of time during which red clients are served in the Markovian queue, we get the queue governed by $Y^\bw$ that can be therefore viewed as the 
blue sub-queue. Namely, if $\mathtt{Blue} $ is the set of times $t$ when a blue client is served, then one sets $\theta_t^{\mathtt{b} , \bw} \! = \! \inf \{ s\ino [0, \infty):\!  \int_0^s\un_{\mathtt{Blue}} (r) dr > t \}$ and the above mentioned embedding is formally given by the fact that a.s.~for all $t\ino [0, \infty)$, 
$(Y^\bw_t, \cH^\bw_t)\! =\! (X^\bw (\theta_t^{\mathtt{b} , \bw}) , 
H^\bw (\theta_t^{\mathtt{b} , \bw}))$ (see Figure \ref{fig:proc_color}).

As proved in a companion paper \cite{BDW2}, the only possible scaling limits of processes like $Y^\bw$ that are relevant for our purpose are the processes already introduced in Aldous \& Limic  
\cite{AlLi98}. These processes $Y$ are parametrized by $(\alpha,\beta,\kappa,{\mathbf c})$, where $\alpha\in \bbR$, $\beta\ge 0$, $\kappa\ge 0$ and ${\mathbf c}= (c_i)_{j\geq 1}$ satisfies 
$c_{j} \! \geq \! c_{j+1}$, and $\sum_{j\geq 1} c_j^3\! <\! \infty$, and they are informally defined by  
\begin{equation}
\label{eq:Y_intro}
Y_t:=-\alpha t \! -\! \frac{_1}{^2}\kappa \beta t^2 + \sqrt{\beta} B_t +  \sum_{j\geq 1} \! c_j (\un_{\{ E_j \leq t \}} \! -\! c_j \kappa t) , \quad t\ino [0, \infty). 
\end{equation}
Here, $(B_t)_{t\in [0, \infty}$ is a standard linear Brownian motion, the $E_j$ are exponentially distributed with parameter $\kappa c_j$ and $B$ and the $E_j$ are independent (see (\ref{AetYdef}) and Remark \ref{Yrepres} for more details). To define the analogue of $\cH^\bw$, we proceed in an indirect way that relies on a continuous version of the embedding of $Y$ into a Markov process. More precisely, the scaling limits of $X^\bw$ are L\'evy processes $(X_t)_{t\in [0, \infty)}$ without negative jump, whose law is characterized by their Laplace exponent of the form 
$$\forall \lambda, t\ino [0, \infty), \quad \psi(\lambda)  \! = \! \log \bE [ \exp ( -\lambda X_t)]   \! = \!   
 \alpha \lambda +\frac{_{_1}}{^{^2}} \beta \lambda^2 +  \!   \sum_{j\geq 1}   \kappa c_j   \big( e^{-\lambda c_j}\! -\! 1\! + \! \lambda c_j \big). $$
To simplify the explanation we focus in this introduction on the (sub)critical cases where $\alpha \! \geq \! 0$. As proved in Le Gall \& Le Jan \cite{LGLJ98} and Le Gall \& D.~\cite{DuLG02}, if $\int^\infty d\lambda / \psi (\lambda) \! <\! \infty$, then there exists a continuous process $(H_t)_{t\in [0, \infty)}$ such that for all $t\ino [0, \infty)$, the following limit holds in probability: 
\begin{equation*}
H_t = \lim_{\varepsilon \rightarrow 0} \frac{1}{\varepsilon} \!  \int_0^{t} \!  \un_{\{  X_s - \inf_{r\in [s, t]} X_r  \leq \varepsilon\}} \, ds \; ,
\end{equation*}
that is a local time version of (\ref{XJHdef}). We refer to $H$ as to 
the height process of $X$ and $(X, H)$ is the continuous analogue of $(X^\bw, H^\bw)$. 
The third main result of the paper (Theorem \ref{Xdefthm} and Theorem \ref{cHdefthm}) 
asserts that, as in 
the discrete setting, there exists a time-change, namely an increasing positive process $\theta^{\mathtt{b}}$ and a continuous process $(\cH_t)_{t\in [0, \infty)}$ that is adapted with respect to $Y$ and such that a.s.~for all $t\ino [0, \infty)$, $(Y_t, \cH_t)\! = \! (X(\theta^{\mathtt{b}}_t) ,  H(\theta^{\mathtt{b}}_t) )$. 

Then continuous versions of multiplicative graphs are obtained as in the discrete setting: we first define $\cT$, the random continuum tree coded by $\cH$, namely $s,t\ino [0, \infty)$ correspond to points in $\cT$ that are at distance $\cH_t+\cH_s \! - \!  2 \!  \min_{r \in [s  \wedge t , s\vee t] }   \cH_r $ and to obtain the $(\alpha, \beta, \kappa, \mathbf{c})$-continuum multiplicative graph $\mathbf{G}$, we identify points in $\cT$ thanks to a Poisson point process on $[0, \infty)^2$ with intensity 
$\kappa \un_{\{ y< Y_t -\inf_{[0, t]} Y  \}}  dt \, dy$ as in the discrete setting: we refer to Sections \ref{codgrasec} and \ref{defICRG} for a precise definition. Specifically, our construction shows that $\mathbf{G}$ can be embedded into the tree coded by $H$, namely a $\psi$-L\'evy tree. It yields a transparent approach to the main geometric properties of 
$\mathbf{G}$: an explicit condition for compactnes and the fact that Hausdorff and packing dimensions of $\mathbf{G}$ and $\psi$-L\'evy trees are the same, as shown in Proposition \ref{fractCIRG}. 
Let us mention that, as shown in the companion paper \cite{BDW2}, the $(\alpha, \beta, \kappa, \mathbf{c})$-continuum multiplicative graphs introduced in this paper are the scaling limits of $\bw$-multiplicative subgraphs.

\medskip

\noi
\textbf{Background and related work.} Let us make briefly the connection with previous articles. 

\smallskip

\noindent 
$\bullet$ \textit{Entrance boundary of the multiplicative coalescent.} The work \cite{AlLi98} of Aldous and Limic already revealed a deep connection between $Y$ as in (\ref{eq:Y_intro}) and the multiplicative coalescent processes. This work also suggests that the family of continuous multiplicative graphs we construct indeed contains all \mm{the} possible limits of random graphs related to the multiplicative coalescent. To be precise, a stochastic process  $(W^{\kappa_{\textrm{AL}}, -\tau_{\textrm{AL}}, \mathbf{c}_{\textrm{AL}}  }_s)_{s\in [0, \infty)}$ is considered in \cite{AlLi98}. Its law is characterised by 
three parameters: $\kappa_{\textrm{AL}}\ino [0, \infty)$, $\tau_{\textrm{AL}} \ino \bbR$ and $ \mathbf{c}_{\textrm{AL}}\ino \elldo_3$. This is actually a rescaled version of the $(\alpha, \beta , \kappa, \mathbf{c})$-process $Y$ defined in (\ref{eq:Y_intro}), simply because 
\begin{equation}
\label{bungabunga}
\forall s\ino [0, \infty), \quad W^{\kappa_{\textrm{AL}}, -\tau_{\textrm{AL}}, \mathbf{c}_{\textrm{AL}}  }_s  = Y_{s/ \kappa}, \; \textrm{where} \quad \kappa_{\textrm{AL}} = \frac{\beta}{\kappa}, \quad \tau_{\textrm{AL}}= \frac{\alpha}{\kappa} \quad \textrm{and} \quad \mathbf{c}_{\textrm{AL}}=\mathbf{c}. 
\end{equation}
Then Aldous and Limic show (Theorem 2 \& Theorem 3 in \cite{AlLi98}) that there is a one-to-one correspondence between the laws of the excursion lengths of $W^{\kappa_{\textrm{AL}}, -\tau_{\textrm{AL}}, \mathbf{c}_{\textrm{AL}}}$ and the marginals of the extremal eternal versions of the multiplicative coalescent. This connection origins from a representation of the (finite-state) coalescent process in terms of the $\cG_{\bw}$ graph, first observed in Aldous \cite{Al97}. The work \cite{Al97} and \cite{AlLi98} is the primary motivation for us to consider these particular classes of random graphs and part of their results (in particular the convergence of the excursion lengths) have figured in our proof of the limit theorems in \cite{BDW2}. See Section 2.3.4 in \cite{BDW2} for more details.

\smallskip

\noindent 
$\bullet$ \textit{Graph limits as instances of continuous multiplicative graphs.} 
As important examples of the family of $(\alpha, \beta , \kappa, \mathbf{c})$-multiplicative graphs, let us
mention that 
\begin{itemize}
\item[--]
{\it in the scaling limit of the Erd\H{o}s--R\'enyi graphs} we find a continuous multiplicative graph with $\alpha \in\bbR$, $\beta\! = \! 1$ and $\mathbf{c}\! = \! 0$.
This result is due to Addario-Berry, Goldschmidt \& B.~\cite{AdBrGo12a}. 
Let us also mention that as shown by Bhamidi, Sen \& X.~Wang in \cite{BhSeWa14}, the same limit object appears in the $\bw_{n}$-multiplicative graphs which are in the basin of attraction of the Erd\H{o}s--R\'enyi graph. 

\item[--]
{\it in the scaling limit of the multiplicative graphs with power-law weights}, we find continuous multiplicative graphs parametrised by $\alpha \ino \bbR$, $\kappa \ino (0, \infty)$, $\beta\! = \! \beta_0\! = \! 0$ and $c_j = q \,  j^{-\frac{1}{\rho}}$, for all $j\! \geq \! 1$. This result is due to Bhamidi, van der Hofstad \& Sen \cite{BhHoSe15}. 
Let us also mention that Conjecture 1.3 right after Theorem 1.2 in \cite{BhHoSe15} is solved by our 
Proposition \ref{fractCIRG} that asserts the following: if $\alpha\ino [0, \infty)$, $\kappa \ino (0, \infty)$, $\beta\! = \! 0$ and 
$c_j\! = \! q \,  j^{-1/\rho}$, then $\eta\! =\! \gamma\! = \!  \rho \!   -\! 1$ (which corresponds to $\tau\! -\! 2$ in \cite{BhHoSe15}) and 
$$ \textrm{$\bP$-a.s.~for all $k\! \geq \! 1$,} \quad \mathrm{dim}_H (\mathbf{G}_k)\! =\! \mathrm{dim}_p (\mathbf{G}_k) =  \frac{\rho \! -\! 1}{\rho \! -\! 2}, $$
where $\mathrm{dim}_H$ and $ \mathrm{dim}_p$ stand respectively for the Hausdorff and for the packing dimensions.

\end{itemize}


 






\medskip
 
\noi
\textbf{Plan of the paper.}  The paper is organized as follows. Section~\ref{sec:results} is dedicated to a detailed exposition of the setting and results: we start with a motivation using the setting of discrete objects, and move towards the continuum objects. The proofs in the discrete setting, namely of the main respresentation theorem and of the embedding of multiplicative graphs into a Galton--Watson forest are found in Sections~\ref{PfThmloiGG} and \ref{EmbedGW}. In Section~\ref{sec:continuous}, we give the proofs concerning the continuum objects together with the relevant background on Lévy trees: here, this is mainly about proving the continuum analog of the embedding into a forest. Finally, Appendix~\ref{Pinfrac} contains auxiliary results used to deal with the fractal dimensions.

\color{black}


\section{Exposition of the main results}\label{sec:results}

\subsection{Exploration of discrete multiplicative random graphs} \label{sec: dis}
We briefly describe the model of discrete random graphs that are considered in this paper and we discuss a combinatorial construction thanks to a LIFO-queue. 
Unless the contrary is specified, all the random variables that we consider are defined on the same probability space $(\Omega, \ccF, \bP)$. 
The graphs $G= (\cV(G), \ccE (G))$ that we consider are \textit{not oriented}, without either loops or multiple edges: $\ccE(G)$ is therefore a set consisting of 
unordered pairs of distinct vertices.  

Let $n\! \geq \! 2$ and let $\bw \! = \! (w_1, \ldots , w_n)$ be a set of \textit{weights}: namely, it is a set of positive real numbers such that 
$ w_1 \! \geq \! w_2 \! \geq \! \ldots \! \geq \! w_n \! >\! 0$. 
We shall use the following notation. 
\begin{equation}
\label{moment}
\forall r\ino (0, \infty), \quad  \sigma_r (\bw)= w^r_1+ \ldots + w^r_n . 
\end{equation} 
The random graph $\cG_\bw$ is said to be \textit{$\bw$-multiplicative} if $\cV (\cG_\bw) \! =\!  \{ 1, \ldots, n\}$ and if 
\begin{equation}
\label{defGwdir}
\textrm{the r.v.~$(\un_{\{ \{ i,j \}\in \ccE (\cG_\bw) \}})_{1\leq  i <  j  \leq n}$, are independent and} \; \bP \big( \{ i,j\}\ino \ccE (\cG_\bw)\big)=1\! -\! 
e^{-w_iw_j/\sigma_1 (\bw)} . 
\end{equation}

In the entire article, for a stochastic process $X$ we use interchangeably $X_t$ and $X(t)$ depending on which is more convenient or readable.

\begin{figure}[tbp]
\includegraphics[page=1,scale=.85]{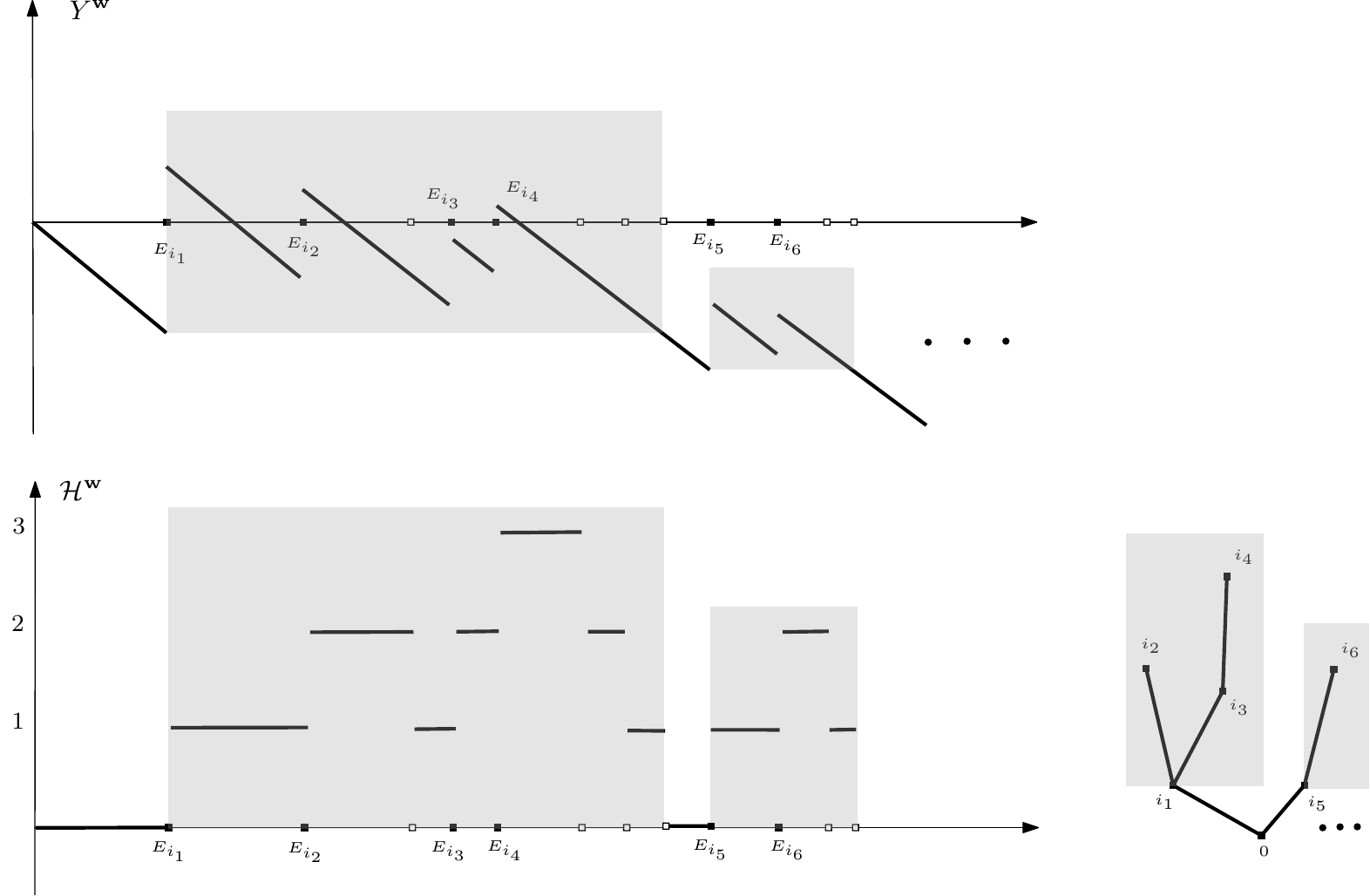}
\caption{{\small An example of $Y^\bw$ and the associated exploration tree. Above, an illustration of $Y^\bw$. The black squares ${\scriptstyle\blacksquare}$ on the  abscissa correspond to  the arrival times of clients, namely, the points of $\{E_i, 1\! \leq i\! \leq \! n\}$. The white squares ${\scriptstyle\square}$ on the  abscissa mark the departures of clients:  By the LIFO rule, the client arriving at time $E_i$ leaves at time $\inf\{t>E_i: Y^\bw_t<Y^\bw_{E_i-}\}$. Below, the exploration tree associated to this queue. Observe that each grey block contains a subtree above the root $0$ and is encoded by an excursion of $Y^\bw-J^\bw$. \cq}  }
\label{fig:LIFO} 
\end{figure}

\subsubsection{A LIFO queueing system exploring multiplicative graphs.}
\label{YLIFOsec}
Let us first explain how to generate a $\bw$-multiplicative graph $\cG_\bw$ thanks to the queueing system that is 
described as follows: there is a \textit{single server}; at most one client is served at a time; the server applies the \textit{Last In First Out} policy (\textit{LIFO}, 
for short). Namely, when a client enters the queue, she/he interrupts the service of the previously served client (if any) and the 
new client is immediately served. When the server completes the service of a client, it comes back to the last arrived client whose service has been interrupted (if there is any).
Exactly $n$ clients will enter the queue; each client is labelled with a distinct integer of $\{ 1, \ldots , n\}$ and $w_i$ stands for the \textit{total amount of time of service that is needed by Client 
$i$} who enters the queue at a time denoted by $E_i$; we refer to $E_i$ as the \textit{time of arrival of Client $i$}; 
we assume that $E_1, \ldots , E_n \ino (0, \infty )$ are distinct. 
For the sake of convenience, we label the server by $0$ and we set $w_0\! = \! \infty$. The single-server 
LIFO queueing system is completely determined by the (always deterministic) 
times of service $\bw$ and the times of arrival $\underline{E}\! = \! (E_1, \ldots, E_n)$, that are random variables whose laws are specified below. We introduce the following processes:
\begin{equation}
\label{Ywdef}
\forall t\ino [0, \infty), \quad Y^{\bw}_t= \! - t + \!\! \sum_{1\leq i\leq n} \!\! w_i \un_{\{ E_i \leq t\}} \quad \textrm{and} \quad \JJ^\bw_t = \inf_{s\in [0, t] } Y^\bw_s \,.
\end{equation}
The \textit{load} at time $t$ (namely the time of service still due by time $t$) is then $Y^\bw_t - \JJ^\bw_t $. We shall sometimes call $Y^{\bw}$ the algebraic load of the queue. 
The LIFO queue policy implies that Client $i$ arriving at time $E_i$ will leave the queue at time $\inf\{t\ge E_i: Y^\bw_t<Y^\bw_{E_i-}\}$, namely the first moment when the service load falls back to the level right before her/his arrival.
We shall refer to the previous queueing system as the \textit{$\bw$-LIFO queueing system}. 

\medskip

\noi
\textbf{The exploration tree.} Denote by $V_t \ino \{ 0, \ldots , n\}$ the label of the 
client who is served at time $t$ if there is one;  
namely, if the server is idle right after time $t$, we set $V_t\! =\! 0$ (see Section 
\ref{PfThmloiGG} for a formal definition). First observe that $V_0\! = \! 0$ and that $t \! \mapsto \! V_t$ is c\`adl\`ag. By convenience, we set $V_{0-}\! = \! 0$. Next note 
that $V_{E_j}\! = \! j$ and that $V_{E_j-}$ is the label of the client who was served when Client $j$ entered the queue. Then, the $\bw$-LIFO queueing system induces an \textit{exploration tree} 
$\cT_{\!\!  \bw}$ with vertex set $\cV(\cT_{\!\! \bw})$ and edge set $\ccE(\cT_{\!\! \bw})$ that are defined as follows: 
\begin{equation}  
 \label{explotree}
\cV( \cT_{\!\!  \bw})=\{ 0, \ldots , n\} \quad \textrm{and} \quad \ccE (\cT_{\!\!  \bw}) = \big\{ \{ V_{E_j-}, j\} ; 1\! \leq \! j \! \leq \! n\big\}.
\end{equation}
The tree $\cT_{\!\!  \bw}$ is rooted at $0$, which allows to view it as a family tree: the ancestor is $0$ (the server) and 
Client $j$ is a child of Client $i$ if Client $j$ enters the queue while Client $i$ is served. 
In particular, the ancestors of Client $i$ are those waiting in queue while $i$ is being served. See Figure \ref{fig:LIFO} for an example.

\medskip

\noi
\textbf{Additional edges.} We obtain a graph $\cG_{\! \bw}$ by 
adding edges to $\cT_{\!\! \bw}$ as follows. Conditionally given $\underline{E}$, let 
\begin{equation}
\label{Poissurpl}
\cP_\bw\! = \!\!\! \! \!\! \sum_{\; 1\leq p\leq \bp_\bw} \!\!\!  \delta_{(t_p, y_p) }
 \; \textrm{be a Poisson pt.~meas.~on $[0, \infty)^2$ 
  with intensity $\frac{{_1}}{{^{\sigma_1 (\bw)}}}  \un_{\{ 0< y <Y^\bw_t -\JJ^\bw_t \}} \, dt\, dy$.}  
\end{equation}
Note that a.s.~the number of atoms $\mathbf{p}_\bw$ is finite since $Y^\bw\! -\! \JJ^\bw$ is null eventually. We set:
\begin{equation}
\label{pin1}
 \Ptt_\bw= \big( (s_p, t_p)\big)_{1\leq p \leq \bp_\bw} \quad \textrm{where} \quad  s_p\! = \! \inf \big\{ s\ino [0, t_p] :  \inf_{u\in [s, t_p]} Y^\bw_u \! - \! \JJ^\bw_{u} > y_p   \big\}, \; 1\! \leq \! p \! \leq \! \bp_\bw \; .
\end{equation}
Note that 
$s_p$ is well-defined since $y_p<Y^\bw_{t_p}-\JJ^\bw_{t_p}$. We then derive $\cG_\bw$ from $\cT_{\!\! \bw}$ and $\Ptt_\bw$ by setting: $ \cV(\cG_\bw) \! =\! 
\{ 1, \ldots, n\} $ and $ \ccE (\cG_\bw) \! = \!  \cA  \sqcup   \cS$, where 
\begin{equation}
\label{defGwE}
 \cA\! = \! \big\{ \{ i,j\}  \ino \ccE(\cT_{\!\! \bw}) : i,j\! \geq \! 1 \big\} \quad \textrm{and} \quad    
 \cS \! = \! \big\{ \{ V_{s_p}, V_{t_p} \}; 1\! \leq \! p \! \leq \! \bp_\bw      \big\} \backslash \cA . 
\end{equation} 
Note that $V_{s_p}$ is necessarily an ancestor of $V_{t_p}$; in other words, $V_{s_p}$ is in the queue at time $t_p$. 
Moreover, we have $V_{s_p}\ne 0$ a.s, since  $Y^\bw_{s_p}-J^\bw_{s_p}\ge y_p>0$.  It follows that the endpoints of an edge belonging to $\cS$ necessarily belong to the same connected component of 
$\cT_{\!\!  \bw}\setminus\{0\}$. Note that $0$ is not a vertex of $\cG_{\! \bw}$; we call $\cS$the set of \textit{surplus edges}. When $\underline{E}$ is suitably distributed, $\cG_{\! \bw}$ is distributed as a 
$\bw$-multiplicative graph: this is the content of the following theorem that is the key to our approach. 
\begin{thm}
\label{loiGG} Keep the previous notation; suppose that $E_1, \ldots, E_n$ are independent exponential r.v.~such that 
$\bE [E_j]\! =\! \sigma_1 (\bw)/w_j$, for all $j\ino \{ 1, \ldots, n\}$.    
Then, $\cG_{\bw}$ is a $\bw$-multiplicative random graph as specified in (\ref{defGwdir}). 
\end{thm}
\noi
\textbf{Proof:} see Section \ref{PfThmloiGG}.  \cqfd 

\smallskip

\noi
\textbf{Height process of the exploration tree.}
For all $t\ino [0, \infty)$, let $\cH^\bw_t$ 
be the number of clients waiting in the line at time $t$. Recall that by the LIFO rule, a client entered at time $s$ is still in the queue at time $t$ if and only if $\inf_{s\le u\le t}Y^\bw_u>Y^\bw_{s-}$.
In terms of $Y^\bw$, $\cH^\bw_t$ is defined by 
\begin{equation}
\label{JHdef}
\cH^\bw_t \! = \# \cJ_t, \;   \textrm{where} \;     \cJ_t\! = \! \big\{ s\ino [0, t] \!  : \!  J^{\bw,s -}_{t} \! <\! J^{\bw,s }_{t} \big\}   \; \textrm{and where} \;  \forall s\ino [0, t], \; J^{\bw,s }_t\! = \! \inf_{r\in [s, t]} \! Y^\bw_r . 
\end{equation}
We refer to $\cH^\bw$ as the \textit{height process associated with $Y^\bw$}. 
Note that $\cH_t^\bw$ is also the height of the vertex $V_t$ in the exploration tree $\cT_{\!\! \bw}$.
Actually, this process is a specific \textit{contour} of the exploration tree $\cT_{\!\! \bw}$ and we easily check that it (a.s.) codes its graph-metric $d_{\cT_{\!\! \bw}}$ 
as follows:
\begin{equation}
\label{gloupii}
\forall s, t\ino [0, \infty), \quad d_{\cT_{\!\! \bw} } (V_s, V_t) = \cH^\bw_t+\cH^\bw_s - \, 2 \!\!\!\!\!\!  \min_{\quad r \in [s  \wedge t , s\vee t] }  \!\!\!\!\! \!    \cH^\bw_r  \; .
\end{equation}
See Figure \ref{fig:LIFO} for more details.

\smallskip

\noi
\textbf{The connected components of the $\bw$-multiplicative graph.}
\label{conGw1}
The above LIFO-queue construction of the $\bw$-multiplicative graph $\cG_\bw$ has the following nice property: the vertex sets of the connected components of $\cG_\bw$ coincide with the vertex sets of the connected components of $\cT_{\!\! 
\bw} \setminus\{0\}$,  since surplus edges from $\cS$ are only added inside the connected components of the latter. 
We equip $\cG_\bw$ with the measure 
$\bm_\bw \! = \! \sum_{1\leq j\leq n} w_j \delta_j $ that is the push-forward measure of the Lebesgue measure via the map 
$V$ restricted to the set of times $\{ t\ino [0, \infty) \! :\!  V_t \! \neq \! 0\}$.  
Denote by $\bq_\bw$ the number of connected components of $\cG_{ \bw}$ that are denoted 
by $\cG^\bw_{\! 1}$, $\ldots$, $\cG_{\! \bq_\bw}^\bw$; here indexation is such that 
$$\bm_\bw (\cV( \cG_{\! 1}^\bw) )\! \geq \! \ldots \!  \geq \! \bm_\bw (\cV( \cG_{\! \bq_\bw}^\bw)  ) \; .$$  
Note that for all $k\ino \{ 1, \ldots , \bq_\bw\}$, $\cG_{\! k}^\bw$ corresponds to a connected component $\cT_{\! \! k}^\bw$ of $\cT_{\!\! \bw}\backslash \{ 0\}$ such that 
$$ \cV( \cG_{\! k}^\bw) = \cV( \cT_{\! \! k}^\bw)  \quad \textrm{and} \quad \cE (\cG_{\! k}^\bw) =  \cE (\cT_{\! \! k}^\bw)  \sqcup \cS_k \quad \textrm{where} \quad \cS_k = 
\big\{ \{ i,j\} \ino \cS : i,j\ino \cV( \cT_{\! \! k}^\bw)  \big\} . $$
 
Let $d_{\cG_{\! k}^\bw}$ and  $d_{\cT_{\!\! k}^\bw}$ be the respective graph-metrics of $\cG_{\! k}^\bw$ and of $\cT_{\!\! k}^\bw$ and denote by $\bm^\bw_k$ the restriction of $\bm_\bw$ to $\cV(\cG_{\! k}^\bw)$. 
Then, $\cH^\bw$ and $\Ptt_\bw$ completely encode the sequence 
$((\cG_{\! k}^\bw  , d_{\cG_{\! k}^\bw} , \bm_k^\bw))_{1\leq k\leq \bq_\bw}$ of connected components  viewed as measured metric spaces. 
Indeed, we will see that each excursion of $\cH^\bw$ above zero corresponds to a connected component $\cT^{ \bw}_{\!\! k}$ of 
$\cT_{\!\! \bw} \backslash \{ 0\}$, the length of the excursion interval is $\bm_\bw (\cV(\cT_{\!\! k}^\bw))$ and 
$\cS_k$ corresponds to pinching times that fall in this excursion interval. A formal description of this requires some preliminary work on the measured metric spaces and is therefore postponed to Section \ref{codgrasec}.

\subsubsection{Embedding the exploration tree into a Galton--Watson tree}
\label{XLIFOsec}

\noi
\textbf{A Markovian LIFO queueing system.} 
We embed the $\bw$-LIFO queueing system governed by $Y^\bw$ into the following Markovian LIFO queueing system: 

\smallskip

\begin{compactenum}

\item[]\textit{A single server receives in total an infinite number of clients; it applies the LIFO policy; 
clients arrive at unit rate; each client has a} \textrm{type} \textit{that is an integer ranging in $\{ 1, \ldots, n\}$; the amount of service required by a client of type $j$ is $w_j$; types are i.i.d.~with law $\nu_\bw \! = \! \frac{1}{\sigma_1 (\bw)}\sum_{1 \leq j\leq n} w_j \delta_{j} $.}
\end{compactenum} 

\medskip

\noi
Let $\tau_k$ be the arrival-time of the $k$-th client and let $\Jtt_k$ be the type of the $k$-th client. Then, the Markovian LIFO queueing system is entirely characterised by $\sum_{k\geq 1} \delta_{(\tau_k , \Jtt_k)}$ that is a Poisson point measure on $[0, \infty) \! \times \! \{ 1, \ldots, n\}$ with intensity 
$\ell \! \otimes \! \nu_\bw$, where $\ell $ stands for the Lebesgue measure on $[0, \infty)$. 
We also introduce the following.  
\begin{equation}
\label{Xwdef}
 \forall t\ino [0, \infty), \quad X^\bw_t  =  -t + \sum_{k\geq 1} w_{\Jtt_k}\un_{[0, t]} (\tau_k) \;  \quad \textrm{and} \quad I^\bw_t \! = \! \inf_{s\in [0, t]} X^\bw_s . 
\end{equation} 
Then, $X^\bw_t \! -\! I^\bw_t$ is the load of the Markovian LIFO-queueing system and $X^{\bw}$ is called the algebraic load of the queue. 
Note that $X^\bw$ is a spectrally positive L\'evy process with initial value $0$ whose law is characterized by 
its Laplace exponent 
$\psi_\bw\! : \! [0, \infty) \! \rightarrow \! \bbR$ given for all $t, \lambda \ino [0, \infty)$ by: 
\begin{multline}
\label{psibwdef}
 \bE \big[ e^{-\lambda X^\bw_t}\big]\! = \! e^{t\psi_\bw (\lambda)} \quad \textrm{where} \\
 \psi_\bw (\lambda)  =  \alpha_\bw \lambda + \!\!\!\!   \sum_{1\leq j\leq n}\!  \! \frac{_{w_j} }{^{\sigma_1 (\bw)}} \big( e^{-\lambda w_j}\! -\! 1\! + \! \lambda w_j \big) \quad \textrm{and} \quad \alpha_\bw \! := \! 1\! -\! \frac{_{\sigma_2 (\bw)}}{^{\sigma_1 (\bw)}} \; .
\end{multline}
Here, recall from (\ref{moment}) that $\sigma_2 (\bw)\! = \! w_1^2 + \ldots + w_n^2$. 
Note that if $\sigma_2 (\bw) / \sigma_1 (\bw)\! \leq \!  1$, then $\alpha_\bw \! \geq \! 0$, a.s.~$\liminf_{t\rightarrow \infty} X^\bw_t \!= \! -\infty$ and the queueing system is recurrent: all clients are  served completely; if $\sigma_2 (\bw) / \sigma_1 (\bw)\! >\!  1$, then $\alpha_\bw \! <\! 0$, a.s.~$\lim_{t\rightarrow \infty} X^\bw_t \!= \! \infty$
and the queueing system is transient: the load tends to $\infty$ and infinitely many 
clients are not served completely. In what follows we shall refer to the following cases: 
\begin{equation}
\label{soususcri}
\textrm{supercritical:} \; \sigma_2 (\bw)\! >\! \sigma_1 (\bw), \quad   \textrm{critical:} \;  \sigma_2 (\bw)\! =\! \sigma_1 (\bw) ,\quad   \textrm{subcritical:} \; \sigma_2 (\bw)\! < \! \sigma_1 (\bw). 
\end{equation}
%
%
%
%
\begin{figure}[tbp]
\centering
\includegraphics[scale=.85]{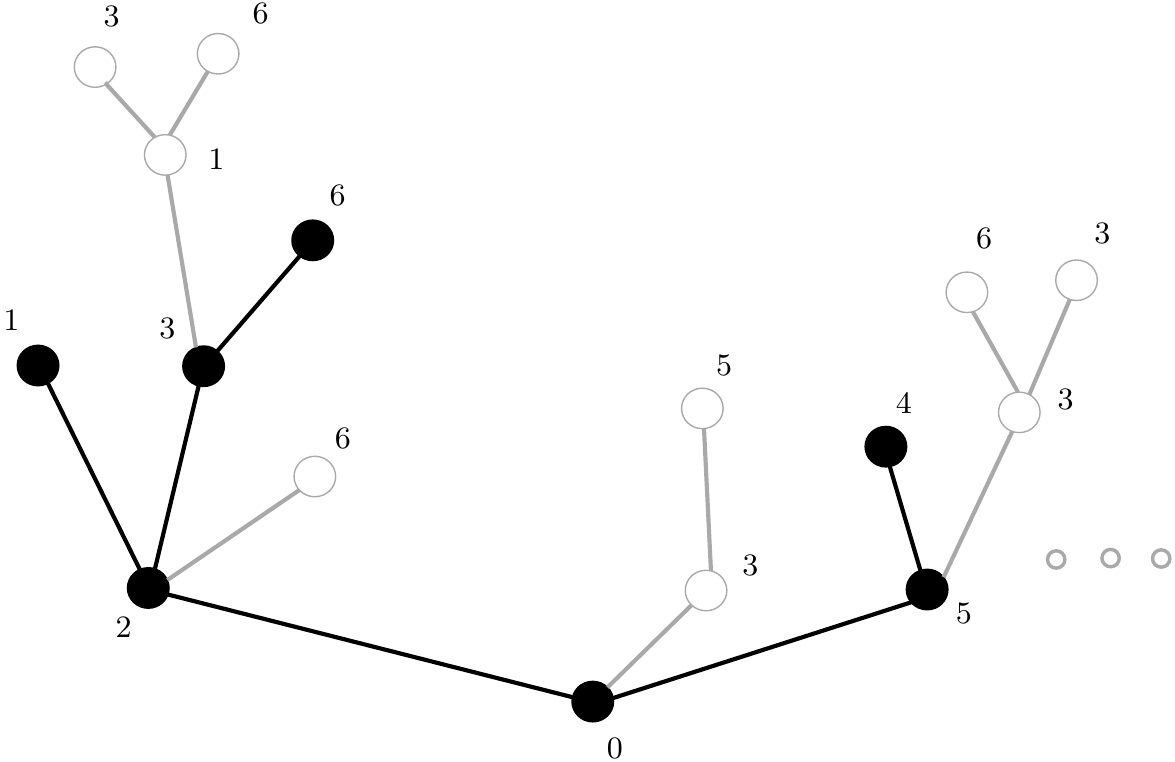}
\caption{{\small Colouring the clients of the Markovian LIFO queue. In this example, we use the exploration tree representation of  the queue. Clients correspond to nodes in the tree; their types are the numbers next to them.  The lexicographic order of the tree (bottom to top, left to right in the picture) corresponds to the order of arrival of the clients.  Applying the colouring rules, we color the clients one by one in this order: blue clients are depicted by $\bullet$, red ones by $\circ$.  Observe that the blue clients form a subtree of the initial tree. Also observe in this example that the first blue client of type $6$ is not the first type-$6$ client in the queue: there is one previous to it, which has been coloured in red because of a red parent. \cq }}
\label{fig:tree_pruning}
\end{figure}

\noi
\textbf{Colouring the clients of the Markovian queueing system.} In critical or subcritical cases, we recover the $\bw$-LIFO queueing system governed by $Y^\bw$ from the Markovian one by colouring each client in the following recursive way.   

\smallskip

\begin{compactenum}
\item[] \texttt{Colouring rules.} \textit{Clients are coloured in red or blue. If the type $\Jtt_k$ of the $k$-th client already 
appeared among the types of the blue clients who previously entered the queue, then the $k$-th client is red.  
Otherwise the $k$-th client inherits her/his colour from the colour of the client who is 
currently served when she/he arrives (and this colour is blue if there is no client served when she/he arrives: namely, we consider that the server is blue).}
\end{compactenum}

\smallskip

\noi
Note that the colour of a client depends in an intricate way on the types of the clients who entered the queue previously. 
For instance, a client who is the first arriving of her/his type is not necessarily coloured in blue; see Figure \ref{fig:tree_pruning} for an example. In critical or subcritical cases, one can check that exactly $n$ clients are coloured in blue and their types are necessarily distinct. While 
a blue client is served, note that her/his ancestors (namely, the other clients waiting in the line, if any) are blue too. Actually, we will see that the sub-queue constituted by the 
blue clients corresponds to the previous $\bw$-LIFO queue in critical and subcritical cases. 

\begin{figure}[tb]
\centering
\includegraphics[scale=.7]{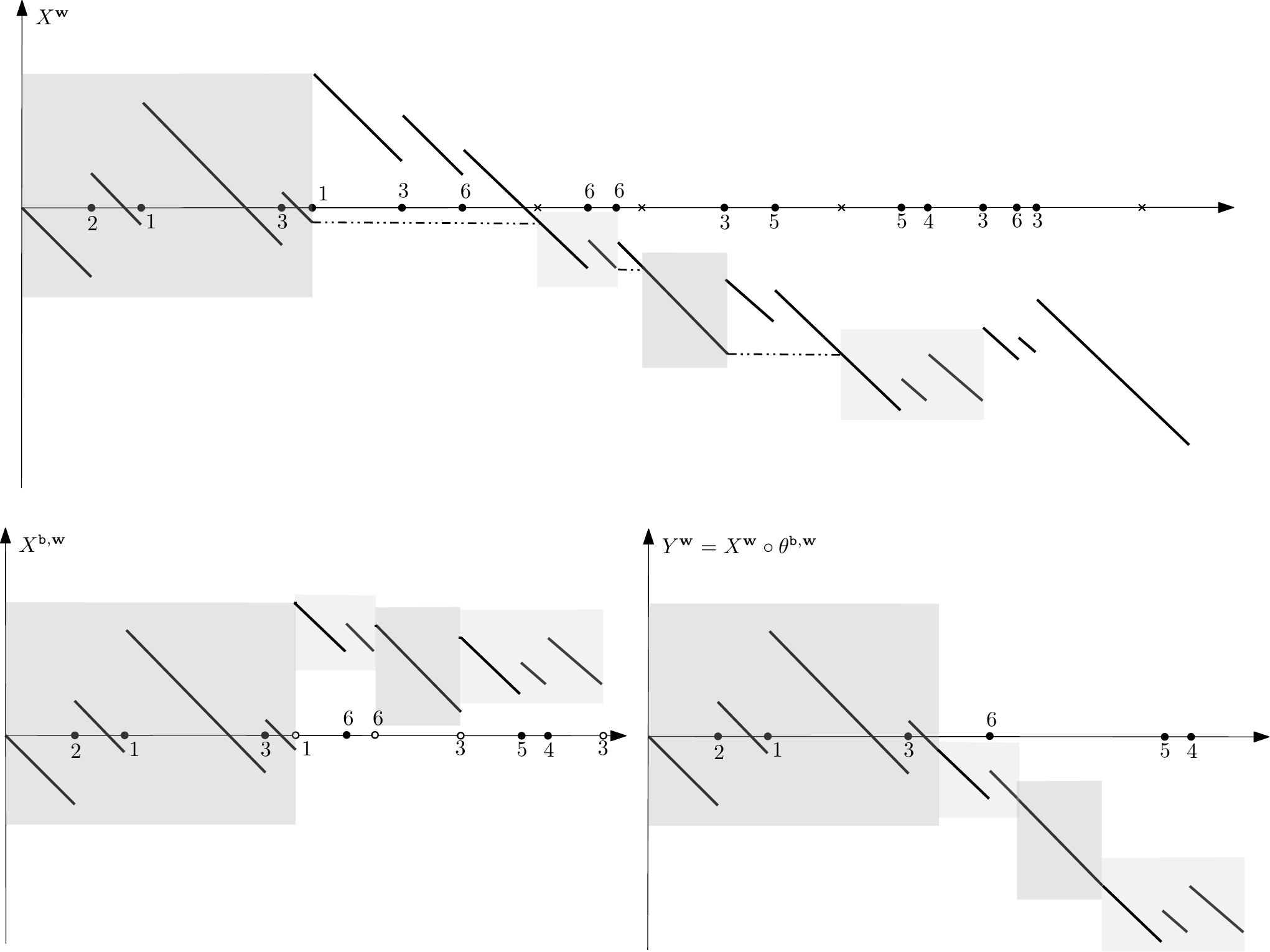}
\caption{{\small Decomposition of $X^\bw$ into $X^{\mathtt b, \mathbf w}$ and $X^{\mathtt r, \mathbf w}$. We take the same example as in Figure \ref{fig:tree_pruning}. Above, the process $X^\bw$: clients are in bijection with its jumps; their types are the numbers next to the jumps.  The grey blocks correspond to the set $\mathtt{Blue}$. 
Concatenating these blocks yields the blue process $X^{\mathtt b, \mathbf w}$. The remaining pieces give rise to the red process $X^{\mathtt r, \mathbf w}$. Concatenating the grey blocks but \textbf{without} the final jump of each block yields $Y^\bw$. Alternatively, we can obtain $Y^\bw$ by removing the temporal gaps between the grey blocks in $X^{\bw}$: this is the graphic representation of  $Y^\bw=X^\bw\circ\theta^{\mathtt b, \mathbf w}$.
Observe also that each connected component of $\mathtt{Red}$ begins with the arrival of a client whose type is a repeat among the types of the previous blue ones,  and ends with the departure of this red client, marked by  ${\scriptstyle\times}$ on the abscissa. \cq}}
\label{fig:proc_color} 
\end{figure}

In supercritical cases, however, we could end up with (strictly) less than $n$ blue clients so the blue sub-queue is only a part of the $\bw$-LIFO queue governed by $Y^\bw$. To deal with this problem and to get a definition of the blue/red queue in a way that can be extended to the continuous setting, we proceed as follows. We first introduce the two following independent random point measures on $[0, \infty) \! \times \! \{ 1, \ldots , n\} $: 
\begin{equation}
\label{ccXbrwdef}\ccX_\bw^{\mathtt{b}}\! = \! \sum_{k\geq 1} \delta_{(\tau^{\mathtt{b}}_k , \Jtt^{\mathtt{b}}_k)} \quad \textrm{and} \quad  \ccX_\bw^{\mathtt{r}}\! = \! \sum_{k\geq 1} \delta_{(\tau^{\mathtt{r}}_k , \Jtt^{\mathtt{r}}_k)}\,,
\end{equation}
that are Poisson point measures with intensity $\ell \! \otimes \! \nu_\bw$, where $\ell $ stands for the Lebesgue measure on $[0, \infty)$ and where $\nu_\bw \! = \! \frac{1}{\sigma_1 (\bw)}\sum_{1 \leq j\leq n} w_j \delta_{j} $. From $ \ccX_\bw^{\mathtt{b}}$, we extract the $\bw$-LIFO queue (without repetition) that generates the desired graph and we explain below how to 
mix $ \ccX_\bw^{\mathtt{b}}$ and  $\ccX_\bw^{\mathtt{r}}$ in order to get the coloured Markovian queue seen at the beginning of the section. To that end, we first set:   
 \begin{equation}
\label{Xbrwdef}
X^{\mathtt{b}, \bw}_t  =  -t + \sum_{k\geq 1} w_{\Jtt^{\mathtt{b}}_k}\un_{[0, t]} (\tau^{\mathtt{b}}_k), \quad   X^{\mathtt{r}, \bw}_t  =  -t + \sum_{k\geq 1} w_{\Jtt^{\mathtt{r}}_k}\un_{[0, t]} (\tau^{\mathtt{r}}_k)  \quad \textrm{and} \quad I^{\mathtt{r}, \bw}_t = \inf_{s\in [0, t]} X^{\mathtt{r}, \bw}_s\; .
\end{equation}
Consequently, $X^{\mathtt{b}, \bw}$ and $X^{\mathtt{r}, \bw}$ are two independent spectrally positive L\'evy processes, both with Laplace exponent $\psi_\bw$ given by (\ref{psibwdef}). 
For all $j\ino \{ 1, \ldots , n\}$ and all $t\in[0, \infty)$, we next set:  
\begin{equation}
\label{Njbdef}
 N^{\bw}_j (t)\! = \! \ccX_\bw^{\mathtt{b}} \big( [0, t] \! \times \! \{ j\} \big) \quad \textrm{and} \quad E^\bw_j = \inf \big\{ t\ino [0, \infty) : \ccX^{\mathtt{b}}_\bw ([0, t] \! \times \! \{ j\})\! = \! 1 \big\}. 
 \end{equation} 
Thus, the $N^\bw_j$ are independent homogeneous Poisson processes with jump-rate $w_j/ \sigma_1 (\bw)$ and the r.v.~$(\frac{w_j}{\sigma_1 (\bw)} E^\bw_j)_{1\leq j\leq n}$ are i.i.d.~exponentially distributed r.v.~with unit mean. Note that $X^{\mathtt{b}, \bw}_t \! =\! -t + \sum_{1\leq j\leq n} w_j N^{\bw}_j (t)$. Thanks to the "blue" r.v.~contained in $\ccX^\mathtt{b}_\bw$, we next define the processes coding the $\bw$-LIFO queue (without repetition) that generates the exploration tree of the graph:  
\begin{equation}
\label{YwAwSigw}
 Y^\bw_t \! = \! -t  \, +\!\!\!   \sum_{1\leq j \leq n} \!\!\!  w_j \un_{\{ E^\bw_j \leq t \}}
 \quad \textrm{and} \quad A^\bw_t =   X^{\mathtt{b}, \bw}_t -Y^\bw_t  =  
  \! \sum_{1\leq j\leq n} \!\! w_j (N^\bw_j (t)\! -\! 1)_+ \; .
\end{equation}
Thanks to $(Y^\bw, A^\bw)$ and $X^{\mathtt{r}, \bw}$,  we reconstruct the Markovian LIFO queue as follows: we first define the "blue" time-change that is the 
increasing c\`adl\`ag process $ \theta^{\mathtt{b}, \bw}$ defined for all $t\ino [0, \infty)$ by: 
\begin{equation}
\label{thetabw}
 \theta^{\mathtt{b}, \bw}_t \!\! = t + \subo^{\mathtt{r}, \bw}_{A^\bw_t}, \;   \textrm{where for all $x\ino [0, \infty)$, we have set:} \;\;  \subo^{\mathtt{r}, \bw}_x \!  = \! \inf \big\{ t\ino [0, \infty ) \! :  X^{\mathtt{r} , \bw}_t \!\!\!  < \! - x \big\}, 
\end{equation}
with the convention that $\inf \emptyset \! = \! \infty$. Note that $\subo^{\mathtt{r}, \bw}_x\! < \! \infty$ iff $x\! < \! \! -I^{\mathtt{r}, \bw}_\infty\! = \! \lim_{t\rightarrow \infty} \! - I^{\mathtt{r}, \bw}_t$ that is a.s.~finite in supercritical cases (and $ -I^{\mathtt{r}, \bw}_\infty$ is a.s.~infinite in critical and subcritical cases).    
Standard results on spectrally positive L\'evy processes (see e.g.~Bertoin's book \cite{Be96} Ch.~VII.) assert that 
$(\subo^{\mathtt{r}, \bw}_x)_{x\in [0, \infty)}$ is a subordinator (that is defective in supercritical cases) whose Laplace exponent is given for all 
$\lambda \ino [0, \infty)$ by: 
\begin{equation}
\label{gamwexpo}
\bE \big[ e^{-\lambda \subo^{\mathtt{r}, \bw}_x } \big]= e^{-x\psi^{-1}_\bw (\lambda)} \quad  \textrm{where} \quad  \psi^{-1}_\bw (\lambda)\! = \! \inf \big\{ u \ino [0, \infty ) : \psi_\bw (u) \! >\! \lambda \big\}.  
\end{equation}
We set $\varrho_{\bw}\! = \!  \psi^{-1}_\bw(0)$ that is the largest root of $\psi_\bw$. Note that 
$\psi_\bw $ has at most two roots since it is convex: in subcritical or critical cases, $\varrho_{\bw}\! = \! 0$ is the only root of $\psi_\bw$ and in supercritical cases, the roots of $\psi_\bw$ are $0$ and $\varrho_{\bw}\! >\! 0$. 
Note that $\psi^{-1}_\bw$ is continuous and strictly increasing and that 
it maps $[0, \infty)$ onto $[\varrho_\bw , \infty)$. As a consequence of (\ref{gamwexpo}), 
in supercritical cases $-I^{\mathtt{r}, \bw}_\infty$ is exponentially distributed with parameter $\varrho_\bw$. We then set: 
\begin{equation}
\label{T*wdef} 
T^*_\bw\! = \! \sup \{ t\ino [0, \infty)\! :  \theta^{\mathtt{b}, \bw}_t \! < \infty \}=   \sup \{ t\ino [0, \infty)\! : A^\bw_t \! < \! - I^{\mathtt{r}, \bw}_\infty \} \; . 
\end{equation}
In critical and subcritical cases, $T^*_\bw\! = \! \infty$ and $ \theta^{\mathtt{b}, \bw}$ only takes 
finite values. In supercritical cases, a.s.~$T^*_\bw \! < \! \infty$ and we check that $ \theta^{\mathtt{b}, \bw} (T^*_\bw-) \! < \! \infty$. 
Imagine that the bi-coloured Markovian LIFO queue has a set of two clocks which never run simultaneously (as in chess clocks): one clock for the blue queue and one for the red. Then  $ \theta^{\mathtt{b}, \bw}_t$ is the (global) time that has been spent when the clock for the blue queue shows $t$. If we denote by $\mathtt{Blue}$ (resp.~$\mathtt{Red}$) the set of times when blue (resp.~red) clients are served, then we can derive these service times from $ \theta^{\mathtt{b}, \bw}$ as follows (recalling that the server is considered as a blue client) 
\begin{equation} 
\label{blureddef}
\mathtt{Red}\! = \! \bigcup_{t\in [0, \infty)} \big[  \theta^{\mathtt{b}, \bw}_{t-}, 
\theta^{\mathtt{b}, \bw}_t \big) \quad \textrm{and} \quad \mathtt{Blue} \! = \! [0, \infty) \backslash \mathtt{Red} . 
\end{equation}
Note that $\mathtt{Red}$ is a countably infinite union of intervals in critical and subcritical cases and that it is a finite union  
in supercritical cases since
$\big[  \theta^{\mathtt{b}, \bw} (T^*_\bw-),  \theta^{\mathtt{b}, \bw} (T^*_\bw))\! = \! [\theta^{\mathtt{b}, \bw} (T^*_\bw-), \infty)$. We next introduce the inverse time-changes $\Lambda^{\mathtt{b}, \bw}$ and $\Lambda^{\mathtt{r}, \bw}$ as follows: 
\begin{equation}
\label{Lambdefi}
\Lambda^{\mathtt{b}, \bw}_t \!  = \! \int_0^t \!\! \un_{\mathtt{Blue}} (s) \, ds =  \inf \big\{ s\ino [0, \infty) \! :  \theta^{\mathtt{b}, \bw}_s \! >\! t \}  \quad \textrm{and} \quad 
 \Lambda^{\mathtt{r}, \bw}_t \!\! = t \! -\!  \Lambda^{\mathtt{b}, \bw}_t \! \! = \! \! \int_0^t \!\! \un_{\mathtt{Red}} (s) \, ds. 
\end{equation}
The processes $\Lambda^{\mathtt{b}, \bw}$ and $\Lambda^{\mathtt{r}, \bw}$ are continuous, nondecreasing and a.s.~$\lim_{t\rightarrow \infty}  \Lambda^{\mathtt{r}, \bw}_t \! = \! \infty$. 
In critical and subcritical cases, we also get a.s.~$\lim_{t\rightarrow \infty}  \Lambda^{\mathtt{b}, \bw}_t \! = \! \infty$ and $\Lambda^{\mathtt{b}, \bw} (\theta^{\mathtt{b}, \bw}_t)\! = \! t$ for all $t\ino [0, \infty)$.  However, in supercritical cases, $ \Lambda^{\mathtt{b}, \bw}_t \! = \! T^*_\bw$ for all $t\ino [\theta^{\mathtt{b}, \bw} (T^*_\bw-), \infty)$ and a.s.~for all $t\ino [0, T^*_\bw)$, 
$\Lambda^{\mathtt{b}, \bw} (\theta^{\mathtt{b}, \bw}_t)\! = \! t$. 
We next derive the load of the Markovian queue $X^\bw$ from $\ccX^{\mathtt{b}}_\bw$ and $\ccX^{\mathtt{r}}_\bw$ as follows.   
\begin{prop}
\label{Xwfrombr} Let $\ccX^{\mathtt{b}}_\bw$ and $\ccX^{\mathtt{r}}_\bw$ be as in (\ref{ccXbrwdef}). Let $X^{\mathtt{b}, \bw}$ and $X^{\mathtt{r}, \bw}$ be defined by (\ref{Xbrwdef}) and let $\Lambda^{\mathtt{b}, \bw}$ and $\Lambda^{\mathtt{r}, \bw}$  be given by (\ref{Lambdefi}). We define the process $X^\bw$ by:   
\begin{equation}
\label{redblumix}
\forall t\ino [0, \infty) , \qquad X^{\bw}_t =  X^{\mathtt{b}, \bw}_{ \Lambda^{\mathtt{b}, \bw}_t } + X^{\mathtt{r}, \bw}_{ \Lambda^{\mathtt{r}, \bw}_t } \; .
\end{equation}
Then, $X^{\bw}$ has the same law as $X^{\mathtt{b}, \bw}$ and $X^{\mathtt{r}, \bw}$: namely, it is a spectrally positive L\'evy process with Laplace exponent $\psi_\bw$ as defined in (\ref{psibwdef}). Furthermore, recall from (\ref{YwAwSigw}) the definition of $Y^\bw$ and 
recall from (\ref{T*wdef}) the definition of the time $T^*_\bw$; then, we also get:  
\begin{equation}
\label{YXtheta}
\textrm{a.s.}\; \forall\,t\ino [0, T^*_\bw) , \quad Y^\bw_t \! = \! X^\bw_{\theta^{\mathtt{b}, \bw}_t } . 
\end{equation}
\end{prop}
\noi
\textbf{Proof.} See Section \ref{EmbedGW} for a proof and 
see Figure \ref{fig:proc_color} for an explanation. \cqfd 
 
\medskip 

\medskip 
 
\noi
\textbf{Tree embeddings.} Recall from (\ref{JHdef}) the definition of the height process $\cH^\bw$: 
\begin{equation}
\label{JHdef2}
\cH^\bw_t \! = \# \cJ_t, \;   \textrm{where} \;     \cJ_t\! = \! \big\{ s\ino [0, t] \!  : \!  J^{\bw,s -}_{t} \! <\! J^{\bw,s }_{t} \big\}   \; \textrm{and where} \;  \forall s\ino [0, t], \; J^{\bw,s }_t\! = \! \inf_{r\in [s, t]} \! Y^\bw_r . 
\end{equation}
Recall that $\cH^\bw_t$ is the number of clients waiting in the $\bw$-LIFO queue (without repetition) governed by $Y^\bw$ and that it is the contour process of the exploration tree $\cT_{\!\! \bw}$. Similarly, 
we denote by $H^\bw_t$ the number of clients waiting in the Markovian 
LIFO queue governed by the process $X^\bw$ given in Proposition \ref{Xwfrombr}: $H^\bw_t$ is defined as follows. 
\begin{equation}
\label{XJHdef}
H^\bw _t \! = \# \cK_t, \;   \textrm{where} \;    \cK_t \! =\!  \big\{ s\ino [0, t]\!  : \! I^{\bw, s-}_{t} \! <\! I^{\bw, s}_{t} \big\}   \;    \textrm{and where} \;     \forall s\ino [0, t], \;  I^{\bw, s}_t= \inf_{r\in [s, t]}X^\bw_r . 
\end{equation}
We shall refer to $H^\bw$ as the \textit{height process associated with $X^{\bw}$.} The process $H^\bw$ is the contour process of a tree that is given as follows: \textit{vertices are the clients and the server is viewed as the root; the $k$-th client to enter the queue is a child of the $l$-th one if the $k$-th client enters when the $l$-th client is served}. It is easy to check that in critical and subcritical cases, this tree is made of a forest of   
i.i.d.~Galton--Watson trees whose roots are all joined to a common vertex representing the server, and 
the common offspring distribution is given by: 
\begin{equation}
\label{mupoissw}
 \forall k \ino \bbN, \quad \mu_\bw (k) \! =\!  \sum_{1\leq j\leq n} \frac{{w_j^{k+1} e^{-w_j} }}{{ \sigma_1 (\bw) \,  k!}}  \; .
\end{equation} 
Observe that $\sum_{k\geq 0} k\mu_\bw (k) \! = \! \sum_{1\leq j\leq n} w_j^2/ \sigma_1 (\bw)\! = \! \sigma_2 (\bw) / \sigma_1 (\bw) $.  
Thus, in critical and subcritical cases, the Galton--Watson trees are a.s.~finite and $H^\bw$ fully explores the whole tree. 
In supercritical cases, we can still think of a sequence of i.i.d.~Galton--Watson trees whose offspring distribution $\mu_\bw$ is given by (\ref{mupoissw}), but we only see in $H^{\bw}$ a subsequence of these trees up to (part of) the first infinite member. 
Namely, the exploration of $H^\bw$ does not go beyond the first infinite line of descent. 
The embedding of the tree $\cT_{\!\! \bw}$ coded by $\cH^\bw$ into the Galton--Watson 
tree coded by $H^\bw$ is given by the following lemma. 
\begin{lem}
\label{Hthetalem} Let $\ccX^{\mathtt{b}}_\bw$ and $\ccX^{\mathtt{r}}_\bw$ be as in (\ref{ccXbrwdef}).
Let $\theta^{\mathtt{b}, \bw}$ be defined by (\ref{thetabw}) and let $T^*_\bw$ be given by (\ref{T*wdef})
Let $\cH^\bw$ and $H^\bw$ be defined resp.~by (\ref{JHdef2}) and (\ref{XJHdef}). 
Then, 
\begin{equation}
\label{HYenHX}
\textrm{a.s.} \; \forall t\ino [0, T^*_\bw), \quad \cH^\bw_t = H^\bw_{\theta^{\mathtt{b} , \bw}_t}  \; . 
\end{equation}
%
%
\end{lem}  
\noi
\textbf{Proof.} See Section \ref{EmbedGW}. \cqfd 

\medskip

\noi
Although the law of $\cT_{\!\! \bw}$ is complicated, (\ref{HYenHX}) allows to define its height process in a tractable way which can then be passed to the limit. Note that the embedding is not complete in supercritical cases; however, it is sufficient to characterise the law of $\cT_{\!\! \bw}$ in terms of $X^\bw$. 
\begin{rem}
\label{HXbpasHY}
Note that the height process of $X^{\mathtt{b} , \bw}$ is actually distinct from $\cH^\bw $. 
Although related to $\cH^\bw$, the tree coded by $X^{\mathtt{b} , \bw}$ is not relevant to our purpose.  \cq    
\end{rem}

\subsection{The multiplicative graph in the continuous setting} \label{sec: con}
\subsubsection{The continuous exploration tree and its height process} \label{hautcontset}

\noi
\textbf{Notations and conventions.}
Recall that $\bbN$ stands for the set of nonnegative integers and that 
$\bbN^*\! = \! \bbN \backslash \{ 0\}$. 
We denote by 
$\elldo_{\infty}\! =\!  \big\{ (w_j)_{j\geq 1} \ino [0, \infty)^{\bbN^*} \! \! \!\! : \,  w_j \! \geq \! w_{j+1} \big\}$ the set of \textit{weights}. By an obvious extension of Notation (\ref{moment}), for all $r\ino (0, \infty)$ and 
all $\bw\! = \!  (w_j)_{j\geq 1}\ino \elldo_\infty$, we set $\sigma_r (\bw)\! = \!  \sum_{j\geq 1} w_j^r \ino [0, \infty]$. 
We also introduce the following: 
$$ \elldo_{{r}} = \big\{ \bw \ino \elldo_\infty : \sigma_r (\bw) \! < \! \infty \big\},  
 \quad \textrm{and} \quad  \elldo_{{\! f}}= 
\big\{ \bw \ino \elldo_\infty  :  \exists j_0 \! \geq \! 1 : w_{j_0} \! = \! 0  \big\} .$$

Let  $(\ccF_t)_{t\in [0, \infty)}$ be a filtration on $(\Omega, \ccF)$ that is specified further. 
A process $(Z_t)_{ t\in [0, \infty)}$ is said to be a
$(\ccF_t)$-L\'evy process with initial value $0$ if 
a.s.~$Z$ is c\`adl\`ag, $Z_0 \! = \! 0$ and if for all a.s.~finite $(\ccF_t)$-stopping time $T$, the process 
$Z_{T+ \, \cdot}\! -\! Z_{T}$ is independent of 
$\ccF_{T}$ and has the same law as $Z$. 

Let $(M_j(\cdot))_{j\geq 1}$ 
be a sequence of c\`adl\`ag  $(\ccF_t)$-martingales that are in $L^2$-summable and orthogonal: namely, 
for all $t\ino [0, \infty)$, $\sum_{j\geq 1} \bE \big[ M_j(t)^2\big] \! < \! \infty$ and $\bE [M_j (t)M_k(t)]\! = \! 0$ if $k\! >\! j$. Then $\sum_{j\geq 1}^{_{\perp}} M_j $ stands for the (unique up to indistinguishability) c\`adl\`ag $(\ccF_t)$-martingale $M (\cdot)$ such that  for all $j\! \geq \! 1$ and all $t\ino [0, \infty)$, 
$\bE \big[ \sup_{s\in [0, t]} \big| M (s)  \! -\! \sum_{1\leq k\leq j} M_k (s) \big|^2 \big] \! \leq \! 4 \sum_{l>j} \bE [M_l (t)^2] $, by Doob's inequality.  

\medskip

\noi
\textbf{Blue processes.}
\label{blupro} We fix the following parameters. 
\begin{equation}
\label{parconing}
\alpha \ino \bbR , \quad \beta \ino [0, \infty) , \quad \kappa \ino (0, \infty) , \quad \mathbf{c}\! = \! (c_j)_{j\geq 1} \ino \elldo_3 \; .
\end{equation}
These quantites are the \textit{parameters of the continuous multiplicative graph}: $\mathbf{c}$ plays the same role as $\bw$ in the discrete setting, $\alpha$ is a drift coefficient similar to $\alpha_\bw$, $\beta$ is a Brownian coefficient 
and the interpretation of $\kappa$ is explained later.

  Next, let $(B_t)_{t \in [0, \infty)}$, $(N_j (t))_{t\in [0, \infty)}$, $j\! \geq \! 1$ be processes that satisfy the following. 
\begin{compactenum}

\smallskip

\item[]  $(b_1)$ \textit{$B$ is a $(\ccF_t)$-real valued standard Brownian motion with initial value $0$}.

\smallskip

\item[] $(b_2)$ \textit{For all $j\! \geq \! 1$, $N_j $ is a $(\ccF_t)$-homogeneous Poisson process with jump-rate $\kappa c_j$}. 

\smallskip

\item[] $(b_3)$ \textit{The processes $B$, $N_j$, $j\! \geq \! 1$ are independent}.

\smallskip 
\end{compactenum}

\noi
The \textit{blue} L\'evy process is then defined by 
\begin{equation}
\label{Xblue}
\forall t\ino [0, \infty) , \quad X^{\mathtt{b}}_t = -\alpha t + \sqrt{\beta} B_t + \sum_{j\geq 1}  \!\!\! \,^{\perp}\,  c_j \big( N_j (t) \! -\! c_j \kappa t \big) . 
\end{equation}
Clearly $X^{\mathtt{b}}$ is a $(\ccF_t)$-spectrally positive L\'evy process with initial value $0$ whose law is 
characterized by the Laplace exponent $\psi \! : \! [0, \infty)  \! \rightarrow \! \bbR$ given for all $t, \lambda \ino [0, \infty)$ by: 
\begin{equation}
\label{psidefi}
 \bE \big[ e^{ - \lambda X^{\mathtt{b}}_t } \big] \!  = \!  e^{ t\psi (\lambda) }, \;  \textrm{where} \quad \psi(\lambda)  \! = \!   
 \alpha \lambda +\frac{_{_1}}{^{^2}} \beta \lambda^2 +  \!   \sum_{j\geq 1}   \kappa c_j   \big( e^{-\lambda c_j}\! -\! 1\! + \! \lambda c_j \big). 
 \end{equation}
If $\alpha \! \geq 0$, then a.s.~$\liminf_{t\rightarrow \infty} X^\mathtt{b}_t \! = \! -\infty$ and if $\alpha \! < \! 0$, then a.s.~$\lim_{t\rightarrow \infty} X^\mathtt{b}_t \! = \! \infty$. By analogy with the discrete situation of the previous section, we distinguish the following cases:
\begin{equation}
\label{mmbjmboo}
\textrm{supercritical cases:} \; \alpha \! < \! 0, \quad   \textrm{critical cases:} \;  \alpha \! = \! 0,  \quad   \textrm{subcritical cases:} \; 
 \alpha \! >\! 0. 
\end{equation}
Most of the time, we shall assume that: 
\begin{equation}
\label{varinfinie}
\textrm{either $\beta >0$ or $\sigma_2 (\mathbf{c})\! = \! \infty$.}
\end{equation}
This assumption is equivalent to the fact that $X^{\mathtt{b}}$ has infinite variation sample paths.  

  We next introduce the analogues of $A^\bw$ and $Y^\bw\! $ defined 
in (\ref{YwAwSigw}). To that end, note that $\bE [ c_j (N_j (t) \! -\! 1)_+ ]\! = \!   c_j   \big( e^{- c_j \kappa  t}\! -\! 1\! + \!  c_j \kappa t \big)Â \! \leq \! \frac{_1}{^2} (\kappa t)^2 c^3_j $. Since $\sigma_3 (\mathbf{c})\! < \! \infty$, it makes sense to set the following. 
\begin{equation}
\label{AetYdef}
\forall t\ino [0, \infty) , \quad A_t = \frac{_{_1}}{^{^2}} \kappa \beta t^2 +  \sum_{j\geq 1}  c_j \big( N_j (t) \! -\! 1 \big)_+ \quad \textrm{and} \quad Y_t = X^{\mathtt{b}}_t \! -\! A_t . 
\end{equation}

\vspace{-5mm}

\begin{rem}
\label{Yrepres}
To view $Y$ as in (\ref{Ywdef}), set $E_j \! = \! \inf \{ t\ino [0, \infty)\! :\!  N_j (t)\! = \! 1 \}$, 
note that $ c_j  (N_j (t) \! -\! c_j \kappa t)\! -\!  c_j (N_j(t)\! -\! 1)_+\! = \!  c_j (\un_{\{ E_j \leq t \}} \! -\! c_j \kappa t)$ and 
check that $ c_j (\un_{\{ E_j \leq t \}} \! -\! c_j \kappa t)\! = \! M^\prime_j (t) \! -\!    \kappa c^2_j (t\! -\! E_j)_+$ where $M^\prime_j$ is a centered $(\ccF_t)$-martingale such that 
$\bE [M^\prime_j (t)^2]\! = \! c_j^2 (1\! -\! e^{-c_j \kappa t}) \! \leq \!    \kappa t c_j^3 $ (\textit{indeed,} $M^\prime_j(t)\! = \! M_j (t\wedge E_j)$, where $M_j(t)\! =\! c_j(N_j(t)\! -\! \kappa c_j t )$). Since 
$\bE [ \kappa c^2_j (t\! -\! E_j)_+] \leq \kappa t c_j^2 (1\! -\! e^{-\kappa c_j t})\! \leq \! \kappa^2 t^{2} c_j^3$, it 
makes sense to write for all $t\ino [0, \infty)$: 
\begin{eqnarray}
\label{Ydef}
 Y_t  &=  &  - \alpha t \! -\! \frac{_1}{^2}\kappa  \beta t^2 + \sqrt{\beta} B_t + \sum_{j\geq 1}  \!\!\! \,^{\perp} c_j \big( \un_{\{ E_j \leq t \}} \! -\! \kappa c_j (t\! \wedge E_j)\big)
 - \!\! \sum_{j\geq 1}  \kappa c^2_j (t\! -\! E_j)_+ \nonumber \\
 &  \overset{\textrm{(informal)}}{=} &
  -\alpha t \! -\! \frac{_1}{^2}\kappa \beta t^2 + \sqrt{\beta} B_t +  \sum_{j\geq 1} \! c_j (\un_{\{ E_j \leq t \}} \! -\! c_j \kappa t).
 \end{eqnarray}
Namely the jump-times of $Y$ are the $E_j$ and $\Delta Y_{E_j}\! = \! c_j$. 
Let us quickly mention that the process $Y$ appears in the Aldous--Limic's approach to the exit boundary problem for general multiplicative processes \cite{AlLi98}. 
\cq 
\end{rem}
\begin{lem}
\label{AYnontriv} Let $(\alpha, \beta, \kappa, \mathbf{c})$ be as in (\ref{parconing}). 
Assume that (\ref{varinfinie}) holds. Recall from (\ref{AetYdef}) the definiton of $A$ and $Y$.  
Then, a.s.~the process $A$ is strictly increasing and the process $Y$ has infinite variation sample paths. 
\end{lem}
\noi
\textbf{Proof.} See Section \ref{proAsec}. \cqfd  

\bigskip

\noi
\textbf{Red and bi-coloured processes.}  We next introduce the \textit{red process} $X^{\mathtt{r}}$ that satisfies the following. 
\begin{compactenum}

\smallskip

\item[] $(r_1)$ $X^{\mathtt{r}}$ is a $(\ccF_t)$-spectrally positive L\'evy process starting at $0$ and whose Laplace exponent is $\psi$ as in (\ref{psidefi}). 

\smallskip 

\item[] $(r_2)$ $X^{\mathtt{r}}$ is independent of the processes $B$ and $(N_j)_{j\geq 1}$.   

\smallskip 

\end{compactenum}

\noi
To keep the filtration $(\ccF_t)$ minimal, we may assume that $\ccF_t$ is the completed sigma-field generated by 
$B_s$, $(N_j (s))_{j\geq 1}$ and $X^{\mathtt{r}}_s$, $s\ino [0, t]$. We next introduce the following processes: 
\begin{equation}
\label{thetabdef}
 \forall x, t\ino [0, \infty), \quad \subo^{\mathtt{r}}_x = \inf \{ s \ino [0, \infty):    X^{\mathtt{r}}_s \! < \! -x \} \quad \textrm{and} \quad \theta_t^{\mathtt{b}}= t + \subo^{\mathtt{r}}_{A_t} ,
\end{equation}
 with the convention: $\inf \emptyset \! = \! \infty$. For all $t\ino [0, \infty)$, we set $I^\mathtt{r}_t \! = \! \inf_{s \in [0, t]} X^\mathtt{r}_s $ and 
 $I^{\mathtt{r}}_\infty\! = \! \lim_{t\rightarrow \infty} \! I^{\mathtt{r}}_t $
 that is a.s.~finite in supercritical cases and that is a.s.~infinite in critical or subcritical cases. 
Note that $\subo^{\mathtt{r}}_x\! < \! \infty$ iff $x\! < \! \! -I^{\mathtt{r}}_\infty $.  
Again, standard results on spectrally positive L\'evy processes (see e.g.~Bertoin's book \cite{Be96} Ch.~VII.) assert that 
$(\subo^{\mathtt{r}}_x)_{x\in [0, \infty)}$ is a subordinator (that is defective in supercritical cases) whose Laplace exponent is given for all 
$\lambda \ino [0, \infty)$ by: 
\begin{equation}
\label{gamexpo}
\bE \big[ e^{-\lambda \subo^{\mathtt{r}}_x } \big]= e^{-x\psi^{-1}(\lambda)}  \quad  \textrm{where} \quad  \psi^{-1} (\lambda)\! = \! \inf \big\{ u \ino [0, \infty ) : \psi (u) \! >\! \lambda \big\}.  
\end{equation}
We set $\varrho \! = \!  \psi^{-1}(0)$ that is the largest root of $\psi$. Again, note that $\psi$ 
has at most two roots since it is convex: in subcritical or critical cases, $\varrho\! = \! 0$ is the only root of $\psi$ and in supercritical cases, the roots of $\psi$ are $0$ and $\varrho\! >\! 0$. 
Note that $\psi^{-1}$ is continuous and strictly increasing and that 
it maps $[0, \infty)$ onto $[\varrho , \infty)$. As a consequence of (\ref{gamexpo}), 
in supercritical cases $-I^{\mathtt{r}}_\infty$ is exponentially distributed with parameter $\varrho$. We then set: 
\begin{equation}
\label{T*def} 
T^*\! = \! \sup \{ t\ino [0, \infty)\! :  \theta^{\mathtt{b}}_t \! < \infty \}=   \sup \{ t\ino [0, \infty)\! : A_t \! < \! - I^{\mathtt{r}}_\infty \} \; . 
\end{equation}
In critical and subcritical cases, $T^*\! = \! \infty$ and $ \theta^{\mathtt{b}, \bw}$ only takes 
finite values. In supercritical cases, a.s.~$T^* \! < \! \infty$ and we check that $ \theta^{\mathtt{b}} (T^*-) \! < \! \infty$. We next introduce the following. 
\begin{equation}
\label{2lambda}
\forall t \ino [0, \infty), \quad \Lambda^{\mathtt{b}}_t =  \inf \{ s \ino [0, \infty):    \theta^{\mathtt{b}}_s \! > \! t \} \quad \textrm{and} \quad \Lambda^{\mathtt{r}}_t \! = \! t-    \Lambda^{\mathtt{b}}_t . 
\end{equation} 
The process $\Lambda^{\mathtt{b}}$ is continuous, nondecreasing and 
in critical and subcritical cases, we also get a.s.~$\lim_{t\rightarrow \infty}  \Lambda^{\mathtt{b}}_t \! = \! \infty$ and $\Lambda^{\mathtt{b}} (\theta^{\mathtt{b}}_t)\! = \! t$ for all $t\ino [0, \infty)$. 
However, in supercritical cases, a.s.~$ \Lambda^{\mathtt{b}}_t \! = \! T^*$ for all $t\ino[ \theta^{\mathtt{b}} (T^*-), \infty)$ and a.s.~for all $t\ino [0, T^*)$, 
$\Lambda^{\mathtt{b}} (\theta^{\mathtt{b}}_t )\! = \! t$. We next prove the analogue of Proposition \ref{Xwfrombr}. 
\begin{thm}
\label{Xdefthm} Let $(\alpha, \beta, \kappa, \mathbf{c})$ be as in (\ref{parconing}). 
Assume that (\ref{varinfinie}) holds. Recall from (\ref{Xblue}) the definition of $X^\mathtt{b}$ and recall from (\ref{2lambda}) the definition of $\Lambda^{\mathtt{b}}$ and $\Lambda^{\mathtt{r}}$. 
Then, $\Lambda^{\mathtt{r}}$ is continuous, nondecreasing and a.s.~$\lim_{t\rightarrow \infty}  \Lambda^{\mathtt{r}}_t \! = \! \infty$. Furthermore, if we set 
\begin{equation}
\label{Xdef}
\forall t \ino [0, \infty) , \quad X_t = X^{\mathtt{b}}_{\Lambda^{\mathtt{b}}_t } + X^{\mathtt{r}}_{\Lambda^{\mathtt{r}}_t }\, ,  
\end{equation} 
then, $X$ is a spectrally positive L\'evy process with initial value $0$ and Laplace exponent $\psi$ as in (\ref{psidefi}). 
Namely, $X$, $X^\mathrm{b}$ and  $X^\mathrm{r}$ have the same law. Moreover, 
\begin{equation}
\label{YfromX}
\textrm{a.s.}\; \forall\,t\ino [0, T^*) , \quad Y_t \! = \! X_{\theta^{\mathtt{b}}_t } . 
\end{equation}    
\end{thm}
\noi
\textbf{Proof.} See Section \ref{Xdefthmpf} for the first statement and see Section  \ref{proXYsec} for the second one.   \cqfd 

\medskip

\noi
\textbf{Height process and pinching points.} We next define the analogue of $H^\bw$. To that end, we assume that $\psi$ defined in (\ref{psidefi}) satisfies the following: 
\begin{equation}
\label{contH}
 \int^\infty \frac{d\lambda}{\psi (\lambda)} <\infty . 
\end{equation}
Note that in particular \eqref{contH} entails \eqref{varinfinie}.

Le Gall \& Le Jan \cite{LGLJ98} (see also Le Gall \& D.~\cite{DuLG02}) prove that there exists a \textit{continuous} process $H\! = \! (H_t)_{t\in [0, \infty)}$ such that the following limit holds true for all $t\ino [0, \infty)$ in probability :
\begin{equation}
\label{approHdef}
H_t = \lim_{\varepsilon \rightarrow 0} \frac{1}{\varepsilon} \!  \int_0^{t} \!  \un_{\{  X_s - \inf_{r\in [s, t]} X_r  \leq \varepsilon\}} \, ds \; . 
\end{equation}
Note that (\ref{approHdef}) is a local time version of (\ref{XJHdef}). We refer to $H$ as 
the \textit{height process of $X$}. 
\begin{rem}
\label{supercr} Let us mention that in Le Gall \& Le Jan \cite{LGLJ98} and 
Le Gall \& D.~\cite{DuLG02}, the height process $H$ is introduced only for critical and subcritical spectrally positive processes. However, it easily extends to supercritical cases thanks to the following 
fact that is proved e.g.~in Bertoin's book \cite{Be96} Ch.~VII:  denote by $\bD ([0, \infty), \bbR)$ the space of c\`adl\`ag functions equipped with Skorokhod's topology; then for all $t\ino [0, \infty)$ and for all nonnegative measurable functional $F \! : \! \bD ([0, \infty), \bbR) \! \rightarrow \! \bbR$, 
\begin{equation}
\label{supsub}
\bE \big[ F( X_{\cdot \wedge t} ) ] \! =\! \bE \big[ \exp (\varrho \overline{X}_t ) \,  F( \overline{X}_{\cdot \wedge t} )\big]  , 
\end{equation}
where $\overline{X}$ stands for a \textit{subcritical} spectrally L\'evy process with Laplace exponent $\psi (\varrho + \cdot)$. \cq 
\end{rem}

Our main result, given below, introduces the analogue of $\cH^{\bw}$ in the continuous setting. 
\begin{thm}   
\label{cHdefthm} Let $(\alpha, \beta, \kappa, \mathbf{c})$ be as in (\ref{parconing}). 
Assume that (\ref{contH}) holds, which implies (\ref{varinfinie}). Recall from (\ref{AetYdef}) the definition of $Y$ and $A$; recall from (\ref{thetabdef}) the definition of $\theta^\mathtt{b}$; recall from (\ref{T*def}) the definition of  $T^*$; recall from (\ref{Xdef}) the definition of $X$ and recall from (\ref{approHdef}) how $H$ is derived from $X$.

Then, there exists a continuous process $(\cH_t)_{t\in [0, \infty)}$ such that for all $t\ino [0, \infty)$, $\cH_t$ is a.s.~equal to a measurable functional of $(Y_{\cdot \wedge t}, A_{\cdot \wedge t}) $ and such that 
\begin{equation}
\label{cHfromX}
\textrm{a.s.}\; \forall\,t\ino [0, T^*) , \quad \cH_t \! = \! H_{\theta^{\mathtt{b}}_t } . 
\end{equation}    
We refer to $\cH$ as the {\rm height process associated with $Y$}. 
\end{thm}
\noi
\textbf{Proof.} See Section \ref{cHdefthmpf}. \cqfd  

\medskip

We next define the pinching times as in (\ref{pin1}): we first set 
\begin{equation}
\label{Jdef}
\forall t\ino [0, \infty), \quad  
\JJ_t \! = \! \inf_{s\in [0, t]} Y_s \; .
\end{equation} 
Then, conditionally given $Y$, let 
\begin{equation}
\label{Poisurcon}
\cP\! = \! \sum_{\; p\geq 1} \!  \delta_{(t_p, y_p) }
 \; \textrm{be a Poisson pt.~meas.~on $[0, \infty)^2$ 
 with intensity $ \kappa\un_{\{ 0< y <Y_t -J_t \}} \, dt\, dy$.} 
\end{equation}
Then, set  
\begin{equation}
\label{pinchset}
\Ptt\! = \! \big( (s_p, t_p) \big)_{p\geq 1} \quad \textrm{where} \quad s_p\! = \! \inf \big\{ s\ino [0, t_p] :  \inf_{u\in [s, t_p]} 
Y_u\! -\! J_u > y_p   \big\}, \; \,  p\! \geq \! 1. 
\end{equation}
In the next section, after some necessary setup on the coding of graphs, we will see that the triple $(Y, \cH, \Ptt)$ completely characterises the continuous version of the multiplicative graph. 
%

\subsubsection{Coding graphs}
\label{codgrasec}

Our primary goal here is to introduce the notions on measured metric spaces which will allow us to build the  continuum graph from the processes $(Y, \cH, \Ptt)$. We will actually consider a slightly more general situation, which permits to include the coding for discrete graphs as well. 

\medskip
\noi
\textbf{Coding trees.}
\label{codtreesec} First let us briefly recall how functions (not necessarily continuous) code trees. Let $h\! :\!  [0, \infty) \! \rightarrow \! [0, \infty)$ be c\`adl\`ag and such that 
\begin{equation}
\label{zetadef}
\zeta_h\! = \! \sup \{ t\ino [0,  \infty) \! : \! h(t) \! >\! 0 \}< \infty \; .
\end{equation}
For all $s, t\ino [0, \zeta_h)$, we set 
\begin{equation}
\label{pseudometric}
b_h(s,t)=\!\!\!\!\!\!\!\!\! \inf_{\quad r\in[s\wedge t,s\vee t]} \!\!\!\!\!\!\!\!\! h(r) \qquad {\rm and} \qquad d_h(s,t)=h(s)+h(t)-2b_h(s,t).
\end{equation}
Note that $d_h$ satisfies the four-point inequality: for all $s_1, s_2, s_3, s_4 \! \in \!  [0, \zeta_h)$, one has
\[d_h(s_1,s_2) + d_h(s_3, s_4) \! \leq \! \big(d_h(s_1, s_3) + d_h(s_2, s_4)\big) \! \vee \!  \big( d_h(s_1, s_4) + d_h(s_2, s_3)  \big).\] Taking $s_3\! = \! s_4$ shows that $d_h$ is a 
pseudometric on $[0, \zeta_h)$.  We then denote by $s\! \sim_h \! t$ the equivalence relation
 $d_h(s,t) \! =\! 0$ and we set 
\begin{equation}
\label{codef}
T_h= [0, \zeta_h) / \! \sim_h \; .  
\end{equation}
Then, $d_h$ induces a true metric on the quotient set $T_h$ that we keep denoting by $d_h$ and we denote by $p_h\! :\!  [0, \zeta_h) \!  \rightarrow \! T_h $ 
the \textit{canonical projection}. Note that  $p_h$ is not necessarily continuous.

\vspace{-2mm}

\begin{rem}
\label{pataques}
The metric space $(T_h, d_h)$ is tree-like but in general it is not necessarily connected or compact. However, we shall consider the following cases.  

\smallskip

\begin{compactenum}
\item[$(a)$]  $h$ takes finitely many values. 

\smallskip

\item[$(b)$] $h$ is continuous. 
\end{compactenum}

\smallskip

\noi
\textit{In Case $(a)$}, $T_h$ is not connected but it is compact; $T_h$ is in fact formed by a finite number of points. In particular, $\mathcal H^\bw$ is in this case: by (\ref{gloupii}),  the exploration tree $\cT_{\!\! \bw}$ as defined in (\ref{explotree}) is actually isometric to $T_{\cH^\bw}$, that is the tree coded by the height process $\cH^{\bw}$ that is derived from $Y^\bw$ by (\ref{JHdef}). 

\smallskip

\noi
\textit{In Case $(b)$}, $T_h$ is compact and connected. Recall that real trees are metric spaces such that any pair of points is joined by a unique injective path that turns out to be a geodesic and recall that real trees are exactly the connected metric spaces satisfying the four-point condition (see Evans \cite{Ev08} for more references on this topic). Thus, $T_h$ is a compact real tree in this case. \cq
\end{rem}

The coding function provides two additional features: a distinguished point $\rho_h\! = \! p_h (0)$ that is called the \textit{root} of $T_h$ and the 
\textit{mass measure} $m_h$ that is the pushforward measure of the Lebesgue measure on $[0, \zeta_h)$ induced by $p_h$ on $T_h$: for any Borel measurable function $f \! : \! T_h \! \rightarrow \!  [0, \infty) $, 
\begin{equation}
\label{masmea}
 \int_{T_h} \!\!\!   f(\sigma ) \, m_h (d\sigma) = \! \int_0^{\zeta_h}  \!\!\! f(p_h(t)) \, dt \; .    
\end{equation}

\medskip

\noi
\textbf{Pinched metric spaces.}
\label{pinmetpar} 
We next briefly explain how additional ``shortcuts''  will modify the metric of a graph.
Let $(E, d)$ be a metric space and let $\Ptt\! = \! ((x_i, y_i))_{1\leq i\leq p}$ where $(x_i, y_i)\ino E^2$, $1\! \leq \! i\! \leq \! p$, are pairs of \textit{pinching points}. Let $\epp \ino [0, \infty)$ that is interpreted as the length of the edges that are added to $E$ (if $\epp \! = \! 0$, then each $x_i$ is identified with $y_i$). 
Set $A_E \! =\! \{ (x,y) ; x,y\ino E\}$ and for all $e\! = \! (x,y)\ino A_E$, 
set $\underline{e}\! = \! x$ and $\overline{e}\! = \! y$. A path $\gamma$ joining $x$ to $y$ is a sequence of $e_1, \ldots , e_q \ino A_E$ such that $\underline{e}_1\! = \! x$, $\overline{e}_q\! = \! y$ and $\overline{e}_i\! = \! \underline{e}_{i+1}$,  for all $1\! \leq \! i \! < \! q$. 
For all  $e\! = \! (x, y)\ino A_E$, we then define its length by $l_e \! =\!  \epp \! \wedge \! d(x_i, y_i)$ if $(x,y)$ or $(y,x)$ is equal to some $(x_i,y_i)\ino \Ptt$; otherwise we set $l_e \! = \! d(x, y)$. The length of a path $\gamma\! = \! (e_1, \ldots , e_q)$ is given by 
$l(\gamma)\! = \! \sum_{1\leq i\leq q} l_{e_i}$, and we set:  
\begin{equation}
\label{pinchdist}
 \forall x, y\ino E, \quad d_{\Ptt, \epp}  (x,y)= \inf \big\{ l (\gamma) ;\;  \textrm{$\gamma$ is a path joining $x$ to $y$} \big\} \; .
\end{equation}
It can be readily checked by standard arguments that $d_{\Ptt, \epp}$ is a pseudo-metric. We refer to Section~\ref{Pinfrac} for more details. 
Denote the equivalence relation $d_{\Ptt, \epp}  (x,y)\! = \! 0$ by $x \equiv_{\Ptt, \epp} y$; 
the \textit{$(\Ptt, \epp)$-pinched metric space associated with $(E, d)$} is then the quotient space $E/\!\! \equiv_{\Ptt, \epp}$ equipped with $d_{\Ptt, \epp}$. First note that if $(E, d)$ is compact or connected, so is the associated $(\Ptt, \epp)$-pinched metric space since the canonical projection $\varpi_{\Ptt, \epp}\! : \! E \! \rightarrow \!  E/\!\! \equiv_{\Ptt, \epp}$ is $1$-Lipschitz. Of course when $\epp \! >\! 0$, $d_{\Ptt , \epp}$ on $E$ is a true metric, $E\! = \!  E/\!\! \equiv_{\Ptt, \epp}$ and $\varpi_{\Ptt, \epp}$ is the identity map on $E$. 

\medskip

\noi
\textbf{Coding pinched trees.} 
\label{codpintre} Let $h\! :\!  [0, \infty) \! \rightarrow \! [0, \infty)$ be a c\`adl\`ag function that satisfies (\ref{zetadef}) and $(a)$ or $(b)$ 
in Remark \ref{pataques}; let $\Pi\! = \! ((s_i, t_i))_{1\leq i \leq p}$ where $0 \leq s_i \! \leq \! t_i \! < \! \zeta_h$, for all $1\! \leq \! i \! \leq \! p$ and let $\epp \ino [0, \infty) $. 
Then, the \textit{compact measured metric space coded by $h$ and the pinching setup $(\Pi, \epp)$} is the $(\Ptt, \epp)$-pinched metric space associated with $(T_h, d_h)$ and the pinching points 
$\Ptt\! = \!  ((p_h (s_i), p_h (t_i)))_{1\leq i\leq p}$, where $p_h \! : \! [0, \zeta_h ) \! \rightarrow \! T_h$ stands for the canonical projection. We shall use the following notation:  
\begin{equation}
\label{defgrgraf}
G (h, \Pi, \epp)= \big( G_{h, \Pi, \epp} , d_{h, \Pi, \epp}, \varrho_{h, \Pi, \epp} , m_{h, \Pi, \epp} \big).
\end{equation}
We shall denote by $p_{h, \Pi, \epp}$ the composition of the canonical projections 
$\varpi_{\Ptt, \epp} \circ  p_h \! :  \! [0, \zeta_h) \! \rightarrow \! G_{h, \Pi, \epp}$;  then 
$\varrho_{h, \Pi, \epp}\! = \! p_{h, \Pi, \epp} (0)$ and 
$m_{h, \Pi, \epp}$ stands for the pushforward measure of the Lebesgue on $[0, \zeta_h )$ via $p_{h, \Pi, \epp}$. 
%
%
%
%
%
%
%
%
%
\medskip

\noi   
\textbf{Coding $\bw$-multiplicative graphs.} 
\label{discrICRG}
Recall from Section \ref{conGw1} that $(\cG_k^\bw)_{1 \leq  k  \leq  \bq_\bw}$, are the connected components of $\cG_\bw$. Here, $\bq_\bw$ is the total number of connected components of  $\cG_\bw$; $\cG_{\! k}^\bw$ 
is equipped with its graph-metric $d_{\cG_{\! k}^\bw}$ and with the restriction $\bm^\bw_k$ 
of the measure $\bm_\bw \! = \! \sum_{1\leq j\leq n} w_j \delta_{j}$ on $\cG_{\! k}^\bw$;  
the indexation satisfies $\bm_\bw (\cG_{\! 1}^\bw )\! \geq \! \ldots \!  \geq \! \bm_\bw (\cG_{\! \bq_\bw}^\bw )$. Let us briefly explain how the excursions of $\cH^\bw$ above $0$ 
code the measured metric spaces $\cG_{\! k}^\bw$. 

First, denote by $(l^\bw_k, r^\bw_k)$, $1\! \leq \! k \! \leq \! \bq_\bw$, the excursion intervals of $\cH^\bw$ above $0$, that are exactly the excursion intervals of $Y^\bw$ above its infimum process $J^\bw_t \! = \! \inf_{s\in [0, t]} Y^\bw_s$. Namely, 
\begin{equation}
\label{excdibor}
 \bigcup_{1\leq k\leq \bq_\bw} [l^\bw_k, r^\bw_k) = \big\{  t\ino [0, \infty) : \cH^\bw_t >0  \big\}= \big\{  t\ino [0, \infty) : Y^\bw_t > \JJ^\bw_t  \big\}\,. 
\end{equation} 
Here, we set $\zeta^\bw_k\! = \! r^\bw_k \! -\! l^\bw_k\! = \! \bm^\bw_k (\cG_{\! k}^\bw)$ and thus 
$\zeta^\bw_1 \! \geq \! \ldots \! \geq \! \zeta^\bw_{\bq_\bw}$; moreover, if $\zeta^\bw_k\! = \! \zeta^\bw_{k+1}$, then we agree on the convention that $l^\bw_k \! < \! l^\bw_{k+1}$; \textit{excursions processes} are then defined as follows: 
\begin{equation}
\label{excdiscret}
 \forall k \ino \{ 1, \ldots , \bq_\bw\}, \; \forall t\ino [0, \infty), \qquad 
\Htt^{\bw}_k (t)= \cH^{\bw}_{(l^\bw_k + t)\wedge r^\bw_k} \quad \textrm{and} \quad  \Ytt^{\bw}_k (t)= Y^{\bw}_{(l^\bw_k + t)\wedge r^\bw_k}- \JJ^\bw_{l^\bw_k} . 
\end{equation}
We next define the sequences of \textit{pinching points of the excursions}: to that end, 
recall from (\ref{Poissurpl}) and (\ref{pin1}) 
the definition of $\Ptt_\bw\! = \! \big( (s_p, t_p)\big)_{1\leq p\leq \bp_\bw}$ the sequence of 
pinching points of $\cG_\bw$; observe that if $t_p \ino [l^\bw_k, r^\bw_k]$, 
then $s_p \ino [l^\bw_k, r^\bw_k]$; then, it allows to define the following for all $k\ino \{ 1, \ldots, \bq_\bw\}$: 
\begin{multline}
\label{defPiwk}
\Ptt_{k}^{\bw}\! = \! \big( (s^k_p, t^k_p)\big)_{1\leq p\leq \bp^\bw_k} \; \textrm{where $(t^k_p)_{1\leq p\leq \bp^\bw_k}$ increases and where} \\ \textrm{ the $(l^\bw_k +s^k_p, l^\bw_k+t^k_p)$ are exactly the terms $(s_{p^\prime} , t_{p^\prime})$ of $\Ptt_\bw$ such that $t_{p^\prime}\ino [l^\bw_k, r^\bw_k]$.}
\end{multline}
%
%
%
%
Then, for all $k\ino \{ 1, \ldots , \bq_\bw \}$, we easily see that $ \cG_{\! k}^\bw$ is coded by $( \Htt^{\bw}_k,  \Ptt_k^{ \bw}, 1)$ as defined in (\ref{defgrgraf}). Namely, 
\begin{equation}
\label{disctry}
G ( \Htt^{\bw}_k,  \Ptt_k^{ \bw}, 1) \; \, \textrm{is isometric to} \; \, \cG_{\! k}^\bw \; .
\end{equation}
Here, \textit{isometric} means that there is a bijective isometry from $G( \Htt^{\bw}_k,  \Ptt_k^{ \bw}, 1)$ onto 
$\cG_{\! k}^\bw$ sending $\bm^\bw_k$ to $\bm_{\Htt^{\bw}_k ,  \Ptt_k^{ \bw}, 1}$. 
This implies in particular that $\zeta_k^\bw\! = \! \bm_{k}^\bw (\cG^\bw_k)$.

\subsubsection{The continuous multiplicative random graph. Fractal properties.} 
\label{defICRG} 
We fix $(\alpha, \beta, \kappa, \mathbf{c})$ as in (\ref{parconing}) 
and we assume that (\ref{contH}) holds true. By analogy with the discrete coding, we now define the \textit{$(\alpha, \beta, \kappa, \mathbf{c})$-continuous multiplicative random graph}, the continuous version of $\bw$-multiplicative graph. 

First, recall from (\ref{AetYdef}) the definition of $Y$; recall from Theorem \ref{cHdefthm} 
the definition of $\cH$, the height process associated with $Y$; 
recall from (\ref{Jdef}) the notation, $J_t \! = \! \inf_{s\in [0,t]} Y_s$, $t\ino [0, \infty)$. 
Lemma \ref{AHeuer} (see further in Section \ref{cHdefthmpf}) 
asserts that the excursion intervals of $\cH$ above $0$ and the excursion intervals of $Y\! -\! J$ above $0$ are 
the same; moreover Proposition 14 in Aldous \& Limic \cite{AlLi98} (recalled further in Proposition \ref{AldLim1}, Section \ref{proXYsec}), asserts that $J_t \! \rightarrow -\infty$ and that 
these excursions can be indexed in the decreasing order of their lengths. Namely, 
\begin{equation}
\label{excuHY}
\big\{  t\ino [0, \infty) : \cH_t >0  \big\}= \big\{  t\ino [0, \infty) : Y_t > \JJ_t  \big\}= \bigcup_{k\geq 1} (l_k , r_k) \,,
\end{equation}
where the sequence $\zeta_k\! = \! l_k \! -\! r_k$, $k\! \geq \! 1$ decreases. This proposition also implies that 
$\{ t\ino [0, \infty)\! : \! \cH_t \! = \! 0\}$ has no isolated point, that
$\bP (\cH_t \! = \! 0)\! = \! 0$ for all $t\ino [0, \infty)$ and that the continuous function $t \! \mapsto \! -J_t$ 
can be viewed as a sort of local-time for the set of zeros of $\cH$. 
We refer to Sections \ref{proXYsec} and \ref{cHdefthmpf} for more details. These properties allow to 
define the \textit{excursion processes} as follows. 
\begin{equation}
\label{exccont}
 \forall k \! \geq \! 1, \; \forall t\ino [0, \infty), \qquad 
\Htt_{k}(t)= \cH_{(l_k + t)\wedge r_k} \quad \textrm{and} \quad  \Ytt_{k}(t)= Y_{(l_k + t)\wedge r_k}- \JJ_{l_k} . 
\end{equation}
The \textit{pinching times} are defined as follows: recall from (\ref{Poisurcon}) and (\ref{pinchset}) 
the definition of $\Ptt\! = \! \big( (s_p, t_p)\big)_{p\geq 1}$. 
If $t_p \ino [l_k, r_k]$, then note that $s_p \ino [l_k, r_k]$, by definition of $s_p$. 
For all $k\! \geq \! 1$, we then define:  
\begin{multline}
\label{defPik}
\Ptt_{k}\! = \! \big( (s^k_p, t^k_p)\big)_{1\leq p\leq \bp_k} \; \textrm{where $(t^k_p)_{1\leq p\leq \bp_k}$ increases 
and where} \\ \textrm{ the $(l_k +s^k_p, l_k+t^k_p)$ are exactly the terms $(s_{p^\prime} , t_{p^\prime})$ of $\Ptt$ such that $t_{p^\prime}\ino [l_k, r_k]$.}
\end{multline}
%
%
%
The \textit{connected components of the $(\alpha, \beta, \kappa, \mathbf{c})$-continuous 
multiplicative random graph} are then 
defined as the sequence of random compact measured metric spaces coded by the excursions $\Htt_k$ and the pinching setups $(\Ptt_k, 0)$. Namely, we shall use the following notation: for all $k\! \geq \! 1$, 
\begin{equation}
\label{defCRG}
\mathbf{G}_{k}:=  \big( \mathbf{G}_{k} , \mathrm{d}_{k}, \varrho_{k} , \bm_{k} \big) \; \textrm{stands for} \; G (\Htt_k, \Ptt_k, 0) \; \textrm{as defined by (\ref{defgrgraf}).}
\end{equation}

The above construction of the continuous 
multiplicative random graph via the height process $\cH$ highlights the intimate connection between the graph and  the L\'evy trees that are the continuum trees coded by the excursions of $H$ above $0$, where $H$ is the height process associated with the L\'evy process $X$. We conclude this section with a study of the fractal properties of the graphs where this connection with the trees plays an essential role.  
Indeed, as a consequence of Theorem \ref{cHdefthm},
each component $\mathbf{G}_k$ of the graph can be 
\textit{embedded in a L\'evy tree} whose branching mechanism $\psi$ is derived from $(\alpha, \beta, \kappa, \mathbf{c})$ by (\ref{psidefi}); roughly speaking 
the measure $\bm_k$ is the restriction to $\mathbf{G}_k$ of the \textit{mass measure} of the L\'evy tree; this measure enjoys specific fractal 
properties and as a consequence of Theorem 5.5 in Le Gall \& D.~\cite{DuLG05}, we 
get the following result.  
\begin{prop}
\label{fractCIRG} Let $(\alpha, \beta, \kappa, \mathbf{c})$ as in (\ref{parconing}).  
Assume that (\ref{contH}) holds true, which implies (\ref{varinfinie}). Let $(\mathbf{G}_k)_{k\geq 1}$ 
be the connected components of the continuous $(\alpha, \beta, \kappa, \mathbf{c})$-multiplicative random graph as defined in (\ref{defCRG}). We denote by $\mathrm{dim}_{\mathrm{H}}$ the Hausdorff dimension and by $\mathrm{dim}_{\mathrm{p}}$ the packing dimension. Then, the following assertions hold true a.s.~for all $k\! \geq \! 1$,

\smallskip

\begin{compactenum}
\item[$(i)$] If $\beta \! \neq \! 0$, then $\mathrm{dim}_{\mathrm{H}} (\mathbf{G}_k)\! = \! \mathrm{dim}_{\mathrm{p}} (\mathbf{G}_k)\! = \! 2 $. 

\smallskip

\item[$(ii)$] Suppose $\beta \! = \! 0$, which implies $\sigma_2 (\mathbf{c})\! = \! \infty$ by (\ref{varinfinie}). Then, set:  
$$\hspace{-8mm}\forall x\ino (0, 1) , \quad J (x) \! =\!  \frac{_1}{^x}\! \! \sum_{j: c_j \leq x }\!\!\!  \kappa c_j^3 \; + \! \sum_{j: c_j > x } \!\!\! \kappa c_j^2 \, = \sum_{j\geq 1} \, \kappa c^2_j \,  \big(1 \!  \wedge  \! (c_j/x) \big) $$ 
that tends to $\infty$ as $x\! \downarrow \! 0$. We next define the following exponents:  
\begin{multline}
\label{expone} 
\gamma  =  1+  \sup \!  \big\{ r \ino [0, \infty)\!  :\!   \lim_{\;\;  x\rightarrow 0+} \! \!\! x^{r} \! J(x) \! = \! \infty  \big\}  \\ 
 \textrm{and} \quad  \eta  =  1+  \inf \!  \big\{ r \ino [1, \infty) \! : \!  \lim_{\;\;  x\rightarrow 0+} \! \!\! x^{r} \! J(x)  \! = \! 0  \big\}. \quad \quad 
\end{multline}
In particular, if $\mathbf c$ varies regularly with index 
$-\rho^{-1}\ino (-\frac{1}{2}, -\frac{1}{3})$, then $\gamma\! =\! \eta\! =\! \rho \! -\! 1$.

\noi
Then, if $\gamma \! >\! 1$, we a.s.~get 
$$\mathrm{dim}_{\mathrm{H}} (\mathbf{G}_k)\! = \frac{\gamma}{\gamma-1}\quad  \textrm{and} \quad  
\mathrm{dim}_{\mathrm{p}} (\mathbf{G}_k)\! = \! \frac{\eta}{\eta -1} \; .$$ 
\end{compactenum}
\end{prop}
\noi
\textbf{Proof.} See Section \ref{embedGinT}. \cqfd 


\subsection{Limit theorems for the multiplicative graphs}

%
 

%
%
%
%
The family of continuous multiplicative random graphs, as defined previously, appears as the scaling limits of the $\bw$-multiplicative graphs. This convergence is proved in \cite{BDW2}. We give here a short description of the main results of this paper. 

Recall from Section \ref{sec: dis} the various coding processes for the discrete graphs: $X^{\bw}, H^{\bw}$ for the Markovian queue/Galton--Watson tree, $Y^{\bw}$ and $\cH^{\bw}$ for the graph, and $\theta^{\mathtt{b}, \bw}$ is the time-change so that $\cH^{\bw}=H^{\bw}\circ\theta^{\mathtt{b}, \bw}$. Then recall from \ref{sec: con} their analogues $X, H, Y, \cH, \theta$ in the continuous setting. 
Let $\bw_n\! = \! (w^{_{(n)}}_{^j})_{j\geq 1}\in  \elldo_{{\! f}}, n\ge 1$, be a sequence of weights. The asymptotic regime considered in \cite{BDW2} is determined by two sequences of positive numbers $a_{n}, b_{n}\to\infty$. More precisely,  we prove in \cite{BDW2} (Theorem 2.4) the following convergence of the coding processes:
\begin{equation}
\label{glurniglup}
\textit{If} \quad \big( \frac{_{_1}}{^{^{a_n}}} X^{\bw_n}_{b_n \cdot }\,  , \frac{_{_{a_n}}}{^{^{b_n}}} H^{\bw_n}_{b_n \cdot } \big)\;  
\underset{n\rightarrow \infty}{-\!\!\! -\!\!\! -\!\!\! \longrightarrow} \; \big( X,  H \big) 
\end{equation}
\textit{weakly on $\bD([0, \infty), \bbR)\! \times \! \bC ([0, \infty) , \bbR)$ equipped with the product of the Skorokhod and the continuous topologies, then the following joint convergence }
 \begin{equation}
\label{hurkoko}
\big( \frac{_{_1}}{^{^{a_n}}} X^{\bw_n}_{b_n \cdot }\,  , \frac{_{_{a_n}}}{^{^{b_n}}} H^{\bw_n}_{b_n \cdot } \, ,  
\frac{_{_1}}{^{^{a_n}}} Y^{\bw_n}_{b_n \cdot } 
,  \frac{_{_{a_n}}}{^{^{b_n}}} \cH^{\bw_n}_{b_n \cdot}  \big) 
\;  \underset{n\rightarrow \infty}{-\!\!\! -\!\!\! -\!\!\! \longrightarrow} \; \big( X, H, 
Y , \cH \big) 
\end{equation}
\textit{holds weakly on $\bD([0, \infty), \bbR)\! \times \! \bC ([0, \infty) , \bbR)\! \times \!  \bD([0, \infty), \bbR)\! \times \! 
\bC ([0, \infty) , \bbR) $ equipped with the product topology.} 

For \eqref{glurniglup} to hold, it is necessary and sufficient that (see Le Gall \& D.~\cite{DuLG02})
\begin{equation}
\label{hurkiki}
 (A):  \frac{_{_1}}{^{^{a_n}}} X^{\bw_n}_{b_n } \overset{\textrm{(weakly)}}{\underset{n\rightarrow \infty}{-\!\!\! -\!\!\! -\!\!\! \longrightarrow}} X_1 \ \textrm{ and } \  (B): \ 
 \exists\, \delta \in \! (0, \infty) , \quad \liminf_{n\rightarrow \infty} \bP \big( Z^{\bw_n}_{\lfloor b_n \delta /a_n \rfloor} = 0 \big) >0 \;. 
 \end{equation} 
Above, $(Z^{\bw_n}_k)_{k\in \bbN}$ stands for a Galton--Watson Markov chain with offspring distribution $\mu_{\bw_n}$ (as defined in (\ref{mupoissw})) and with initial state $Z^{\bw_n}_0\! = \! \lfloor a_n \rfloor$. In \cite{BDW2}, for each $X$ satisfying the conditions \eqref{parconing} and \eqref{contH}, we construct examples of $(a_n, b_n, \bw_n)_{n\in \bbN}$ for  (\ref{hurkiki}) to take place, so that every member of the continuous multiplicative graph family appears in the limits of discrete graphs. Our construction there employs a sufficient condition for $(B)$ in (\ref{hurkiki}) which can be readily checked in a wide range of situations. Moreover, in the near critical regime where $\sigma_1 (\bw_n) \asymp \sigma_2 (\bw_n)$ and when large weights persist to the limit, \cite{BDW2} shows that $(\alpha, \beta, \kappa, \mathbf{c})$-spectrally L\'evy processes 
are the only possible scaling limits of the $X^\bw$. Thus, in some sense,  $(\alpha, \beta, \kappa, \mathbf{c})$-continuum multiplicative graphs as introduced in this paper, are the only possible non-degenerate scaling limits of near critical multiplicative graphs.

In terms of the convergence of the graphs, (\ref{hurkoko}) implies the following. Recall from \eqref{disctry} the graph $(\cG_{\! k}^{\bw_n}, d_{k}^{\bw_n}, \varrho_k^{\bw_n}, \bm_k^{\bw_n})$ encoded by 
$( \Htt^{\bw}_k,  \Ptt_k^{ \bw}, 1)$. Put another way,  $(\cG_{\! k}^{\bw_n})_{1\le k\le\mathbf{q}_{\bw_n}}$ forms the sequence of the connected components of $\cG_{\bw_n}$ sorted in decreasing order of their $\bm^{\bw_n}$-masses: $\bm^{\bw_n} (\cG_{\! 1}^{\bw_n}) \! \geq \! \ldots  \! \geq \! 
\bm^{\bw_n} (\cG_{\! \mathbf{q}_{\bw_n} }^{\bw_n})$. Analogously, recall from \eqref{defCRG} the connected components $(\mathbf G_{k})_{k\ge 1}$ of the $(\alpha, \beta, \kappa, \mathbf{c})$-continuous 
multiplicative random graph (ranked in decreasing order of their $\bm$-masses). 
Completing the finite sequence with null entries, we have that (Theorem 2.8 in \cite{BDW2})  
\begin{equation}
\label{hurkuku}
\big(\big( \cG_{\! k}^{\bw_n} ,  \frac{_{_{a_n}}}{^{^{b_n}}}d_{k}^{\bw_n} , \varrho_k^{\bw_n}, \frac{_{_{1}}}{^{^{b_n}}}\bm_k^{\bw_n}  \big) \big)_{k\geq 1}
\;  \underset{n\rightarrow \infty}{-\!\!\! -\!\!\! -\!\!\! \longrightarrow} \; \big(\big( \mathbf{G}_{k} , \mathrm{d}_{k}, \varrho_{k} , \bm_{k} \big) \big)_{k\geq 1}
\end{equation} 
holds weakly with respect to the product topology induced by the Gromov--Hausdorff--Prokhorov distance. Moreover, the same type of convergence holds if we replace in (\ref{hurkuku}) the weight-measure $\bm^{\bw_n}$ by the counting measure $\# \! = \! \sum_{1\leq j \leq \mathbf{j}_n} \delta_j$, where 
$\mathbf{j}_n \! := \! \sup \{ j\! \geq \! 1 \! : \! w^{_{(n)}}_{^j}\! >\! 0\}$. We can also list the connected components of $\cG_{\bw_{n}}$ in the decreasing order of their numbers of vertices. 
Then again, Theorem 2.13 in \cite{BDW2} asserts that under the additional assumption that $\sqrt{\mathbf{j}_n}/ b_n \! \rightarrow \! 0$, we have this (potentially) different list of connected components converging to the same limit object. Theorem~2.8 in the companion paper \cite{BDW2} completes the previous scaling limits of truely inhomogeneous graphs due to: Bhamidi, Sen \& X.~Wang \cite{BhSeWa14}, Bhamidi, Sen, X.~Wang \& B.~\cite{BhBrSeWa14}, Bhamidi, van der Hofstad \& van 
Leeuwaarden \cite{BhHoLe12} and 
Bhamidi, van der Hofstad \& Sen \cite{BhHoSe15}. We refer to \cite{BDW2} for a more precise discussion on the previous works dealing with the scaling limits of inhomogeneous graphs.

\section{Proof of Theorem~\ref{loiGG}.}
\label{PfThmloiGG}
 Let $G\! = \! (\cV(G), \ccE (G))$ be a graph with $\cV (G) \! \subset \! \bbN\backslash \{ 0\}$.
We suppose that $G$ has $q$ connected components $\mathtt{C}_{^1}^{_G}, \ldots ,\mathtt{C}_{^q}^{_G}$ that are listed in an arbitrary (deterministic way): let us say that \textit{the components are listed in the increasing order of their least vertex}; namely, $\min \mathtt{C}_{^1}^{_G} \! < \! \ldots \! <\!  \min \mathtt{C}_{^q}^{_G}$. 
Let $\bw\! = \! (w_j)_{j\in \cV(G)}$ be a set of strictly positive weights; we set 
$\rbm \! = \! \sum_{j\in \cV (G)} w_j \delta_j$ that is a measure on $\cV(G)$. 

We call the set of the atoms of a Poisson point measure \textit{a Poisson random set}. 
We then define a law $\Lambda_{G, \rbm}$ on $((0, \infty) \! \times \!  \cV(G))^q$ as follows.  
Let $(\Pi_j)_{j\in \cV(G)}$ be independent Poisson random subsets of $(0, \infty)$ with rate $w_j/ \sigma_1 (\bw)$   
for all non-empty subset $S \! \subset \! \cV(G)$ we set $ \Pi (S)\! := \! \bigcup_{j\in S} \Pi_j $ (in particular  $\Pi (\{ j\}) \! = \! \Pi_j$, for all $j\ino \cV (G)$). 
Then $\Pi (S)$ is a Poisson random set with rate $\rbm (S)/ \sigma_1(\bw)$. 
For all $k \ino  \{ 1, \ldots , q\}$, we then 
define $(T_k , U_k) \! : \! \Omega \! \rightarrow \! (0, \infty) \! \times \!  \cV(G)$ by:  
$$ 
T_k \! = \! \inf \Pi (\mathtt{C}^{_G}_{^{\mathbf{s} (k)}} \big) \! = \! \inf \Pi_{U_k} \; \textit{where the permutation $\mathbf{s}$ is such that} \;  \inf \Pi \big(\mathtt{C}^{_{G}}_{^{\bs (1)}} \big) \! <\!  \ldots  \! < \!  \inf \Pi \big( \mathtt{C}^{_{G}}_{^{\bs (q)}} \big).
$$
Namely, $T_k$ is the $k$-th order statistic of $( \inf \Pi (\mathtt{C}^{_G}_{^1}), \ldots ,  \inf \Pi ( \mathtt{C}^{_G}_{^q} ))$. 
%
We denote by $\Lambda_{G, \rbm}$ the joint law of $((T_k, U_k))_{1\leq k \leq q}$ and we easily 
check:   

\begin{eqnarray}
\label{Lambformu}
\hspace{-10mm}\Lambda_{G, \rbm} ( dt_1 \ldots dt_q  \!\!\!\! \!\!   & ;& \!\!\!\! \!\! j_1, \ldots , j_q)=  \bP \big( T_1\ino dt_1; \ldots ; T_q \ino dt_q; U_1 \! =\!  j_1; \ldots ; U_q \! = \! j_q \big) \nonumber \\
\!\!\!\! & =&\!\!\!\!  \un_{\{ 0\leq t_1 \leq \ldots \leq t_q \} } \frac{w_{j_1}}{\sigma_1 (\bw)} \ldots  \frac{w_{j_q}}{\sigma_1 (\bw)} \exp \Big(\! \!  -\!  \frac{_1}{^{\sigma_1 (\bw)}}\!\! \sum_{\quad 1 \leq k \leq q} \!\!\! t_k  \rbm \big( \mathtt{C}^{_G}_{^{s (k)}}\big) \Big) dt_1 \ldots dt_q \,,
\end{eqnarray}
where $s$ is the unique permutation of $\{1, \ldots , q\}$ such that $j_l \ino \cV(\mathtt{C}^{_G}_{^{s(l)}})$, for all $l\ino \{ 1, \ldots, q\}$. The following lemma, whose elementary proof is left to the reader, 
provides a description of the law of $((T_k, U_k))_{2\leq k \leq q }$ conditionally given 
$(T_1, U_1)$. Let us mention that it is formulated with specific notation for further use. 
\begin{lem}
\label{1compopo} Let $G^o$ be a finite graph with $q^o$ connected components; let $\bw^o\! = \! (w^o_j)_{j \in \cV(G^o)}$ be strictly positive weights; let $\rbm^o\! = \! \sum_{j \in \cV(G^o)} w^o_j \delta_j$. We fix $k^* \ino \{ 1, \ldots , q^o\} $.  
Then, we set $G^\prime\! = \! G^o \backslash \mathtt{C}_{^{k^*}}^{_{G^o}}$ 
and $a=(\sigma_{1}(\bw^{o})-\rbm^{o}(\mathtt{C}_{^{k^*}}^{_{G^o}}))/\sigma_{1}(\bw^{o})$; we equip $G^\prime$  with the set of weights $w^\prime_j\! = \! a w^o_j$, $j\ino \cV(G^\prime)$ and we set 
$\rbm^\prime\! = \! \sum_{j \in \cV(G^\prime)} w^\prime_j \delta_j$. 
Let $T$ and $(T^\prime_k, U^\prime_k)_{1\leq k\leq q^o-1}$ be independent r.v.~such 
that $T$ is exponentially distributed with unit mean and such that 
$(T^\prime_k, U^\prime_k)_{1\leq k\leq q^o-1}$ has law $\Lambda_{G^\prime, \rbm^\prime}$. We set 
$$ T^o_1= T \quad \textrm{and} \quad \forall k \ino \{ 1, \ldots , q^o\! -\! 1\}, \quad T^o_{k+1}= T+ \frac{_{_1}}{^{^a}} T^\prime_{k} \quad \textrm{and} \quad U^o_{k+1} \! = \! U^\prime_k \; .$$  
Then, for all $j^*\ino \mathtt{C}_{^{k^*}}^{_{G^o}}$,
$$ \frac{w^o_{j^*}}{\sigma_1 (\bw^o)} \bP \big( T^o_1 \ino dt_1; \ldots ; T^o_{q^o} \ino dt_{q^o}; U^o_2\! =\! j_2; \ldots ; U^o_{q^o}\! = \! j_{q^o} \big) \! = \! 
\Lambda_{G^o\! , \rbm^o} (dt_1 \ldots dt_{q^o}; j^*\! , j_2, \ldots , j_{q^o} ) .$$  
\end{lem}

Next we briefly recall how to derive a graph from the LIFO queue (and from an additional point process) as discussed in Section \ref{YLIFOsec}. 
 Let $\cV \! \subset \! \bbN \backslash \{ 0 \} $ a finite set of vertices (or finite set of labels of clients) associated with a set of strictly positive weights denoted by 
$\bw\! = \! (w_j)_{i\in \cV}$ (recall that the total amount of service of Client $j$ is $w_j$); let 
$\underline{E}\! = \! (E_j)_{j\in \cV}$ be the times of arrival of the clients. We assume that the clients arrive at distinct times and that no client enters exactly when another client leaves the queue. 
These restrictions correspond to a Borel subset of 
$(0, \infty)^{\# \cV}$ for $\underline{E}$ that has a full Lebesgue measure. We next set 
\begin{equation}
\label{psoque}
Y_t \! = \! -t + \sum_{j\in \cV} w_j \un_{\{ E_j \leq t \}}
\quad \textrm{and} \quad J_t \! = \! \inf_{r\in [0, t]} Y_r \; .
\end{equation}
We then define $V\! : \! [0, \infty) \! \rightarrow \! \cV$ such that $V_t$ is the label of the client who is served at time $t$: since $Y$ only increases by jumps, for all $t\ino [0, \infty)$, we get the following: 
\begin{compactenum}
\item[$\bullet$] either $\{ s\ino [0, t]\! : \! Y_{s-} < \inf_{[s, t]}Y \}$ is empty and we set $V_t\! = \! 0$, 
\item[$\bullet$] or there exists $j\ino \cV$ such that $E_j \! = \! 
\max \{ s\ino [0, t]\! : \! Y_{s-} < \inf_{[s, t]}Y \}$ and we set $V_t\! = \! j$. 
\end{compactenum} 
Note that $V_t \! = \! 0$ if the server is idle. Observe that $V$ is c\`adl\`ag. 
As mentioned in Section \ref{YLIFOsec}, the LIFO-queue yields an exploration forest 
whose set of vertices is $\cV$ and whose set of edges are 
$$\cA\! = \! \big\{ \{i,j \}\! : \! i,j \ino \cV  \, \textrm{and} \;  V_{E_j-}= i  \big\} \; .$$ 
Additional edges are created thanks to a finite set of points 
$\Pi\! = \! \{ (t_p, y_p) \, ;  1\! \leq \! p \! \leq \! \mathbf{p}_\bw \} $ in $D\! =\! \{ (t, y)\ino (0, \infty)^2 \! : \! 0 \! < \! y \! < \! Y_t\! -\! J_t \}$ as follows. For all $(t, y)\ino D$, define 
$\tau (t,y) \! =\!  \inf \big\{ s\ino [0, t] : \inf_{u\in [s, t]} Y_u \!  >\!  y+ J_t \big\}$. Then, the set of additional edges is defined by 
$$\cS\! = \! \big\{  \{ i,j\} : 
 \! i,j \ino \cV  \, \textrm{distinct and $\exists (t,y) \ino \Pi$ such that} \; V_{\tau (t,y)}\! = \! i \; \textrm{and} \;  V_t \! = \! j \big\} \; .$$ 
Then the graph produced by $\underline{E}$, $\bw$ and $\Pi$ is $ \cG \! = \! \big( \cV (\cG)\! = \! \cV \, ;  \,  \ccE (\cG)\! = \! \cA \cup \cS  \big) $. 

\smallskip

Theorem \ref{loiGG} asserts that if $\underline{E}$ and $\Pi$ have the appropriate distribution, 
then $\cG$ is a $\bw$-multiplicative graph, whose law is denoted by 
$M_{\cV\! , \bw}$ given as follows: for all graphs $G$ such that $\cV(G)\! = \! \cV$, 
\begin{equation}
\label{gongorisme}
M_{\cV\! , \bw} (G) = \prod_{\{ i,j\} \in \ccE (G)} \big(1\! -\! e^{-w_jw_j/ \sigma_1 (\bw)} \big)  \prod_{\{ i,j\} \notin \ccE (G)} e^{-w_jw_j/ \sigma_1 (\bw)} \; , 
\end{equation}
We actually prove a result that is slightly more general than Theorem \ref{loiGG} and that involves additional features derived from the LIFO queue. More precisely, denote by $\bq$ the number of excursions of $Y$ strictly above its infimum and denote by $(l_k, r_k)$, $k \ino \{ 1, \ldots , \bq\}$ the corresponding excursion intervals listed \textbf{in the increasing order of their left endpoints: $l_1 \! < \! \ldots  \! < \! l_{\bq}$} (of course, this indexation does not necessarily coincide with the indexation of the connected components of the graph $\cG$ in the increasing order of their least vertex, nor with the decreasing order of their $\bm_\bw$-measure).  Then, we set:  
\begin{equation}
\label{syllepse} \forall k \ino \{ 1, \ldots , \bq \}, \quad T_k = -J_{l_k} \quad \textrm{and} \quad U_k \ino \cV \; \textrm{is such that} \; E_{U_k}= l_k \; .
\end{equation}

From the definition of $\cG$ as a deterministic function of $(\underline{E}, \bw, \Pi)$ as recalled above, we easily check the following: $\cG$ has $\bq$ 
connected components $\mathtt{C}_{^1}^{_{\cG}}, \ldots ,\mathtt{C}_{^{\bq}}^{_{\cG}}$ (recall that they are listed in the increasing order of their least vertex: namely, $\min \mathtt{C}_{^1}^{_{\cG} } \! < \! \ldots \! <\!  \min 
\mathtt{C}_{^{\bq}}^{_{\cG}}$). Then, we define the permutation $\mathbf{s}$ on $\{ 1, \ldots , \bq\}$ that satisfies $U_k \ino \mathtt{C}_{^{\mathbf{s} (k)}}^{_{\cG}}$ for all $k\ino \{ 1, \ldots , \bq\}$. Observe that 
$r_k \! -\! l_k \! = \!  \mathtt{m} \big(\mathtt{C}_{^{\mathbf{s} (k)}}^{_{\cG}} \big)$ and that 
the excursion $(Y_{t+l_k} \! -\! J_{l_k})_{t\in [0, r_k-l_k]}$ 
codes the connected component $\mathtt{C}_{^{\mathbf{s} (k)}}^{_{\cG}}$. 
The quantity $T_k$ is actually the total amount of time during which the server is idle before 
the $k$-th connected component is visited, and $U_k$ is the first visited 
vertex of the $k$-th component.  
We denote by $\Phi$ the (deterministic) function that associates $(\underline{E}, \bw, \Pi)$ to 
$(\cG, Y, J, (T_k,U_k)_{1 \leq k \leq \bq})$: 
\begin{equation}
\label{Phigoudas}
 \Phi \big(\underline{E}, \bw, \Pi \big)  = \big(\cG, Y, J, (T_k,U_k)_{1 \leq k \leq \bq} \big) \; .
 \end{equation}
We next prove the following theorem that implies Theorem \ref{loiGG}. 
\begin{thm}
\label{+th} We keep the notation from above. We assume that $\underline{E}\! = \! (E_j)_{j\in \cV}$ are independent exponentially distributed r.v.~such that $\bE [E_j]\! = \! \sigma_1 (\bw)/ w_j$, for all $j\ino \cV$.  
We assume that conditionally given $\underline{E}$, $\Pi$ is a Poisson random subset of 
$D\! =\! \{ (t, y)\ino (0, \infty)^2 \! : \! 0 \! \leq  \! y \! < \! Y_t\! -\! J_t \}$ with intensity $\sigma_1 (\bw)^{-1} \un_{D} (t,y) \,  dtdy$. Let $\big(\cG, (T_k,U_k)_{1 \leq k \leq \bq} \big)$ be derived from 
$\big(\underline{E}, \bw, \Pi \big)$ by (\ref{Phigoudas}).  
Then, for all graphs $G$ whose set of vertices is $\cV$ and that have $q$ connected components, we get 
\begin{eqnarray}
\label{+loi}
\bP (\cG\! = \! G \, ; \, T_1\ino dt_1; \ldots ; T_q \ino dt_q \!\!\! \!\! &; &  \!\!\!  \! U_1 \! =\!  j_1 ;  \ldots ; U_q \! = \! j_q \big) \nonumber \\ 
& & = M_{\cV\! , \bw} (G) \,  \Lambda_{G, \, \rbm} (dt_1, \ldots , dt_q; j_1, \ldots , j_q) \; .
\end{eqnarray}
where $M_{\cV\! , \bw}$ is defined by (\ref{gongorisme}) and $\Lambda_{G, \, \rbm}$ is defined by (\ref{Lambformu}). 
\end{thm}
\textbf{Proof.} We proceed by induction on the number of vertices of $\cG$. When $\cG$ has only one vertex, then (\ref{+loi}) is obvious. We fix an integer $n\! \geq \! 1$ and we 
assume that (\ref{+loi}) holds for all $\cV \! \subset \! \bbN \backslash \{ 0\}$ such that $\# \cV\! = \! n$ and all sets of positive weights $\bw\! = \! (w_j)_{j\in \cV}$. 

  Then, we fix $\cV^o \! \subset \! \bbN \backslash \{ 0\}$ such that $\# \cV^o\! = \! n+1$; 
we fix $\bw^o\! = \! (w^o_j)_{j \in \cV^o}$ that are strictly positive weights and 
we also fix $\underline{E}^o\! = \! (E_j^o)_{j\in \cV^o}$ in 
$(0, \infty)^{n+1}$; we assume that in the corresponding LIFO queue, clients arrive at distinct times and that no client enters exactly when another client leaves the queue. 
We next set 
$Y^o_t \! = \! -t+ \sum_{j\in \cV^o} w_j^o \un_{\{ E^o_j \leq t\}}$ and $J^o_t\! = \! \inf_{[0, t]} Y^o$. Let $\bq^o$ be the number of excursions of $Y^o$ strictly above its infimum process $J^o$; let 
$(T^o_k, U^o_k)_{1\leq k\leq \bq^o}$ be as in (\ref{syllepse}): 
namely $T^o_k\! =\!  -\! J^{o}_{{l^{o}_{k}}}$ and $E^{o}_{{U^o_k}} \! = \! l^o_k$, where $(l^o_k, r^o_k)$ is the 
$k$-th excursion interval of $Y^o$ strictly above $J^o$ listed in the increasing order of their 
left endpoint (namely,  $l^o_1 \! < \! \ldots  \! < \! l^o_{\bq^o}$). 

The main idea for the induction is to shift the LIFO queue at the time of arrival $T^o_1$ of the first 
client (with label $U^o_1$) 
and to consider the resulting graph. More precisely, we set the following. 
\begin{equation}
\label{truchement}
\cV\! := \! \cV^o\backslash \{ U^o_1\}, \; \;  a\! :=\!  \frac{\sigma_1 (\bw^o)\! -\! w^o_{U^o_1}}{\sigma_1 (\bw^o)} , \; \;  \forall j\ino \cV, \; w_j\! = \!a w^o_j \;\,  \textrm{and} \; \,  E_j \! =\!  a (E^o_j \! -\! T^o_1)  \; .
\end{equation}
Let $Y$ and $J$ be derived from $\underline{E}\! :=\! (E_j)_{j \in \cV}$ and $\bw\! := \! (w_j)_{j\in \cV}$ as in (\ref{psoque}). Then observe that 
\begin{equation}
\label{kermes}
a \big( Y^o_{ a^{-1}t+ T^o_1} -Y^{o}_{T^o_1} \big)= Y_t \! = \! -t + \sum_{j\in \cV} w_j \un_{\{ E_j \leq t \}} \; , \quad t \ino [0, \infty) \; .
\end{equation} 
Note that $T^o_1= \min_{j\in \cV^o} E^o_j$. Then, the alarm clock lemma implies the following. 

\smallskip

\begin{compactenum} 
\item[(I)] \textit{Suppose that  $(E^o_j)_{j\in \cV^o}$ 
are independent exponentially distributed r.v.~such that 
 $\bE [E^o_j] \! = \! \sigma_1 (\bw^o)/w^o_j$, for all $j\ino \cV^o$. Then, 
$T^o_1$ is an exponentially distributed r.v.~with unit mean, 
$\bP (U^o_1\! = j) \! = w^o_j / \sigma_1 (\bw^o)$, for all $j\ino \cV^o$, $T^o_1$ and $U^o_1$ are independent and under the conditional probability 
$\bP (\, \cdot \, | U^o_1\! = \! j^*)$, $(E_j)_{j\in \cV}$, as 
defined in (\ref{truchement}), are independent exponentially distributed r.v.~such that  $\bE [E_j] \! = \!  \sigma_1 (\bw)/w_j$, for all $j\ino \cV$ (where a.s.~$\cV \! = \! \cV^o\backslash \{ j^*\}$). }    
\end{compactenum}

\smallskip

\begin{figure}[htp]
\centering
\includegraphics[height=4cm]{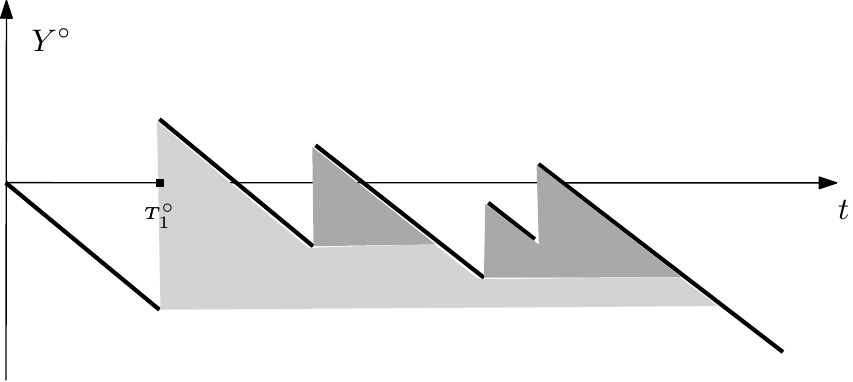}
\caption{{\small An example of $Y^{\circ}$ and $D^{\circ}$. Darkly shaded area corresponds to $f_{1}(D_{1})$; lightly shaded area corresponds to $f_{2}(D_{2})$. Together they form a partition of $D^\circ$}. \cq  }
\label{fig:Poisson} 
\end{figure}

We next fix $\cQ_1$ and $\cQ_2$, two discrete (i.e.~without limit-point) subsets of $[0, \infty)^2$  
that yield the 
additional edges in a specific way that is explained later; let us mention that eventually $\cQ_1$ and $\cQ_2$ are taken random and distributed as independent 
Poisson point processes whose intensity is Lebesgue measure on $(0, \infty)^2$. We first set: 
\begin{eqnarray}
\label{+Pi12}
D_1\! := \! \big\{ (t, y)\ino [0, \infty)^2 \! : \! 0 \! \leq  \! y \! < \! Y_t\! -\! J_t \big\} , \!\!\! & &\!\!\!  \Pi \! := \! \cQ_1 \cap D_1 \nonumber  \\  
 \!\!\! \!\!\!\!\!  & & \!\!\! \!\!\!\!\!  \!\!\! \!\!\!\!\!    \!\!\! \!\!\! \textrm{and} \quad 
D_2 \! = \!  \big\{ (t, y)\ino [0, \infty)^2 \! : \! -J_t \!  < \! y \! \leq  \! aw^o_{U^o_1}  \big\} , \; \, \Pi^2 \! := \! \cQ_2 \cap D_2 . 
\end{eqnarray}
We next define $f_1$ and $f_2$ from $[0, \infty)^2$ to $[0, \infty)^2$ and we define a set of points $\Pi^o$ by: 
\begin{eqnarray}
\label{+Pio12}
 f_1 (t,y)=  \big( \frac{_{_1}}{^{^a}}t +T^o_1 \, , \,  \frac{_{_1}}{^{^a}}\big( y + (aw^o_{U^o_1} + J_t)_+\big) \big) , & & \hspace{-10mm} f_2 (t, y)  = 
 \big( \frac{_{_1}}{^{^a}} t +T^o_1 \, , \,  w^o_{U^o_1} \! -\! \frac{_{_1}}{^{^a}}y \big)  
 \nonumber  \\  
& \textrm{ and} &  \quad    \Pi^o\! = \! f_1(\Pi) \cup f_2 (\Pi^2) \,.
\end{eqnarray}
 We check the following (see also Fig.~\ref{fig:Poisson}). 

\smallskip

\begin{compactenum} 
\item[(II)] \textit{Fix $\underline{E}^o$; suppose that $\cQ_1$ and $\cQ_2$ are two independent Poisson random subsets of $[0, \infty)^2$ with intensity $\frac{_1}{^{\sigma_1 (\bw)}} \, dtdy$. Then, $\Pi^o$ is a Poisson random subset of  
$D^o\! = \! \{ (t,y) \! \in \! [0, \infty)^2\! : \! 0 \! \leq  \!  y \! < \! Y^o_t\! -\! J^o_t \}$ with intensity $\frac{_1}{^{\sigma_1 (\bw^o)}} \un_{D^o} (t,y)\, dtdy$.  }    
\end{compactenum}

\smallskip

\noi
\textit{Indeed,} observe that $f_1 (D_1)$ and $f_2 (D_2)$ form a partition of $D^o$. Then, note 
that $f_1$ is piecewise affine with slope $1/a$ in both coordinates on the 
excursion intervals of $Y\! -\! J$ strictly above $0$, and that that $f_2$ is affine with slope $1/a$ in both coordinates. Standard results on Poisson subsets entail that $\Pi^o$ is a Poisson random subset on 
$D^o$ with intensity $\frac{_{a^2}}{^{\sigma_1 (\bw)}} \un_{D^o} (t,y)\, dtdy$ and by (\ref{truchement}) $\frac{_{a^2}}{^{\sigma_1 (\bw)}}\! = \! \frac{_{1}}{^{\sigma_1 (\bw^o)}}$, which implies (II). \cq 

\smallskip

\smallskip

Recall notation (\ref{Phigoudas}). We next introduce the two following graphs: 
\begin{equation}
\label{meloe}
 \Phi \big(\underline{E}^o, \bw^o, \Pi^o \big)  \! =\!  \big(\cG^o\! , Y^o\! , J^o\! , 
 (T^o_k\! ,U^o_k)_{1 \leq k \leq \bq^o} \big) \; \,  \textrm{and} \; \,   \Phi \big(\underline{E}, \bw, \Pi \big)  \! =\!  \big(\cG, Y, J, (T_k,U_k)_{1 \leq k \leq \bq} \big)  .
 \end{equation}
Then, the previous construction of $\Pi^o$ combined with (\ref{kermes}) easily implies the following:

\smallskip

\begin{compactenum} 
\item[(III)] \textit{Fix $\underline{E}^o, \cQ_1$ and $\cQ_2$. Then, 
$\cG$ is obtained by removing the vertex $U^o_1$ from $\cG^o\! $: namely, $\cV(\cG)\! = \! \cV \! = \! \cV^o\backslash \{U^o_1 \}$ and $\ccE(\cG)\! = \! \big\{ \{ i,j\} \ino \ccE (\cG^o): i,j\ino \cV(\cG)\big\}$.}    
\end{compactenum}

\smallskip

We next consider which connected components of $\cG$ 
are attached to $U^o _1$ in $\cG^o\! $.  
To that end, recall that the $\mathtt{C}^{_{\cG^o}}_{^{l}}$ (resp.~the $\mathtt{C}^{_{\cG}}_{^{l}}$) are the connected components of $\cG^o$ (resp.~of $\cG$) listed in the increasing order of their least vertex; recall that 
$\mathbf{s}^o$ (resp.~$\mathbf{s}$) is the permutation on $\{ 1, \ldots , \bq^o\}$ (resp.~on $\{ 1, \ldots , \bq\}$) such that 
$U^o_k \ino \mathtt{C}^{_{\cG^o}}_{^{\mathbf{s}^o(k)}}$ for all $k\ino \{ 1, \ldots , \bq^o\}$ (resp.~$U_k \ino \mathtt{C}^{_{\cG}}_{^{\mathbf{s}(k)}}$ for all $k\ino \{ 1, \ldots , \bq \}$). 
We first introduce 
\begin{equation}
\label{pieridae}
 \cG^\prime\! :=\!  \cG^o\backslash  \mathtt{C}^{_{\cG^o}}_{^{\mathbf{s}^o(1)}} \quad \textrm{and} \quad 
K \! :=\!  \sup \big \{  k \ino \{ 1, \ldots , \bq\} : T_k \leq a w^o_{U^o_1} \big\} .
\end{equation}
with the convention $\sup \emptyset \! = \! 0$. The graph $\cG^\prime$ is the graph $\cG^o$ where the first (in the order of visit)  connected component has been removed. 
Note that $\cG^\prime$ is possibly empty: it has $\bq^o\! -\! 1$ connected components.   
We easily check that $\bq^o \! = \! \bq \! -\!  K + 1$. 
Here again, we recall that the connected components of $\cG^\prime$ are listed in the increasing order of their least vertex: $\min \mathtt{C}^{\cG^\prime}_1 \! < \! \ldots \! < \! \min \mathtt{C}^{\cG^\prime}_{\mathbf{q}^o-1}$. 
Then, when $\bq^o\! \geq \! 2$, denote by $\mathbf{s}^\prime$ the permutation of $\{ 1, \ldots , \bq^o-1 \}$ such that 
\begin{equation}
\label{papilionidae}
\forall k\ino \{ 1, \ldots , \bq^o\! -\! 1\}, \quad 
\mathtt{C}^{{\cG^\prime}}_{{\mathbf{s}^\prime(k)}}\! = \!  \mathtt{C}^{{\cG^o}}_{{\mathbf{s}^o(k+1)}}\; .
\end{equation}
We also set:   
\begin{equation}
\label{nymphalidae}
\forall  k \ino \{ 1, \ldots , \bq^o\! -\! 1\}, \quad T^\prime_k \! =\! 
 T_{K + k} \! -\! a w^o_{U^o_1}  \quad \textrm{and} \quad U^\prime_{k}\! = \! U_{K +k} \; .
\end{equation}
Suppose that  $\underline{E}^o, \cQ_1, \cQ_2$ are fixed, then we also check that 
\begin{equation}
\label{erebia}
\forall k \ino \{ 1, \ldots , \bq^o\! -\! 1\}, \quad  T^o_{k+1}= T^o_1 + \frac{_{_1}}{^{^a}} T^\prime_k  , \quad U^o_{k+1}\! = \! U^\prime_{k}   \quad 
\textrm{and} \quad \mathtt{C}^{{\cG^\prime}}_{{\mathbf{s}^\prime(k)}}\! = \! \mathtt{C}^{{\cG}}_{
{\mathbf{s}(K+k)}}\; .
\end{equation}

%

 We now explain how additional edges are added to connect $\cG$ to $U^o_1$.  
For all $j\ino \cV$, let $I_j\! = \! \{ t\ino [0, \infty) \! : \! V_t \! = \! j\}$; $I_j$ is 
the set of times during which Client $j$ is served; we easily check that 
 $I_j$ is a finite union of disjoint intervals of the form $[x, y)$ whose total Lebesgue measure 
 is $w_j$: namely, $\ell  (I_j)\! =\! w_j$.
 We also set: 
  $$\Pi_j^* \! =\!  \big\{ y\ino [0, \infty)\! : \! \exists t\ino I_j \; \textrm{such that} \; (t,y)\ino \cQ_2 \, \textrm{and} \; y \! > \! -J_t  \big\} \,.
$$
Note that if $j\ino  \mathtt{C}^{_{\cG}}_{^{\mathbf{s}(k)}}$ and $t\ino I_j$, then $-J_{t}\! = \! T_k$. Combined with elementary results on Poisson random sets, it implies the following.

\smallskip

\begin{compactenum} 
\item[(IV)] \textit{Fix $\underline{E}^o$ and $\cQ_1$; 
suppose that $\cQ_2$ is a Poisson random subset of $[0, \infty)^2$ 
with intensity 
$\frac{_1}{^{\sigma_1 (\bw)}} \, dtdy$. Then, the $(\Pi^*_j)_{j\in \cV}$ are independent and $\Pi^*_j$ is a Poisson random subset of $(T_k, \infty)$ with rate $w_j / \sigma_1 (\bw)$, where $k$ is such that $j\ino \mathtt{C}^{_{\cG}}_{^{\mathbf{s} (k)}}$.  
 }    
\end{compactenum}

\smallskip

\noi
We next introduce the following.

\smallskip

\begin{compactenum}

\smallskip

\item[$\ast$] We set $\Pi_j \! := \! \{ T_k \} \cup \Pi_j^*$ if there is $k\ino \{ 1, \ldots, \bq\}$ such that $j\! = \! U_k$.   

\smallskip

\item[$\ast$]  We set $\Pi_j\! := \! \Pi^*_j$ if $j\ino \cV \backslash \{ U_1, \ldots , U_{\bq}\}$.

\smallskip

\item[$\ast$]  We set 
$\Pi_j^\prime\! := \! \big\{ y \! -\! aw^o_{U^o_1} \, ; \, y \! \in \Pi_j \! \cap \! ( aw^o_{U^o_1} , \infty) \big\}$. 

\smallskip

\item[$\ast$] For all non-empty $S \! \subset \! \cV$, we set $\Pi^\prime (S)\! := \! \bigcup_{j\in S}\Pi^\prime_j $. 

\smallskip

\end{compactenum}

\noi
We claim the following.

\smallskip

\begin{compactenum} 
\item[(V)] \textit{Fix $\underline{E}^o, \cQ_1, \cQ_2$. Then,} 

\smallskip

\begin{compactenum} 
\item[(Va)] \textit{for all $j\ino \cV$, $\{ U^o_1, j \} \ino \ccE (\cG^o)$ iff $\, \Pi_j \! \cap\!  [0, aw^o_{U^o_1}] \! \neq \! \emptyset$;}

\smallskip

\item[(Vb)] \textit{for all $k\ino \{ 1, \ldots , \bq^o\! -\! 1\}$, $T_{K+k}\! -\! aw^o_{U^o_1} \! =\! \inf \Pi^\prime ( \mathtt{C}^{{\cG}}_{{\mathbf{s}(K+k)}})\! = \! \inf~\Pi^\prime_{U_{K+k}}$; namely, 
\begin{equation}
\label{lycaenidae}
\forall k\ino \{ 1, \ldots , \bq^o\! -\! 1\}, \quad  T^\prime_k \! = \! \inf \Pi^\prime 
( \mathtt{C}^{{\cG^\prime}}_{{\mathbf{s}^\prime(k)}})=  \inf~ \Pi^\prime_{U^\prime_k} , 
\end{equation}
with the notation of \eqref{pieridae}, \eqref{papilionidae}, \eqref{nymphalidae} and \eqref{erebia}.}
\end{compactenum}
\end{compactenum}

\smallskip

\noi
\textit{Proof of (V).} The argument is deterministic and elementary. Suppose that $\{ U^{o}_{1}, j\} \ino \ccE(\cG^o)$. There are two cases to consider.  

\noi
$\bullet$ \texttt{(Case 1)}: $\{ U^o_1, j\}$ is part of the exploration tree generated by $(\underline{E}^o, \bw^o)$ which means that $V^o_{E^o_j-}\! = \! U^o_1$. Namely, it means that 
$T^o_1\! = \! \sup \{ s\ino [0, E^o_j]: Y^o_{s-} \! < \! \inf_{[s, E^o_j]} Y^o\}$ 
because $E^o_{U^o_1}\! = \! T^o_1$, by definition of $(U^o_1, T^o_1)$; 
since $Y^o_{T^o_1-}+ w^{o}_{U^o_1}\! = \! Y^o_{T^o_1}$ and by (\ref{kermes}), 
it is equivalent to: $-aw^{o}_{U^o_1}<J_{E_j} \! = \! Y_{E_j-}$. It implies that $E_j$ is the left endpoint 
of an excursion of $Y$ strictly above its infimum; therefore, there exists $k\ino \{ 1, \ldots , \bq \} $ such that $j\! = \! U_k$ and  $-J_{E_j} \! = \! T_k$, by (\ref{syllepse}); thus $T_k \ino [0, aw^o_{U^o_1}]$ which implies that $ \Pi_j \! \cap \! [0, aw^o_{U^o_1}] \! \neq \! \emptyset$ (since $T_k \ino \Pi_j$ as we are in the case where $j\! = \! U_k$).   

Conversely, let  $j\! = \! U_k$ be such that  $\Pi_j \cap [0, aw^o_{U^o_1}] \! \neq \! \emptyset$. It implies that $T_k \! \leq \! aw^o_{U^o_1}$ and the previous arguments can be reversed verbatim to prove that $\{ U^o_1, j\}$ is an edge of $\cG^o$ that is part of the exploration tree generated by $(\underline{E}^o, \bw^o)$. 

\smallskip

\noi
$\bullet$ \texttt{(Case 2)}: $\{ U^o_1, j\}$ is an additional edge of $\cG^o\! $. Then, there exists 
$(t^\prime, y^\prime)\ino \Pi^o$ such that $V^{o}_{t^\prime}\! = \! j$ and 
$V^o_{\tau^o (t^\prime, y^\prime)}\! = \! U^o_1$, where $\tau^o (t^\prime,y^\prime) \! =\!  
\inf \big\{ s\ino [0, t^\prime] : \inf_{[s, t^\prime]} Y^o \!  >\!  y^\prime+ J^o_{t^\prime} \big\}$. Note 
that $V^{o}_{t^\prime}\! = \! j$ implies that $t^\prime \! >\! T^o_1$ since $j\! \neq \! U^o_1$ and since $U^o_1$ is the first visited vertex (or the first client). Also observe that $V^o_{\tau^o (t^\prime, y^\prime)}\! = \! U^o_1$ implies $\tau^o (t^\prime,y^\prime)\! =\! T^o_1$. It also implies 
that $t^\prime$ lies in the first excursion interval of $Y^o$ strictly above its infimum, which entails that $J^o_{t^\prime}\! = \! J^o_{T^o_1}\! = \! -T^o_1$. 
Then, we set $t\! = \! a(t^\prime \! -\! T^o_1)$ and, thanks to (\ref{kermes}), we rewrite the previous conditions in terms of $Y, J$ and $V$ as follows: $V_{t}\! = \! j$ and 
$0\!  =\!  \inf \{ s\ino [0, t]\! : \! aw^o_{U^o_1}+ \inf_{[s, t]} Y \! >\! ay^\prime \}$, which is equivalent to: 
$t\ino I_j $ and $y\! :=\! a( w^o_{U^o_1} \! -\! y^\prime) \! > \!  -J_{t}$. 
This proves that there is $(t, y)\ino D_2$ (as defined in (\ref{+Pi12})) such that $(t^\prime, y^\prime)\! = \! f_2 (t,y)$ as defined in (\ref{+Pio12}). Since $f_1(\Pi)$ and $f_2 (\Pi^2)$ form 
a partition of $\Pi^o$, $(t,y) \ino \Pi^2$ and this proves that  
there is $(t, y)\ino \cQ_2$ such that $t\ino I_j $ and 
$aw^o_{U^o_1} \! \geq \! y >\!    -J_t$ 
which implies $\Pi_j \cap [0, aw^o_{U^o_1}] \! \neq \! \emptyset$. 

Conversely, suppose that $\Pi_j  \cap  [0, aw^o_{U^o_1}] \! \neq \! \emptyset$ and that $j\ino \cV \backslash \{ U_1, \ldots, U_{\bq}\}$. Then, 
$\Pi^*_j\!  \cap \! [0, aw^o_{U^o_1}] \! \neq \! \emptyset$ and the previous arguments can be reversed 
verbatim to prove that $\{ U^o_1, j\}$ is an (additional) edge of $\cG^o$, which 
completes the proof of (Va).  

\medskip

Let us prove (Vb). Let $k\ino \{ 1, \ldots, \bq^o\! -\! 1\}$. By definition $U_{K+k}\ino \mathtt{C}^{_{\cG}}_{^{\mathbf{s}(K+k)}}$. Let $y\ino 
\Pi^*_{U_{K+k}}$: namely, there exists $t$ such that $V_t\! = \! U_{K+k}$, $(t,y)\ino \cQ_2$ and $y \! > \! -J_t$. 
But $V_t\! = \! U_{K+k}$ implies that $-J_t\! = \! T_{K+k}$. Since $\Pi_{U_{K+k}} \! = \! \Pi_{U_{K+k}}^* \! \cup \! \{ T_{K+k}\} $, we get $\inf \Pi_{U_{K+k}} \! = \! T_{K+k}$. 
By definition of $K$, $T_{K+k}\! >\! aw^o_{U^o_1}$, which entails $\inf \Pi^\prime_{U_{K+k}} \! = \! T_{K+k} \! -\! aw^o_{U^o_1}$.

Let $j\ino  \mathtt{C}^{_{\cG}}_{^{\mathbf{s}(K+k)}}\backslash \{U_{K+k} \}$ (if any) and let 
$y\ino \Pi_j$. Necessarily, $j\ino \cV \backslash \{U_1, \ldots, U_{\bq}\}$, 
which entails $\Pi_j\! = \! \Pi^*_j$, by definition. Then there exists $t$ such that $V_t\! = \! j$, $(t,y)\ino \cQ_2$ and $y \! > \! -J_t$. 
Note that $I_{j}$ is included in the excursion interval of $Y$ strictly above its infimum whose left endpoint is 
$E_{U_{K+k}}$, which implies that $J_t \! = \! -T_{K+k}$, for all $t\ino I_j$. Thus $y \! >\! T_{K+k}$. This proves that $\inf \Pi^\prime ( \mathtt{C}^{_{\cG}}_{^{\mathbf{s}(K+k)}}) \! = \! \inf \Pi^\prime_{U_{K+k}} $, which completes the proof of (Vb). \cq 
 
\medskip

We now complete the proof of Theorem \ref{+th} as follows: let $(E^o_j)_{j\in \cV^o}$ be 
independent  exponentially distributed r.v.\ such that $\bE [E^o_j] \! = \!  \sigma_1 (\bw^o)/w^o_j$ for all $j\ino \cV^o$. We then fix $j^* \ino \cV^o$ and we work under $\bP_{\! j^*}\! :=\! \bP ( \, \cdot \, | U^o_1 \! = \! j^*)$; under $\bP_{\! j^*}$, 
by (\ref{truchement}), we get  
$\cV\! = \! \cV^o\backslash \{ j^*\}$, $a\! =\!  (\sigma_1 (\bw^o)\! -\! w^o_{j^*})/ \sigma_1 (\bw^o)$ and for all $j\ino \cV$, $w_j\! = \! aw^o_j$ and $E_j\! = \! a(E^o_j\! -\! T^o_1)$. Under $\bP_{\! j^*}$, we take 
$\cQ_1$ and $\cQ_2$ as two independent Poisson random subsets of $[0, \infty)^2$ with intensity $\frac{_1}{^{\sigma_1 (\bw)}} \, dtdy$; $\cQ_1$ and $\cQ_2$ are also supposed 
independent of $\underline{E}^o$. 
Recall from (\ref{+Pi12})  and (\ref{+Pio12}) the definitions of $\Pi$ and 
$\Pi^o$ and from (\ref{meloe}) the definition $(\cG, (T_k, U_k)_{1\leq k\leq \bq})$. By (I), under $\bP_{\! j^*}$ $(\underline{E}, \bw, \Pi)$ has the required distribution so that the 
induction hypothesis applies to 
$(\cG, (T_k, U_k)_{1\leq k\leq \bq})$: namely, under $\bP_{\! j^*}$, $\cG$ has law $M_{\cV\! , \, \bw}$ (as defined in (\ref{gongorisme})) 
and conditionally given $\cG$, $(T_k, U_k)_{1\leq k\leq \bq}$ 
has law $\Lambda_{\cG\! , \, \rbm}$ as defined in 
(\ref{Lambformu}) (namely, (\ref{+loi}) holds).  

By (IV), under $\bP_{\! j^*}$ and conditionally given $(T^o_1, \cG, (T_k, U_k)_{1\leq k\leq \bq})$, the 
$(\Pi^*_j)_{j\in \cV}$ are independent and 
$\Pi^*_j$ is a Poisson random subset of $(T_k, \infty)$ with rate $w_j / \sigma_1 (\bw)$, where $k$ is such that $j\ino \mathtt{C}^{_{\cG}}_{^{\mathbf{s} (k)}}$. Since, conditionnally given 
$(T^o_1, \cG)$, the r.v.~$(T_k, U_k)_{1\leq k\leq \bq}$ has law $\Lambda_{\cG\! , \, \rbm}$, 
the definition of $\Pi_j$ combined with elementary results on Poisson processes imply the following key point.

\smallskip

\begin{compactenum} 
\item[(VI)] \textit{Under $\bP_{\! j^*}$ and conditionally given $(T^o_1, \cG)$, the 
$(\Pi_j)_{j\in \cV}$ are independent and $\Pi_j$ is a Poisson random subset of $(0, \infty)$ with rate $w_j / \sigma_1 (\bw)$; therefore, under $\bP_{\! j^*}$, the $(\Pi_j)_{j\in \cV}$ 
are independent of $(T^o_1, \cG)$ and the very definition of $\Pi_j$ implies that for all 
$\{ 1, \ldots , \bq \}$,} 
$$
T_k \! = \! \inf \Pi (\mathtt{C^{\cG}_{\mathbf{s} (k)}} \big) \! = \! \inf \Pi_{U_k} \; \textit{where $\mathbf{s}$ is such that} \;  \inf \Pi \big(\mathtt{C}^{_{\cG}}_{^{\bs (1)}} \big) \! <\!  \ldots  \! < \!  \inf \Pi \big( \mathtt{C}^{_{\cG}}_{^{\bs (q)}} \big).$$
\end{compactenum}

Consequently, under $\bP_{\! j^*}$, 
by (Va) and (VI), $\cG^o$ only depends on $\cG$ and on the 
$\Pi_j \cap [0, aw^o_{j^*}]$, $j\ino \cV$. This, combined with
elementary results on Poisson processes, implies that  
$\cG^o$ is independent from $(\Pi^\prime_{j})_{j\in \cV}$ and from $T^o_1$; (Va) and (VI) also  
imply that under $\bP_{\! j^*}$, the events $\{ \{ j^*\! , j\} \ino \ccE(\cG^o)\}$, $j\ino \cV$, are independent 
and occur with probabilities $1\! -\! \exp (-w_jaw^o_{j^*}/ \sigma_1 (\bw))$. Then, 
note that $w_jaw^o_{j^*}/ \sigma_1 (\bw)\! = \! w^ o_j w^o_{j^*}/ \sigma_1 (\bw^o)$. Thus, under $\bP_{\! j^*}$, $\cG^o$ has law $M_{\cV^o\! , \, \bw^o}$ and it is independent from 
$(T^o_1; \Pi^\prime_j , j\ino \cV)$.  

Recall from (\ref{pieridae}), (\ref{papilionidae}), (\ref{nymphalidae}) and (\ref{erebia}) notation 
$\cG^\prime$, $\mathbf{s}^\prime$ and $(T^\prime_k, U^\prime_k)_{1\leq k \leq \bq^o-1}$. 
Then, observe that under $\bP_{\! j^*}$ and conditionally given $\cG^o$, (Vb) and (VI) imply that 
$(T^\prime_k, U^\prime_k)_{1\leq k \leq \bq^o-1}$ has conditional law $\Lambda_{\cG^\prime\! , \, \rbm^\prime}$, where $\rbm^\prime\! = \! \sum_{j\in \cV(\cG^\prime)} w_j \delta_j$. Then, under $\bP_{\! j^*}$ and conditionally given $\cG^o$, Lemma \ref{1compopo} applies and (\ref{erebia}) entails that 
\begin{eqnarray*}
 \frac{w^o_{j^*}}{\sigma_1 (\bw^o)} \bP_{j^*} 
 \big( T^o_1 \ino dt_1; \ldots ; T^o_{\bq^o} \ino dt_{\bq^o} \!\!\!\!\!\! &; & \!\!\!\!\!\!  U^o_2\! =\! j_2; \ldots ; U^o_{\bq^o}\! = \! j_{\bq^o} \, \big| \, \cG^o \big) \\ 
 & =  & 
\Lambda_{\cG^o\! ,\,  \rbm^o} (dt_1 \ldots dt_{\bq^o}; j^*\! , j_2, \ldots , j_{\bq^o} ) .  
\end{eqnarray*}
Since $\bP (U^o_1\! = \! j^*)\! = \! w^o_{j^*}/ \sigma_1 (\bw^o)$, it implies that for all graph 
$G^o$ whose set of vertices is $\cV^o$ and that has $q^o$ connected components, we get 
\begin{eqnarray*}
 \bP\big( \cG^o\! =\!  G^o \, ; \, T^o_1 \ino dt_1; \ldots ; T^o_{q^o} \ino dt_{q^o} \!\!\!\!\!\! &; & \!\!\!\!\!\! U^o_1 \! = \! j^*;  U^o_2\! =\! j_2; \ldots ; U^o_{q^o}\! = \! j_{q^o} \big) \\ 
 & =  & 
M_{\cV^o\! , \, \bw^o} (G^o) \, \Lambda_{G^o\! , \, \rbm^o} (dt_1 \ldots dt_{q^o}; j^*\! , j_2, \ldots , j_{q^o} ) .  
\end{eqnarray*}
This completes the proof of Theorem \ref{+th} by induction on the number of vertices. \cqfd

\section{Embedding the multiplicative graph in a GW-tree}
\label{EmbedGW}
We recall from Section \ref{XLIFOsec} the definition of the Markovian algebraic load process 
$X^\bw$ that is derived from the blue and red processes $X^{\mathtt{b}, \bw}$ and $X^{\mathtt{r}, \bw}$ as in (\ref{redblumix}). Namely recall from (\ref{ccXbrwdef}) that 
$\ccX_\bw^{\mathtt{b}}\! = \! \sum_{k\geq 1} \delta_{(\tau^{\mathtt{b}}_k , \Jtt^{\mathtt{b}}_k)}$ and 
$\ccX_\bw^{\mathtt{r}}\! = \! \sum_{k\geq 1} \delta_{(\tau^{\mathtt{r}}_k , \Jtt^{\mathtt{r}}_k)}$ are two  independent Poisson random point measures on $[0, \infty)  \times  \{ 1, \ldots , n\} $ with intensity $\ell  \otimes  \nu_\bw$, where $\ell $ stands for the Lebesgue measure on $[0, \infty)$ and where $\nu_\bw \! = \! \frac{1}{\sigma_1 (\bw)}\sum_{1 \leq j\leq n} w_j \delta_{j} $. Recall from (\ref{Xbrwdef}) the following notations 
$$X^{\mathtt{b}, \bw}_t  =  -t + \sum_{k\geq 1} w_{\Jtt^{\mathtt{b}}_k}\un_{[0, t]} (\tau^{\mathtt{b}}_k), \quad   X^{\mathtt{r}, \bw}_t  =  -t + \sum_{k\geq 1} w_{\Jtt^{\mathtt{r}}_k}\un_{[0, t]} (\tau^{\mathtt{r}}_k)  \quad \textrm{and} \quad I^{\mathtt{r}, \bw}_t = \inf_{s\in [0, t]} X^{\mathtt{r}, \bw}_s\; .$$
 For all $j\ino \{ 1, \ldots , n\}$ and all $t\ino [0, \infty)$, recall from (\ref{Njbdef}) that  
$ N^{\bw}_j (t)\! = \! \ccX_\bw^{\mathtt{b}} \big( [0, t] \! \times \! \{ j\} \big)$ and that 
$E^\bw_j \! = \! \inf \big\{ t\ino [0, \infty) : \ccX^{\mathtt{b}}_\bw ([0, t] \! \times \! \{ j\})\! = \! 1 \big\}$; 
the $N^\bw_j$ are independent homogeneous Poisson processes with jump-rate $w_j/ \sigma_1 (\bw)$ and the r.v.~$(\frac{w_j}{\sigma_1 (\bw)} E^\bw_j)_{1\leq j\leq n}$ are i.i.d.~exponentially distributed r.v.~with unit mean. Note that $X^{\mathtt{b}, \bw}_t \! =\! -t + \sum_{1\leq j\leq n} w_j N^{\bw}_j (t)$. Then recall from (\ref{YwAwSigw}) that 
\begin{equation}
\label{YAglopop}
Y^\bw_t \! = \! -t  \, +\!\!\!   \sum_{1\leq j \leq n} \!\!\!  w_j \un_{\{ E^\bw_j \leq t \}}
 \quad \textrm{and} \quad A^\bw_t =   X^{\mathtt{b}, \bw}_t -Y^\bw_t  =  
  \! \sum_{1\leq j\leq n} \!\! w_j (N^\bw_j (t)\! -\! 1)_+ \; .
 \end{equation} 
Thanks to $(Y^\bw, A^\bw)$ and $X^{\mathtt{r}, \bw}$ we reconstruct the Markovian LIFO queue as follows. First recall from (\ref{thetabw}) the definition of the "blue" time-change 
$ \theta^{\mathtt{b}, \bw}$: for all $t\ino [0, \infty)$,  
\begin{equation}
\label{rethetabw}
 \theta^{\mathtt{b}, \bw}_t \!\! = t + \subo^{\mathtt{r}, \bw}_{A^\bw_t}, \;   \textrm{where for all $x\ino [0, \infty)$, we have set:} \;\;  \subo^{\mathtt{r}, \bw}_x \!  = \! \inf \big\{ t\ino [0, \infty ) \! :  X^{\mathtt{r} , \bw}_t \!\!\!  < \! - x \big\}, 
\end{equation}
with the convention that $\inf \emptyset \! = \! \infty$. Note that $\subo^{\mathtt{r}, \bw}_x\! < \! \infty$ iff $x\! < \! \! -I^{\mathtt{r}, \bw}_\infty\! = \! \lim_{t\rightarrow \infty} \! - I^{\mathtt{r}, \bw}_t$ that is a.s.~finite in supercritical cases (and $ -I^{\mathtt{r}, \bw}_\infty$ is a.s.~infinite in critical and subcritical cases).    
%
We also need to recall from (\ref{T*wdef}) the definition of 
\begin{equation}
\label{grumide}
T^*_\bw\! = \!  \sup \{ t\ino [0, \infty)\! : A^\bw_t \! < \! - I^{\mathtt{r}, \bw}_\infty \} \!  = \! \sup \{ t\ino [0, \infty):  \theta^{\mathtt{b}, \bw}_t \! < \infty \}.
\end{equation} 
In critical and subcritical cases, $T^*_\bw\! = \! \infty$ and $ \theta^{\mathtt{b}, \bw}$ only takes 
finite values. In supercritical cases, a.s.~$T^*_\bw \! < \! \infty$ and we check that $ \theta^{\mathtt{b}, \bw} (T^*_\bw-) \! < \! \infty$. 
We next recall from (\ref{blureddef}) the definition of $\mathtt{Blue}$ and 
$\mathtt{Red}$ that are the sets of times during which blue and red clients are served (recall that the server is considered as a blue client): 
\begin{equation}
\label{gramidos}
\mathtt{Red}\! = \! \bigcup_{t\in [0, T^*_\bw]: \Delta \theta^{\mathtt{b}, \bw}_t >0} \big[  \theta^{\mathtt{b}, \bw}_{t-}, 
\theta^{\mathtt{b}, \bw}_t \big) \quad \textrm{and} \quad \mathtt{Blue} \! = \! [0, \infty) \backslash \mathtt{Red} . 
\end{equation}
Note that the union defining $\mathtt{Red}$ is countably infinite in critical and subcritical cases and that it is a finite union  
in supercritical cases since $\big[  \theta^{\mathtt{b}, \bw} (T^*_\bw-),  \theta^{\mathtt{b}, \bw} (T^*_\bw))\! = \! [\theta^{\mathtt{b}, \bw} (T^*_\bw-), \infty)$. 
We next recall from (\ref{Lambdefi}) the definition of the time-changes $\Lambda^{\mathtt{b}, \bw}$ and $\Lambda^{\mathtt{r}, \bw}$: 
\begin{equation}
\label{gagaoutch}
\Lambda^{\mathtt{b}, \bw}_t \!  = \! \int_0^t \!\! \un_{\mathtt{Blue}} (s) \, ds =  \inf \big\{ s\ino [0, \infty) \! :  \theta^{\mathtt{b}, \bw}_s \! >\! t \}  \quad \textrm{and} \quad 
 \Lambda^{\mathtt{r}, \bw}_t \!\! = t \! -\!  \Lambda^{\mathtt{b}, \bw}_t \! \! = \! \! \int_0^t \!\! \un_{\mathtt{Red}} (s) \, ds. 
 \end{equation}
The processes $\Lambda^{\mathtt{b}, \bw}$ and $\Lambda^{\mathtt{r}, \bw}$ are continuous and nondecreasing and a.s.~$\lim_{t\rightarrow \infty}  \Lambda^{\mathtt{r}, \bw}_t \! = \! \infty$. 
In critical and subcritical cases, we also get a.s.~$\lim_{t\rightarrow \infty}  \Lambda^{\mathtt{b}, \bw}_t \! = \! \infty$ and $\Lambda^{\mathtt{b}, \bw} (\theta^{\mathtt{b}, \bw}_t)\! = \! t$ for all $t\ino [0, \infty)$.  However, in supercritical cases, $ \Lambda^{\mathtt{b}, \bw}$ is eventually constant and equal to 
$T^*_\bw$ and a.s.~for all $t\ino [0, T^*_\bw)$, 
$\Lambda^{\mathtt{b}, \bw} (\theta^{\mathtt{b}, \bw}_t)\! = \! t$. 
We next recall from (\ref{redblumix}) the definition of the 
load of the Markovian queue $X^\bw$:  
\begin{equation}
\label{redblumixbis}
\forall t\ino [0, \infty) , \qquad X^{\bw}_t =  X^{\mathtt{b}, \bw}_{ \Lambda^{\mathtt{b}, \bw}_t } + X^{\mathtt{r}, \bw}_{ \Lambda^{\mathtt{r}, \bw}_t } \; .
\end{equation}

\noi
\textbf{Proof of (\ref{YXtheta}) in Proposition \ref{Xwfrombr}}.
Note that for all 
$t\ino  [0, T_\bw^*)$, $\Lambda^{\mathtt{b}, \bw} (\theta^{\mathtt{b}, \bw}_t)\! = \! t$ and thus 
$\Lambda^{\mathtt{r}, \bw} (\theta^{\mathtt{b}, \bw}_t)\!  = \!  \theta^{\mathtt{b}, \bw}_t \! -\! t\! = \! \gamma^{\mathtt{r}, \bw} (A^\bw_t)$. Therefore, 
$$ X^{\bw}_{\theta^{\mathtt{b}, \bw}_t}\! = \! X^{\mathtt{b}, \bw}_{\Lambda^{\mathtt{b}, \bw} (\theta^{\mathtt{b}, \bw}_t)}+  X^{\mathtt{r}, \bw}_{\Lambda^{\mathtt{r}, \bw} (\theta^{\mathtt{b}, \bw}_t)}= 
 X^{\mathtt{b}, \bw}_t + X^{\mathtt{r}, \bw}_{\gamma^{\mathtt{r}, \bw} (A^\bw_t)}=  X^{\mathtt{b}, \bw}_t -A_t^\bw\! = \! Y^\bw_t, $$
which proves  (\ref{YXtheta}) in Proposition \ref{Xwfrombr}. \cqfd

\medskip

We next prove  Lemma \ref{Hthetalem}. To that end we need the following lemma. 
\begin{lem}
\label{gumodis} Almost surely, for all $b\ino [0, T^*_\bw]$ such that $\theta^{\mathtt{b}, \bw}_{b-}<\theta^{\mathtt{b}, \bw}_{b}$, we get for all $s\ino [\theta^{\mathtt{b}, \bw}_{b-} , \theta^{\mathtt{b}, \bw}_{b} )$
\begin{equation}
\label{flaccidos}
 X^{\bw}_s > X^\bw_{(\theta^{\mathtt{b}, \bw}_{b-})-}\! = \! Y^{\bw}_b , \quad \Delta X^\bw_{\theta^{\mathtt{b}, \bw}_{b-}}= \Delta A^{\bw}_b \quad \textrm{and} \quad X^\bw_{(\theta^{\mathtt{b}, \bw}_{b-})-}\! \! = \! X^\bw_{\theta^{\mathtt{b}, \bw}_{b}}=X^\bw_{(\theta^{\mathtt{b}, \bw}_{b})-} \quad \textrm{if} \quad \theta^{\mathtt{b}, \bw}_{b}\!  < \infty. 
\end{equation}
Thus, a.s.~for all $s\ino [0, \infty)$, $X^\bw_s \! \geq \! Y^\bw(\Lambda^{\mathtt{b}, \bw}_s)$. Moreover, a.s.~for all $s_1, s_2 \ino [0, \infty)$ such that $\Lambda^{\mathtt{b}, \bw}_{s_1} \! < \!  \Lambda^{\mathtt{b}, \bw}_{s_2}$, then 
\begin{equation}
\label{flaccid}
\inf_{b\in [\Lambda^{\mathtt{b}, \bw}_{s_1} ,  \Lambda^{\mathtt{b}, \bw}_{s_2}]} Y^\bw_b= \inf_{s\in [s_1, s_2]} X^\bw_s\; .
\end{equation}
We next introduce the red time-change: 
\begin{equation}
\label{thetarw}
\theta^{\mathtt{r}, \bw}_t = \inf \big\{s \ino [0, \infty): \Lambda^{\mathtt{r}, \bw}_s >t \,  \big\} \,.
\end{equation}
Then, for all $s, t\ino [0, \infty)$, $\theta^{\mathtt{r}, \bw}_{s+t} \! -\! \theta^{\mathtt{r}, \bw}_t \! \geq \! s$ and if $\Delta \theta^{\mathtt{r}, \bw}_t \! >\! 0$, then $\Delta 
X^{\mathtt{r}, \bw}_t \! = \! 0$. 
\end{lem}
\noi

\begin{figure}[tb]
\centering
\includegraphics[scale=.7]{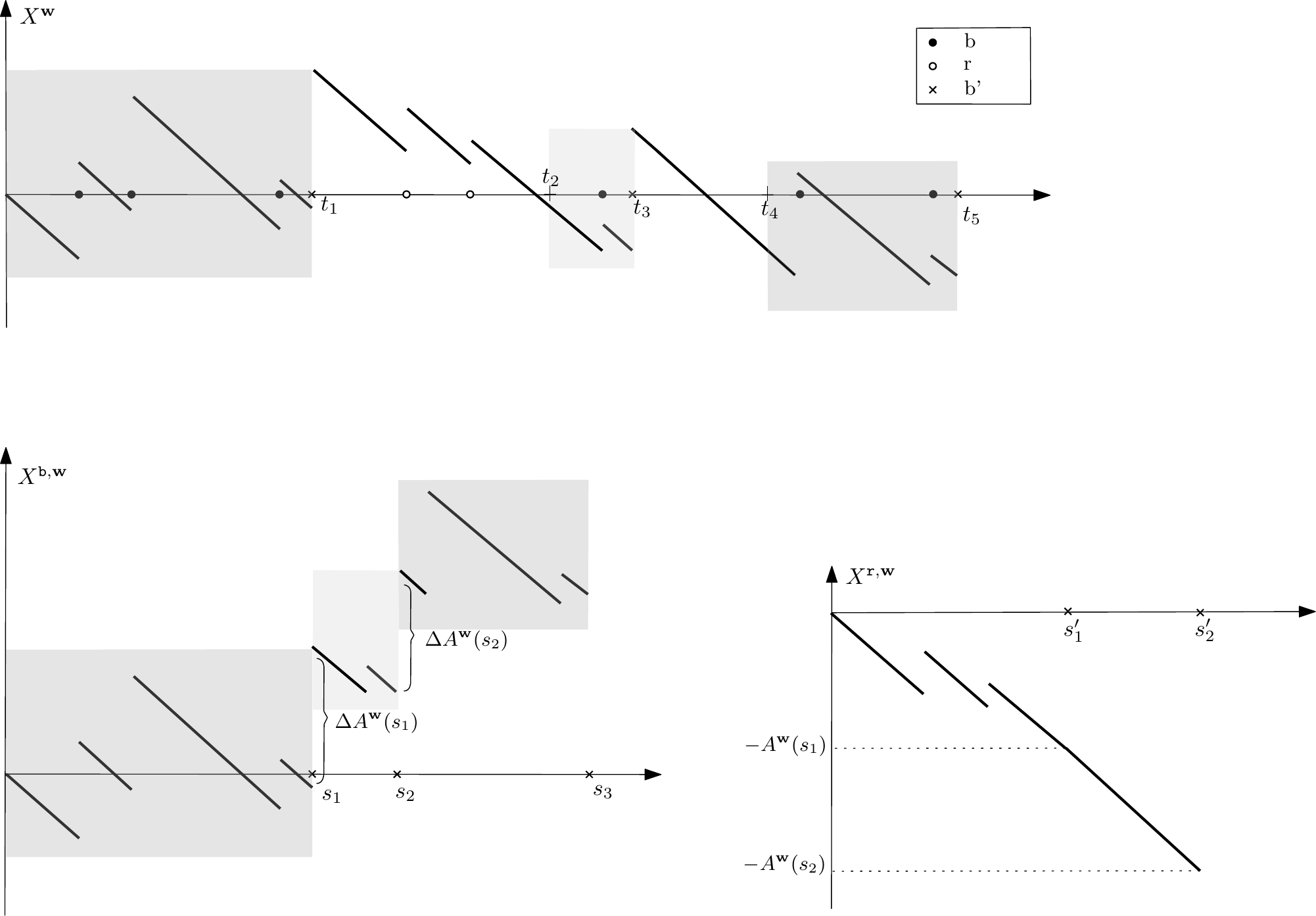}
\caption{{\small Colouring of the Markovian queue.  The grey blocks form the set $\mathtt{Blue}$. Observe that the sequence $(t_l)_{l\in \bbN}$ (i.e.~the moments when the queue switches its colours), marks the boundaries of the blocks. 
Recall the point measure $\overline{\ccX}_\bw$ formed by the jumps of $X^{\bw}$. Jumps of colour $\mathtt{b}$ are marked by $\bullet$ on the abscissa; those of colour $\mathtt{b}'$ marked by ${\scriptstyle \times}$ and those of colour $\mathtt{r}$ by $\circ$\,. Observe also that jumps inside the grey blocks are of colour $\mathtt{b}$, the right endpoint of each grey block corresponds to a jump of colour $\mathtt{b}'$, which is also a jump of $A^{\bw}$. Jumps outside these blocks are of colour $\mathtt{r}$. \cq}  }
\label{fig:colouring} 
\end{figure}

\textbf{Proof.}
 As defined by (\ref{redblumixbis}), the process $X^\bw$ is obtained by concatenating alternatively successive parts of 
$X^{\mathtt{b}, \bw}$ and $X^{\mathtt{r}, \bw}$ as follows: recall from (\ref{YAglopop}) that $A^\bw$ is the sum of the jumps of $X^{\mathtt{b}, \bw}$ that already 
occurred once; it is therefore a nondecreasing process; denote by $s_l$ the $l$-th jump-time 
of $A^\bw$ and to simplify, set $s^\prime_l \! = \! \gamma^{\mathtt{r}, \bw}(A^\bw_{s_l})$ (note that in supercritical cases, $s^\prime_l$ may be infinite). 
Let us agree on saying that the colour of the queue is defined by the colour of the client currently served. We then introduce a sequence $(t_{l})_{l\in \bbN}$ which corresponds to the times when the queue switches its colour. For convenience, we first set $s_0\! = \! s^\prime_0\! = \! 0$ and then we introduce
\begin{equation}
\label{telltell} 
\forall l\ino \bbN, t_{2l}\! := \! s_l+s^\prime_l \quad \textrm{and} \quad t_{2l+1}\! := \! s_{l+1}+ s^\prime_l. \quad 
\end{equation} 
See also Fig.~\ref{fig:colouring}.
Next set $l^*\! = \! \inf \big\{ l \! \geq \! 0: s^\prime_l \! =\!  \infty \big\}$ (that is infinite in the subcritical and the critical cases). 
Then, $t_l \! < \! \infty$ for all $l \! < \! 2l^*$ and $t_{2l^*}\! = \! \infty$. 
Fix $l\! < \! l^*$. Then, note that $t_{2l+1}\! -\! t_{2l}\! = \! s_{l+1} \! -\! s_l$, that 
$t_{2l+2}\! -\! t_{2l+1}\! = \! s^\prime_{l+1} \! -\! s^\prime_l$ and observe that (\ref{redblumixbis}) 
can be rewritten as follows. 
\begin{multline}
\label{mourmansk}
 \forall s\ino [0, t_{2l+1} \! -\! t_{2l} ) , \quad 
X^\bw_{t_{2l}+s}\! -\!  X^\bw_{t_{2l}}\! = \! X^{\mathtt{b}, \bw}_{s_l+s}\! -\! X^{\mathtt{b}, \bw}_{s_l} \; ,  \\
\textrm{and} \quad  \forall s\ino [0, t_{2l+2} \! -\! t_{2l+1} ) , \quad    X^\bw_{t_{2l+1}+s}\! -\!  X^\bw_{t_{2l+1}}\! = \! X^{\mathtt{r}, \bw}_{s^\prime_l+s}\! -\! X^{\mathtt{r}, \bw}_{s^\prime_l} . 
\end{multline}       
Moreover, we get $\mathtt{Red}\! = \! \bigcup_{0\leq l<l^*} [t_{2l+1} , t_{2l+2} )$ and $ \mathtt{Blue}\! = \!  \bigcup_{0\leq l<l^*} [t_{2l} , t_{2l+1} ) $.  
Namely, 
\begin{multline*}
\big\{ \theta^{\mathtt{b} , \bw}_{t-} \, ; \, t \ino [0, T^*_\bw] \! : \Delta \theta^{\mathtt{b}, \bw}_t \! >\! 0 \big\} \! = \l \big\{  t_{2l+1}\,  ; \, 1\! \leq \! l\! <\! l^* \big\} \\ 
\quad \textrm{and} \quad   
\big\{ \theta^{\mathtt{b} , \bw}_{t} \, ;\,  t \ino [0, T^*_\bw] \! : \Delta \theta^{\mathtt{b}, \bw}_t \! >\! 0 \big\} \! = \l \big\{  t_{2l} \, ; 1\! \leq \! l\! <\! l^* \big\} \; . 
\end{multline*}  
Next fix $l\! < \! l^*$. Note that $X^{\mathtt{r}, \bw} (s^\prime_l)\! = \! -A^{\bw} (s_l) \! =\!  -A^\bw (s_{l+1}-)$ and that
$\Delta X^\bw (t_{2l+1}) \! =\! \Delta A^\bw (s_{l+1})\! =\!  \Delta X^{\mathtt{b}, \bw} (s_{l+1})$. 
Thus, for all $s\ino [0, t_{2l+2} \! -\! t_{2l+1} )$,  $X^{\bw} (t_{2l+1}+s) \! - \! X^\bw (t_{2l+1}-)\! = \! X^{\mathtt{r}, \bw} (s^\prime_l+s)\! +\! A^\bw (s_{l+1}) \! >\! 0$. 
Namely, for all $t \ino [t_{2l+1} , t_{2l+2})$, we get $X^\bw (t) \! > \! X^\bw (t_{2l+1} -)$ and if $t_{2l+2}\! < \! \infty$, then we also get $X^{\bw} (t_{2l+2}) \! = \! X^\bw (t_{2l+1}-) $,  
which completes the proof of (\ref{flaccidos}).

Clearly, (\ref{flaccidos}) implies that  a.s.~for all $s\ino [0, \infty)$, $X^\bw_s \! \geq \! Y^\bw(\Lambda^{\mathtt{b}, \bw}_s)$. We next prove (\ref{flaccid}): let  $s_1, s_2 \ino [0, \infty)$ such that $b_1\! :=\! \Lambda^{\mathtt{b}, \bw}_{s_1} \! < \!  \Lambda^{\mathtt{b}, \bw}_{s_2}\! =  : b_2$, which implies that $b_1 \! < \! T^*_{\bw}\! $. We have proved that $Y^\bw_b \! = \! X^\bw(\theta^{\mathtt{b}, \bw}_b)$ for all $b \ino [0, T_\bw^*)$; since by construction, $\Delta Y^\bw (T^*_{\bw})\! = \! 0$, we also get $Y^\bw_{T^*_{\bw}}\! = \! X^\bw (\theta^{\mathtt{b}, \bw}_{T^*_{\bw}\! -}\! -)$. This implies that 
$\inf_{[b_1, b_2]} Y^\bw\! \geq \! \inf_{[\theta^{\mathtt{b}, \bw} (b_1), \theta^{\mathtt{b}; \bw} (b_2-)] }
X^\bw$. Since $s_1 \! \leq\!  \theta^{\mathtt{b}, \bw}_{b_1} \! < \! 
 \theta^{\mathtt{b}, \bw}_{b_2-} \! \leq \! s_2$, we get $\inf_{[b_1, b_2]} Y^{\bw}\! \geq \! 
 \inf_{[s_1, s_2]} X^\bw$. 
But we have proved  also that 
$X^\bw_s \! \geq \! Y^\bw (\Lambda^{\! \mathtt{b}, \bw}_s)$ for all $s\ino [0, \infty)$: therefore, $ \inf_{[s_1, s_2]} X^\bw \! \geq \! \inf_{[b_1, b_2]} Y^\bw$, which entails (\ref{flaccid}). 

We now complete the proof of the lemma. 
Since $\Lambda^{\! \mathtt{b} , \bw}\!  + \Lambda^{\! \mathtt{r} , \bw}$ is the identity map, we get 
$ \theta^{\mathtt{r} , \bw}_{t}\! = \!  \Lambda^{\! \mathtt{b} , \bw} ( \theta^{\mathtt{r} , \bw}_{t})+t$, which shows that $\theta^{\mathtt{r} , \bw}$ is super-linear. 
Next observe that if $\Delta \theta^{\mathtt{r}, \bw}_t \! >\! 0$, then there exists $l\! < \! l^*$such that $\theta^{\mathtt{r}, \bw}_t \! =\!  t_{2l+1}$ and $t\! = \! \Lambda^{\mathtt{r}, \bw} (t_{2l+1})\! = \! s_l^\prime\! = \! \gamma^{\mathtt{r}, \bw} (A^\bw_{s_l})$. Thus, we get $\Delta X^{\mathtt{r}, \bw}_t \! = \! 0$, which completes the proof of the lemma. \cqfd

\medskip

\noi
\textbf{Proof of Lemma \ref{Hthetalem}.} Namely, we prove a.s.~that 
for all $a\ino [0, T^*_{\bw})$, $\cH^{\bw}_a\! = \! H^\bw (\theta^{\mathtt{b}, \bw}_a)$. 
To that end, we fix $a\ino [0, T^*_{\bw})$ and we set $\cJ_a\! = \! \{ r\ino [0, a]\! : \! J^{\bw, r-}_a \! \! < \! J^{\bw, r}_a \}$ where $J^{\bw, r}_a \! = \! \inf_{[r, a]} Y^\bw$. By the definition (\ref{JHdef2}) of $\cH$, we get $\# \cJ_a \! = \! \cH^\bw_a$. For any $t\ino [0, \infty)$, we next set $\cK (t)\! = \! \{ s\ino [0, t]\! : \! 
I^{\bw, s-}_t \! < \! I^{\bw, s}_t \}$ where $I^{\bw, s}_t \! = \! \inf_{[s, t]} X^\bw$.
By the definition (\ref{XJHdef}), we get $\# \cK (t) \! = \! H^\bw_t $. Then, it is sufficient to prove that $\Lambda^{\mathtt{b}, \bw}$ is one-to-one from $\cK (\theta^{\mathtt{b}, \bw}_a)$ onto $\cJ_a$. 

We first prove: $\cJ_a\! \subset \!  \Lambda^{\mathtt{b}, \bw} (\cK (\theta^{\mathtt{b}, \bw}_a))$. First observe that $a\! = \! \Lambda^{\mathtt{b}, \bw} (\theta^{\mathtt{b}, \bw}_a)$,  
which implies that, when $a\in \cJ_a$, we also have $a\ino \Lambda^{\mathtt{b}, \bw} (\cK (\theta^{\mathtt{b}, \bw}_a))$. Let $r\ino \cJ_a$ and suppose that $r\! < \! a$. By (\ref{flaccid}), we get $J^{\bw, r}_a\! = \! \inf \{ X^\bw_s \! ;  
s\ino [\theta^{\mathtt{b}, \bw}_r \! ,\theta^{\mathtt{b}, \bw}_a] \}$. 
Moreover, since $r\ino \cJ_a$, $\Delta Y^{\bw}_r\! > \! 0$. Suppose that $\theta^{\mathtt{b}, \bw}_{r-}\! < \! \theta^{\mathtt{b}, \bw}_{r}$, then (\ref{flaccidos}) would easily entail 
that $\Delta Y^{\bw}_r\! = \! \Delta X^{\bw} (\theta^{\mathtt{b}, \bw}_r)\! = \! 0$. Thus, $\theta^{\mathtt{b}, \bw}_{r-}\! = \! \theta^{\mathtt{b}, \bw}_{r}$ and $ \theta^{\mathtt{b}, \bw}_{r}$ lies therefore 
in the interior of $\mathtt{Blue}$; therefore for all sufficiently small $\epp \! >\! 0$, $\theta^{\mathtt{b}, \bw}_{r-\epp} \! = \! \theta^{\mathtt{b}, \bw}_{r}\!\! - \! \epp $, 
and thus, 
$J^{\bw, r-\epp}_a\! = \! \inf \{ X^\bw_s \! ; s \ino [\theta^{\mathtt{b}, \bw}_r \!\! -\! \epp ,
\theta^{\mathtt{b}, \bw}_a] \}$. 
This implies that $ \theta^{\mathtt{b}, \bw}_{r} \ino   \cK (\theta^{\mathtt{b}, \bw}_a)$ and it completes the proof that $\cJ_a\! \subset \!  \Lambda^{\mathtt{b}, \bw} (\cK (\theta^{\mathtt{b}, \bw}_a))$. 

Conversely, let us prove that $\Lambda^{\mathtt{b}, \bw} (\cK (\theta^{\mathtt{b}, \bw}_a)) \! \subset \! \cJ_a$. Let $s\ino \cK (\theta^{\mathtt{b}, \bw}_a)$. To simplify notation we set $r\! = \!  \Lambda^{\mathtt{b}, \bw}_{s}$. Note that $r\! < \! T^*_{\bw}$. 
We first prove that $s$ lies in the interior of $\mathtt{Blue}$. Indeed, suppose the contrary: namely suppose that $s$ lies in the closure of $\mathtt{Red}$; then we would get 
$\theta^{\mathtt{b}, \bw}_{r-}\!  < \!  \theta^{\mathtt{b}, \bw}_{r}$ and $s\ino [\theta^{\mathtt{b}, \bw}_{r-} ,  \theta^{\mathtt{b}, \bw}_{r}]$, and (\ref{flaccidos}) would entails that $ X^{\bw}_{s-}\! \geq \! 
X^\bw (\theta^{\mathtt{b}, \bw}_r)\! \geq \! I^{\bw , s}_{\theta^{\mathtt{b}, \bw}_a}$ and thus, $I^{\bw , s}_{\theta^{\mathtt{b}, \bw}_a}\! = \! I^{\bw , s-}_{\theta^{\mathtt{b}, \bw}_a}$, which contradicts 
$ I^{\bw , s-}_{\theta^{\mathtt{b}, \bw}_a}\! < \!  I^{\bw , s}_{\theta^{\mathtt{b}, \bw}_a}$. Thus, $s$ lies in the interior of $\mathtt{Blue}$. Consequently, for all sufficiently small $\epp\ino (0, \infty)$, $\Lambda^{\mathtt{b}, \bw}_{s-\epp}\! = \! \Lambda^{\mathtt{b}, \bw}_{s}\! - \! \epp\! = \! r\! -\! \epp$ and (\ref{flaccid}) entails 
$J^{\bw, r-\epp }_a\! = \!  I^{\bw , s-\epp}_{\theta^{\mathtt{b}, \bw}_a}$, which implies that $r\! = \!  \Lambda^{\mathtt{b}, \bw}_{s} \ino \cJ_a$ because
$$  J^{\bw, r- }_a= I^{\bw, s-}_{\theta^{\mathtt{b}, \bw}_a } \! < \!  I^{\bw, s}_{\theta^{\mathtt{b}, \bw}_a }=J^{\bw, r }_a \; .$$
 
We have proved that 
$\Lambda^{\mathtt{b}, \bw} (\cK (\theta^{\mathtt{b}, \bw}_a)) \!=\! \cJ_a$. But recall that we have also proved that $\cK (\theta^{\mathtt{b}, \bw}_a)$ is included in the interior of $\mathtt{Blue}$ on which $\Lambda^{\mathtt{b}, \bw}$ is one-to-one, by definition. Consequently, we get 
$H^\bw (\theta^{\mathtt{b}, \bw}_a)\! = \! \# \cK (\theta^{\mathtt{b}, \bw}_a) \! = \! \# \cJ_a \! = \! \cH^{\bw}_a$, which proves Lemma \ref{Hthetalem}. \cqfd  

\begin{rem}
\label{pluahisv}
By (\ref{flaccidos}), we actually get $H^\bw_t \! = \!  \cH^{\bw}(\Lambda^{\mathtt{b}, \bw}_t)+ H^{\mathtt{r}, \bw} (\Lambda^{\mathtt{r}, \bw}_t)$, where we have set 
$H^{\mathtt{r}, \bw} _t \! =\!  \# \big\{ s\ino [0, t]\!  : \! I^{\prime \, \bw,s-}_{t} \! <\! I^{\prime\, \bw, s}_{t} \big\}$ and $ I^{\prime \, \bw, s}_{t}\! = \! \inf_{[s, t]} X^{\mathtt{r}, \bw}$. \cq
\end{rem}

%
%
%
%
%


\noi
\textbf{End of the proof of Proposition~\ref{Xwfrombr}.} It remains to prove that $X^\bw$ has the same law as $X^{\mathtt{b}, \bw}$ (or $X^{\mathtt{r}, \bw}$). To that end, we first define a point measure that records the successive jump-times of $X^\bw$, along with their types and their colours. Actually, here it will be convenient to distinguish three colours: $\mathtt{b}$ (\textit{blue}) for the jumps of $Y^\bw$, $\mathtt{b}^\prime$ 
(\textit{repeated blue}) for the jumps of $A$ and $\mathtt{r}$ (\textit{red}) for the jumps of $X^{\mathtt{r},\bw}$. More precisely, set $k_* \! =\! \inf \{ k \! \geq \! 1: \theta^{\mathtt{b}, \bw} (\tau^{\mathtt{b}}_k)\! = \! \infty \}$ (that is infinite in the critical or subcritical cases). 
Recall from (\ref{ccXbrwdef}) that 
$\ccX_\bw^{\mathtt{b}}\! = \! \sum_{k\geq 1} \delta_{(\tau^{\mathtt{b}}_k , \Jtt^{\mathtt{b}}_k)}$. 
Let us fix $k \! < \!  k_*$. 

\smallskip

\begin{compactenum}
\item[\texttt{Colour} $\mathtt{b}$:] If $\tau^{\mathtt{b}}_k$ is a jump-time of $Y^\bw$, then $\tau^{\mathtt{b}}_k$ is in the interior of $\mathtt{Blue}$: namely, $\theta^{\mathtt{b}, \bw} (\tau^{\mathtt{b}}_k-)\! = \! \theta^{\mathtt{b}, \bw} (\tau^{\mathtt{b}}_k)$. On the interior of $\mathtt{Blue}$, $\theta^{\mathtt{b}, \bw}$ is continuous. Thus, if $\tau^{\mathtt{b}}_k$ is a jump-time of $Y^\bw$, then $\tau^{\mathtt{b}}_k$ corresponds to the jump-time $\theta^{\mathtt{b}, \bw}(\tau^{\mathtt{b}}_k)$ of $X^\bw$.

\smallskip

\item[\texttt{Colour} $\mathtt{b}^\prime$:]
 If $\tau^{\mathtt{b}}_k$ is a jump-time of $A^\bw$, then $\tau^{\mathtt{b}}_k$ is an endpoint of a connected component of $\mathtt{Blue}$: namely, 
 $\theta^{\mathtt{b}, \bw} (\tau^{\mathtt{b}}_k-)\! < \! \theta^{\mathtt{b}, \bw} (\tau^{\mathtt{b}}_k)$,  
 $\theta^{\mathtt{b}, \bw} (\tau^{\mathtt{b}}_k-)$ is a time  
when the queue changes its colour and 
the corresponding jump-time in $X^\bw$ is $\theta^{\mathtt{b}, \bw}(\tau^{\mathtt{b}}_k-)$. By the definition (\ref{telltell}) 
of the times $t_{l}$, there exists $l$ such that $\theta^{\mathtt{b}, \bw}(\tau^{\mathtt{b}}_k-)\! = \! t_{2l+1}$ and $\theta^{\mathtt{b}, \bw}(\tau^{\mathtt{b}}_k)\! =\! t_{2l+2}$.

\smallskip

\item[\texttt{Colour} $\mathtt{r}$:] Recall from (\ref{ccXbrwdef}) that 
$\ccX_\bw^{\mathtt{r}}\! = \! \sum_{k\geq 1} \delta_{(\tau^{\mathtt{r}}_k , \Jtt^{\mathtt{r}}_k)}$ and recall from (\ref{thetarw}) the red time-change $\theta^{\mathtt{r}, \bw}$. 
Then, for all $k\! \geq \! 1$, $\tau^{\mathtt{r}}_k$ corresponds to the jump-time 
$\theta^{\mathtt{r}, \bw} (\tau^{\mathtt{r}}_k)$ of the process $X^\bw$. 
\end{compactenum}

\smallskip

\noi
We then define as follows the 
increasing sequence of jump-times $\tau_p$, of $X^\bw$, the corresponding types $\Jtt_p$ and their colour $\mathtt{col} (p)\ino \{ \mathtt{b},\mathtt{b}^\prime, \mathtt{r} \}$.
\begin{eqnarray}
 \overline{\ccX}_\bw & = & \!\!\! \!  \sum_{1\leq k<k_*} \!\!\! \un_{\{ \Delta Y^\bw (\tau^{\mathtt{b}}_k)\neq 0 \}} \, 
\delta_{ \big(\theta^{_{\mathtt{b}, \bw}}_{^{\tau^{_{\mathtt{b}}}_{^k}}},\,   \Jtt^{\mathtt{b} }_k \, ,\,  \mathtt{b} \big)} + \!\!\! \sum_{1\leq k<k_*} \!\!\! \un_{\{ \Delta A^\bw (\tau^{\mathtt{b}}_k)\neq 0 \}} \, 
\delta_{ \big(\theta^{\mathtt{b}, \bw}_{\tau^{\mathtt{b}}_k-}, \,  \Jtt^{\mathtt{b}_k} \, ,\,  \mathtt{b}^\prime \big)} 
+ \!  \sum_{k\geq 1}\delta_{ \big(\theta^{\mathtt{r}, \bw}_{\tau^{\mathtt{r}}_k}, \, \Jtt^{\mathtt{r}_k} \, ,\,   \mathtt{r} \big)} \\ \notag
& =:& \sum_{p\geq 1} \delta_{\big( \tau_p , \Jtt_p , \, \mathtt{col} (p)\big)}  \; .
\end{eqnarray}
See Fig.~\ref{fig:colouring} for an example.
We get: 
$$ \forall t\ino [0, \infty), \quad X^\bw_t = -t + \sum_{p\geq 1} w_{\Jtt_p} \un_{[0, t]} (\tau_p) \; .$$
To prove Proposition \ref{Xwfrombr}, it is then sufficient to prove that $\ccX_\bw\! := \! \sum_{q\geq 1} \delta_{(\tau_q,  \Jtt_q)}$ is a Poisson point measure on  $[0, \infty) \! \times \! \{ 1, \ldots , n\}$ with intensity 
$\ell \otimes \nu_\bw$. 

To that end, we introduce the following notation: let $M\! = \! \sum_{q\geq 1} \delta_{(r_q, j_q)}$ be a point measure on $[0, \infty) \! \times \! \{ 1, \ldots , n\}$ and let 
$a\ino [0, \infty)$; then, we denote by $\mathtt{R} (M, a)\! = \! \sum_{q\geq 1} \un_{[0, a ]} (r_q)  \delta_{(r_q, j_q)}$ the \textit{restriction} of $M$ on $[0, a]\! \times \! \{ 1, \ldots , n\}$ and we denote 
by $\mathtt{S} (M, a)\! = \! \sum_{q\geq 1} \un_{(a, \infty )} (r_q)  \delta_{(r_q-a, j_q)}$ the \textit{$a$-time shifted} measure $M$.  
\begin{lem}
\label{kalmarr} We fix an integer $p\! \geq \! 1$. Then, 
conditionally given $\ccZ\! := \! \mathtt{R} (\ccX_\bw^{\mathtt{b}}, \Lambda^{\mathtt{b}, \bw}_{\tau_p})$ and $\ccZ^\prime \! := \!\mathtt{R} (\ccX_\bw^{\mathtt{r}}, \Lambda^{\mathtt{r}, \bw}_{\tau_p})$, the shifted measures 
$\ccY\! := \!\mathtt{S} (\ccX_\bw^{\mathtt{b}}, \Lambda^{\mathtt{b}, \bw}_{\tau_p})$ and $\ccY^\prime\! := \!\mathtt{S} (\ccX_\bw^{\mathtt{r}}, \Lambda^{\mathtt{r}, \bw}_{\tau_p})$ are independent Poisson point measures on  $[0, \infty) \! \times \! \{ 1, \ldots , n\}$ with intensity 
$\ell \otimes \nu_\bw$. 
\end{lem}
\noi
\textbf{Proof.} We introduce $\mathtt{red} (p)\! = \! \# \big\{ l \ino \{ 1, \ldots , p\} \! :\!  \mathtt{col} (l)\! =\!  \mathtt{r} \big\}$ that counts the number of red jump-times of the process $X^\bw$ among its first $p$ jump-times. For all $k \ino \{ 0, \ldots , p\}$, we introduce the events 
$A(k,p)\! = \{   \mathtt{red} (p)\! = \! k ; \mathtt{col} (p)\! = \! \mathtt{r}\}$ and 
$A^\prime(k,p)\! = \{   \mathtt{red} (p)\! = \! k ; \mathtt{col} (p)\! \neq \! \mathtt{r}\}$. 

We first observe that on $A(k,p)$, $\Lambda^{\mathtt{r}, \bw}_{\tau_p}\! = \! \tau^{\mathtt{r}}_k$ and $\Lambda^{\mathtt{b}, \bw}_{\tau_p}\! = \! \tau^{\mathtt{b}}_{p-k}$. Namely, $\ccZ\! = \! \mathtt{R} (\ccX_\bw^{\mathtt{b}},  \tau^\mathtt{b}_{p-k})$, 
$\ccZ^\prime \! = \!\mathtt{R} (\ccX_\bw^{\mathtt{r}},  \tau^\mathtt{r}_{k})$, 
$\ccY\! = \!\mathtt{S} (\ccX_\bw^{\mathtt{b}}, \tau^{\mathtt{b}}_{p-k})$ and 
$\ccY^\prime\! = \!\mathtt{S} (\ccX_\bw^{\mathtt{r}}, \tau^\mathtt{r}_{k})$. 
To simplify, denote by $\ccH$ the sigma-field generated by
$\mathtt{R} (\ccX^\mathtt{b}_\bw, \tau^\mathtt{b}_{p-k})$ and $\mathtt{R} (\ccX^\mathtt{r}_\bw, 
\tau^\mathtt{r}_{k})$. By elementary properties of Poisson point measures, 
conditionally given $\ccH$, $\mathtt{S} (\ccX^\mathtt{b}_\bw, \tau^\mathtt{b}_{p-k})$ and $\mathtt{S} (\ccX^\mathtt{r}_\bw, \tau^\mathtt{r}_{k})$ are independent Poisson point measures on  $[0, \infty) \! \times \! \{ 1, \ldots , n\}$ with intensity $\ell \otimes \nu_\bw$. 
Then, remark that $A(k,p)$ belongs to $\ccH$. Thus, we get 
\begin{equation}
\label{raggouni}
\bE [G(\ccZ, \ccZ^\prime) F( \ccY, \ccY^\prime)\un_{A(k,p)}]\! = \!  \bE [ G(\ccZ, \ccZ^\prime) \un_{A(k,p)}]\bE [F( \ccX^\mathtt{b}_\bw, \ccX^\mathtt{r}_\bw) ] \; .
\end{equation}

Next observe that on $A^\prime (p,k)$, 
$\Lambda^{\mathtt{b}, \bw}_{\tau_p}\! = \! \tau^{\mathtt{b}}_{p-k}$ and 
$\Lambda^{\mathtt{r}, \bw}_{\tau_p}\! = \! \gamma^{\mathtt{r}, \bw} (A^\bw (\tau^\mathtt{b}_{p-k}\! - )) \! < \! \infty $, where $\gamma^{\mathtt{r}, \bw}$ is given by (\ref{rethetabw}). 
Thus, $\ccZ\! = \! \mathtt{R} (\ccX_\bw^{\mathtt{b}},  \tau^\mathtt{b}_{p-k})$, 
$\ccZ^\prime \! = \!\mathtt{R} (\ccX_\bw^{\mathtt{r}}, \gamma^{\mathtt{r}, \bw} (A^\bw (\tau^\mathtt{b}_{p-k}\! - )))$, 
$\ccY\! = \!\mathtt{S} (\ccX_\bw^{\mathtt{b}}, \tau^{\mathtt{b}}_{p-k})$ and 
$\ccY^\prime\! = \!\mathtt{S} (\ccX_\bw^{\mathtt{r}},\gamma^{\mathtt{r}, \bw} (A^\bw (\tau^\mathtt{b}_{p-k}\! -\! )))$. Then, for all $a\ino [0, \infty)$, denote by $\ccG^{\mathtt{r}}_a$ 
(resp.~$\ccG^{\mathtt{b}}_a$) the sigma-field generated by $\mathtt{R} (\ccX_\bw^{\mathtt{r}},a)$ (resp.~$\mathtt{R} (\ccX_\bw^{\mathtt{b}},a)$). Note that for all $x\ino [0, \infty)$, $\gamma_x^{\mathtt{r}, \bw}$ is a$(\ccG^\mathtt{r}_a)$-stopping time and by standard results, conditionally given $\ccG^{\mathtt{r}} (\gamma^{\mathtt{r}, \bw}_x)$ under $\bP (\, \cdot \, |\,  
\gamma_x^{\mathtt{r}, \bw} \! <\! \infty)$, $
\mathtt{R} (\ccX_\bw^{\mathtt{r}}, \gamma_x^{\mathtt{r}, \bw})$ and $
\mathtt{S} (\ccX_\bw^{\mathtt{r}}, \gamma_x^{\mathtt{r}, \bw})$ are independent and 
the conditional law of $
\mathtt{S} (\ccX_\bw^{\mathtt{r}}, \gamma_x^{\mathtt{r}, \bw})$ is that of a Poisson point measure on  $[0, \infty) \! \times \! \{ 1, \ldots , n\}$ with intensity 
$\ell \otimes \nu_\bw$. To simplify, denote by $\ccH^\prime$ the sigma-field generated by $\ccG^\mathtt{b} (\tau^\mathtt{b}_{p-k})$ and $\ccG^{\mathtt{r}} (\gamma^{\mathtt{r}, \bw} (A^\bw (\tau^\mathtt{b}_{p-k}\! - )))$. Note that $A^\prime (p,k)\ino \ccH^\prime$, that $\mathtt{R} (\ccX_\bw^{\mathtt{b}},  \tau^\mathtt{b}_{p-k})$ and $\mathtt{R} (\ccX_\bw^{\mathtt{r}} , \gamma^{\mathtt{r}, \bw} (A^\bw (\tau^\mathtt{b}_{p-k}\! - )))$ are $\ccH^\prime$-measurables and that conditionally given $\ccH^\prime$,  
$\mathtt{S} (\ccX_\bw^{\mathtt{b}}, \tau^{\mathtt{b}}_{p-k})$ and 
$\mathtt{S} (\ccX_\bw^{\mathtt{r}},\gamma^{\mathtt{r}, \bw} (A^\bw (\tau^\mathtt{b}_{p-k}\! -\! )))$
are distributed as two independent 
Poisson point measures on $[0, \infty) \! \times \! \{ 1, \ldots , n\}$ with intensity 
$\ell  \otimes  \nu_\bw$. Thus, it implies that 
$$\bE [G(\ccZ, \ccZ^\prime) F( \ccY, \ccY^\prime) \un_{A^\prime(k,p)}]= 
\bE [ G(\ccZ, \ccZ^\prime)\un_{A^\prime(k,p)}] \bE [F( \ccX^\mathtt{b}_\bw, \ccX^\mathtt{r}_\bw) ] \; .$$
Since the events $A(k,p), A^\prime(k,p)$, $k\ino \{ 0, \ldots , p\}$ form a partition of $\Omega$, the previous equality combined with (\ref{raggouni}) entails the desired result. \cqfd   

\medskip

We fix $p\! \geq \! 1$ and we recall notation $\ccZ, \ccZ^\prime,\ccY, \ccY^\prime$ from Lemma \ref{kalmarr}. 
Denote by $(e^{\mathtt{b}} , J^\mathtt{b})$ the first atom of $\ccY$ and denote by $(e^{\mathtt{r}} , J^\mathtt{r})$ the first atom of $\ccY^\prime$. Recall the notation $\ccX_\bw\! = \! \sum_{q\geq 1} \delta_{(\tau_q, \Jtt_q)}$ and observe that $\mathtt{R}(\ccX_\bw, \tau_p)$ is a deterministic function of $\ccZ$ and $\ccZ^\prime$. Then, by Lemma \ref{kalmarr}, 
\begin{equation}
\label{lubeck}
 (e^{\mathtt{b}} , J^\mathtt{b}), \; (e^{\mathtt{r}} , J^\mathtt{r})\quad  \textrm{and}  
\quad \mathtt{R}(\ccX_\bw, \tau_p)\quad \textrm{are independent} 
\end{equation}
and $\bP (e^{\mathtt{b}}\! >\! t ; 
J^\mathtt{b}\! = \! j)\! =\!\bP (e^{\mathtt{r}}\! >\! t ; 
J^\mathtt{r}\! = \! j)\! =\! e^{-t} \nu_\bw (j)$, for all $t\ino [0, \infty)$ and all $j\ino \{ 1, \ldots, n\}$.    

We next explain how $(\tau_{p+1} \! -\! \tau_p , \Jtt_{p+1})$ (that is the first atom of $\mathtt{S} (\ccX_\bw, \tau_p)$) is derived from $X^\bw_{\cdot \wedge \tau_p}$, $(e^{\mathtt{b}} , J^\mathtt{b})$ and 
$(e^{\mathtt{r}} , J^\mathtt{r})$. First observe that 
$p\! -\! \mathtt{red}(p)$ is the number of jumps of $X^\bw$ whose colour is $\mathtt{b}$ or $\mathtt{b}^\prime$. Therefore, the time $\tau\! := \! \theta^{\mathtt{b}, \bw} (\tau^{\mathtt{b}}_{p-\mathtt{red} (p)}-)$ is the largest non-red jump-time of $X^\bw$ before $\tau_p$ (thus, if $\mathtt{col} (p)\ino \{ \mathtt{b}, \mathtt{b}^\prime \}$, then $\tau\! = \! \tau_p$). We next set $\Delta:= X^\bw (\tau_p)\! -\! X^{\bw} (\tau-)$. By construction of $X^\bw$, we easily check that $\Delta \! >\! 0$ and that the following holds true (see also Fig.~\ref{fig:colouring} for an example): 

\smallskip
 
\begin{compactenum}

\item[$\bullet$] $(\tau_{p+1} \! -\! \tau_p , \Jtt_{p+1})\! = \! (e^\mathtt{b}, J^\mathtt{b})  $ if $\mathtt{col} (p)\! = \! 
\mathtt{b}$.   

\smallskip

\item[$\bullet$] $(\tau_{p+1} \! -\! \tau_p , \Jtt_{p+1})\! = \! (e^\mathtt{r}, J^\mathtt{r}) $ if $\mathtt{col} (p)\ino \{ \mathtt{r} , \mathtt{b}^\prime\}$ and $e^\mathtt{r} \! \leq \! \Delta $.   

\smallskip

\item[$\bullet$] $(\tau_{p+1} \! -\! \tau_p ,  \Jtt_{p+1})\! = \! (\Delta +e^\mathtt{b}, J^\mathtt{b}) $ 
if $\mathtt{col} (p)\ino \{ \mathtt{r} , \mathtt{b}^\prime\}$ and $e^\mathtt{r} \! > \!  \Delta $.  

\smallskip

\end{compactenum}

Since $(\Delta, \mathtt{col} (p))$ are deterministic functions of $(\ccZ, \ccZ^\prime)$, elementary properties of exponentially distributed r.v.~combined with (\ref{lubeck}) entail that $ \mathtt{R}(\ccX_\bw, \tau_p)$ and $(\tau_{p+1}  - \tau_p , J_{p+1})$ are independent and that 
$\bP \big( \tau_{p+1} \! -\! \tau_p \! >\! t ; J_{p+1}\! = \! j) \! =\! e^{-t} \nu_\bw (j)$, for all $t\ino [0, \infty)$ and for all $j\ino \{ 1, \ldots, n\}$. This easily implies that $\ccX_\bw$ is a  Poisson point measure on  $[0, \infty) \! \times \! \{ 1, \ldots , n\}$ with intensity 
$\ell \otimes \nu_\bw$. Consequently, $X^\bw$ has the same distribution as $X^{\mathtt{b}, \bw}$ (or as $X^{\mathtt{r}, \bw}$), which completes the proof of Proposition \ref{Xwfrombr}. \cqfd

\section{Properties of the continuous processes}\label{sec:continuous}
\subsection{The height process of a L\'evy tree}
\label{Heightsec}
In this section, we collect the various properties of the height process associated with a L\'evy process that will be used later. In Section \ref{Heightsec}, there is no new result, the only exception being the technical Lemma \ref{cucueille}.

\subsubsection{Infimum process of a spectrally positive L\'evy process}
\label{infimsec}

We first recall several results on the following specific L\'evy processes : we fix 
$\alpha\ino \bbR$, $\beta \ino [0, \infty) $, $\kappa \ino (0, \infty) $, $\mathbf{c}\! = \! (c_j)_{j\geq 1} \ino \elldo_3$ and we set  
\begin{equation}
\label{repsidefi}
\forall \lambda \ino [0, \infty) , \quad \psi(\lambda)  \! = \!   
 \alpha \lambda +\frac{_{_1}}{^{^2}} \beta \lambda^2 +  \!  \sum_{j\geq 1}   \kappa c_j   \big( e^{-\lambda c_j}\! -\! 1\! + \! \lambda c_j \big). 
 \end{equation}
Let $(X_t)_{t\in [0, \infty)}$ be a spectrally positive L\'evy process with initial state $X_0 \! = \! 0$ and with Laplace exponent $\psi$: namely, 
$ \log \bE [ \exp ( - \lambda X_t )] \! = \!   t\psi (\lambda) $, for all $t, \lambda \ino [0, \infty)$. The L\'evy measure of $X$ is 
$\pi \! = \! \sum_{j\geq 1} \kappa c_j \delta_{c_j}$, $\beta$ is its Brownian parameter and $\alpha$ is its drift. 
Note that it includes the Lévy process $X^{\bw}$ by taking $\mathbf{c}\! = \! \bw\ino \elldo_{\! f}$, 
$\kappa \! = \! 1/\sigma_1 (\bw)$, $\beta\! = \! 0$ and $\alpha\! = \!  1\! -\! 
\frac{{\sigma_2 (\bw)}}{{\sigma_1 (\bw)}}$.

Recall from (\ref{varinfinie}) that $X$ has infinite variation sample paths 
if and only if either $\beta\! > \! 0$ or 
 $\sigma_2 (\mathbf{c})\! = \! \int_{(0, \infty)}\!  r\, \pi (dr) \! = \! \infty$, which is implied by the following condition: 
 \begin{equation}
\label{recontH}
 \int^\infty \frac{d\lambda}{\psi (\lambda)} <\infty \; , 
\end{equation}
that is assumed in various places in the paper.   

Recall from (\ref{thetabdef}) the notation: $\subo_x \! =\!  \inf \{ s \ino [0, \infty):    X_s \! < \! -x \}$, for all $\ino [0, \infty)$, with the convention: $\inf \emptyset \! = \! \infty$.
For all $t\ino [0, \infty)$, we set $I_t \! = \! \inf_{s \in [0, t]} X_s $ and 
 $I_\infty\! = \! \lim_{t\rightarrow \infty} \! I_t$ that is a.s.~finite in supercritical cases and a.s.~infinite in critical or subcritical cases. 
Observe that 
$\subo_x\! < \! \infty$ iff $x\! < \! \! -I_\infty $.  Standard results on spectrally positive L\'evy processes (see e.g.~Bertoin's book \cite{Be96} Ch.~VII) assert that 
$(\subo_x)_{x\in [0, \infty)}$ is a c\`adl\`ag subordinator (a killed subordinator in supercritical cases) whose Laplace exponent is given for all 
$\lambda \ino [0, \infty)$ by: 
\begin{equation}
\label{gamexpobs}
\bE \big[ e^{-\lambda \subo_x } \big]= e^{-x\psi^{-1}(\lambda)} \quad  \textrm{where} \quad  \psi^{-1} (\lambda)\! = \! \inf \big\{ u \ino [0, \infty ) : \psi (u) \! >\! \lambda \big\}.  
\end{equation}
We set $\varrho \! = \!  \psi^{-1}(0)$ that is the largest root of $\psi$. Note that $\varrho \! > \! 0 $ iff $\alpha \! < \! 0$. 
%
Let us also recall from Chapter VII in Bertoin's book \cite{Be96} the following absolute continuity relationship between a supercritical spectrally positive 
L\'evy process and a subcritical one.  More precisely,
%
for all $t\ino [0, \infty)$ and for all nonnegative measurable functional $F \! : \! \bD ([0, \infty), \bbR) \! \rightarrow \! \bbR$, 
\begin{equation}
\label{supsubis}
\bE \big[ F( X_{\cdot \wedge t} ) ] \! =\! \bE \big[ \exp (\varrho \overline{X}_t ) \,  F( \overline{X}_{\cdot \wedge t} )\big]  , 
\end{equation}
where $\overline{X}$ stands for a \textit{subcritical} spectrally L\'evy process with Laplace exponent $\psi (\varrho + \cdot)$. 

\subsubsection{Local time at the supremum}
\label{locsupsec} 
We assume (\ref{varinfinie}), namely that $X$ has infinite variation sample paths. 
For all $t\ino [0, \infty)$, we set $S_t \! = \! \sup_{s\in [0, t]} X_s$. Basic results of fluctuation theory 
entail that $S\! -\! X$ is a strong Markov process in $[0, \infty)$ and that $0$ is regular for 
$(0, \infty)$ and recurrent with respect to this Markov process (see for instance Bertoin \cite{Be96} VI.1). 
We denote by $(L_t)_{t\in [0, \infty)}$ the local time of $X$ at its supremum (namely, the local time of $S\! -\! X$ at $0$), 
whose normalisation is such that for all $t\ino [0, \infty)$ the following holds in probability: 
\begin{equation}
\label{approL}
 L_t= \lim_{\varepsilon \rightarrow 0} \frac{1}{\varepsilon} \int_0^t \!\!  \un_{\{ S_s -X_s \leq \varepsilon  \}} \, ds \; .
\end{equation}
See Le Gall \& D.~\cite{DuLG02} (Chapter 1, Lemma 1.1.3 p.~21) for more details.  If $\beta \! >\! 0$, then standard results about subordinators 
imply that a.s.~for all $t\ino [0, \infty)$, $L_t \! =\!  \frac{2}{\beta} \mathtt{Leb}( \{ S_s ; s\ino [0, t]\})$. 
We also need the following approximation of $L$ that holds when 
$\sigma_2 (\mathbf{c})\! = \! \infty$: for all $\varepsilon \ino (0 , c_1)$, we set 
\begin{equation}
\label{kepsdef}
q(\varepsilon)\! =\!  \int_{(\varepsilon , \infty)}\!\!\!\!\!\!  \!\!\! dx\,  \pi ((x, \infty))  =\! \sum_{j\geq 1} \kappa c_j (c_j \! -\! \varepsilon)_+ \quad \textrm{and} \quad 
\ccL^\varepsilon_t = \big\{ s \ino (0, t] : S_{s-}\!  + \varepsilon  < \! X_s \big\} \; .
\end{equation}
If $\sigma_2 (\mathbf{c})\! = \! \infty$, then the following approximation holds true:
\begin{equation}
\label{L2approx}
\forall x, t\ino (0, \infty) , \quad \bE \Big[  \un_{\{ L_t \leq x \}} \sup_{s\in [0, t]} \big| L_s \! -\!  \frac{_1}{^{q(\varepsilon)}} \# \ccL^\varepsilon_s  \big|^2 \Big] \leq \frac{4x}{q(\varepsilon)}.   
 \end{equation}
This is a standard consequence of the decomposition of $X$ into excursions under its supremum: see Bertoin \cite{Be96}, Chapter VI.


\subsubsection{The height process}
\label{hautruc} 
For all $t\ino (0, \infty)$, we denote by $\widehat{X}^t\! = \! (X_t \! -\! X_{(t-s)-})_{s\in [0, t]}$ that is the process $X$ reversed at time $t$; recall that $\widehat{X}^t$ has the same law as $(X_s)_{s\in [0, t]}$. 
Under (\ref{recontH}), Le Gall \& Le Jan \cite{LGLJ98} (see also Le Gall \& D.~\cite{DuLG02}) 
prove that there exists a \textit{continuous} process $H\! = \! (H_t)_{t\in [0, \infty)}$ such that \begin{equation}
\label{Hdef}
\forall t\ino [0, \infty), \; \textrm{a.s.} \quad H_t = L_t (\widehat{X}^t) \; .
\end{equation}
As previously mentioned in Remark \ref{supercr}, this result has been established only in subcritical or critical cases. By (\ref{supsubis}), it easily extends to the supercritical case. 
In particular, one can show that 
\begin{equation}
 \label{glurniglou}
\textrm{If $\alpha \! < \! 0$, then a.s.} \quad \lim_{t\rightarrow \infty} H_t = \infty \; .
\end{equation}

We next derive from the two approximations of $L$ mentionned in Section \ref{locsupsec}, similar approximations for the height process $H$. 
For all real numbers $t \! \geq \! s  \!  \geq \! 0$, we first introduce the following: 
\begin{equation}
\label{infideff}
  I_t^s \! = \! \inf_{r\in [s, t]} X_r,\quad   I_t \! =\!  I_t^0 \! = \! \! \inf_{r\in [0, t]} \!\! X_r  \quad   \textrm{and} \quad  \ccH_t^\varepsilon \! = \! \big\{s\ino (0, t]: X_{s-} + \varepsilon < I_t^s   \big\} \;.
\end{equation}
%
Then we easily derive from (\ref{Hdef}) and (\ref{approL}) that (\ref{approHdef}) holds true: namely, for all $t\ino [0, \infty)$, the following limit holds in probability:   
$$H_t \! =\!  \lim_{\varepsilon \rightarrow 0} \frac{1}{\varepsilon} \!  \int_0^t \!  \un_{\{  X_s - I_t^s  \leq \varepsilon\}} \, ds \; .$$ 
Of course, we also get the following. 
\begin{equation}
\label{brocueille} \textrm{If $\beta \! >\! 0$, then a.s.~for all $t\ino [0, \infty)$,} \quad 
H_t  \! = \! \frac{_2}{^\beta} \, \mathtt{Leb} \big( \{ I_t^s ; s\ino [0, t]\} \big) \; .  
\end{equation}
If $\sigma_2 (\mathbf{c}) \! = \! \infty$, then  (\ref{Hdef}) and (\ref{L2approx}) easily imply 
that for all $t\ino [0, \infty)$, $H_t \! =\!  \lim_{\varepsilon \rightarrow 0} \frac{_1}{^{q(\varepsilon)}}  \# \ccH^{\varepsilon}_t$ in probability.  
Actually, a closer look at the uniform approximation (\ref{L2approx}) shows the following.  
\begin{multline}
\label{cueille}
\textrm{If $\sigma_2 (\mathbf{c}) \! = \! \infty$, then $\forall t\ino [0, \infty)$, $\exists (\varepsilon_k)_{k\in \bbN}$ 
decreasing to $0$ such that:} \\
\textrm{$\bP$-a.s.~for all $s\ino [0, t]$ such that $X_{s-} \! \leq \! I_t^s$} , \quad H_s = \lim_{k \rightarrow \infty} \frac{_1}{^{q(\varepsilon_k)}}  
\# \ccH^{\varepsilon_k}_s.  
\end{multline}
Thanks to approximations (\ref{brocueille}) and (\ref{cueille}), we next prove the following lemma. 
\begin{lem}  
\label{cucueille} We assume (\ref{recontH}). Then $\bP$-a.s.~for all $t_1 \! >\! t_0$, if for all $t\ino (t_0, t_1)$, $X_t \! > \! X_{t_0-} \! = \! X_{t_1}$, then for all $t\ino (t_0, t_1)$, $H_t \! \geq  \! 
H_{t_0} \! = \! H_{t_1}$. 
\end{lem}
\noi
\textbf{Proof.} Let $t_1 \! >\! t_0$ be such that 
for all $t\ino (t_0, t_1)$, $X_t \! > \! X_{t_0-} \! = \! X_{t_1}$. Since $X$ has only positive jumps, it implies that 
$\Delta X_{t_1}\! = \! 0$; thus, for all $s\ino [t_0, t_1]$, we get $I_{t_1}^s\! = \! X_{t_1}$. Moreover,
for all $s\ino [0, t_0)$ and for all $t\ino [t_0, t_1]$, we get 
$I_{t_1}^s\! = \! I_{t_0}^s\! = \! I_t^s$. It implies for all $t\ino [t_0, t_1]$, that $\{ I_{t_0}^s; s\ino [0, t_0) \}\backslash \{ X_{t_0-} \} \! =\!  \{ I_{t_1}^s; s\ino [0, t_1) \}\backslash \{ X_{t_1} \} \subset \{ I_{t}^s; 
s\ino [0, t) \}$ which entails the desired result when $\beta \! >\! 0$ by (\ref{brocueille}). 

 Suppose next that $\sigma_2 (\mathbf{c}) \! = \! \infty$. By a diagonal argument and (\ref{cueille}), there is a 
sequence $(\varepsilon_k)_{k\in \bbN}$ decreasing to $0$ such that $\bP$-a.s.~for all $t\ino [0, \infty )\! \cap \!  \bbQ$ and for all $s\ino [0, t]$ such that  $X_{s-} \! \leq \! I_t^s$, 
$H_s = \lim_{k \rightarrow \infty} q (\varepsilon_k)^{-1}
\#\ccH^{\varepsilon_k}_s$. First observe that for all $t\ino (t_0, t_1) \cap \bbQ$, we get $X_{t_0-} \! \leq \! I^{t_0}_t$ and that $ \ccH^{\varepsilon_k}_{t_0} \! \subset \! \ccH^{\varepsilon_k}_t$, for all $k$. Consequently, $H_{t_0} \leq H_t$, for all $t\ino (t_0, t_1) \cap \bbQ$, and thus for all $t\ino [t_0, t_1]$  since $H$ is continuous. 

  Let $t \ino (t_1, \infty) \cap \bbQ$. Let $s\ino [t_1, t]$ be such that $X_{s-}\! = \! I_t^{t_1}$. Then observe that $X_{s-} \! \leq \! I_t^s$ and that 
$\ccH^{\varepsilon_k}_{s} \! \subset \! \ccH^{\varepsilon_k}_{t_0}$ for all $k$. Consequently, $H_s \! \leq H_{t_0}$. Since $s$ can be arbitrarily close to $t_1$, the continuity of $H$ entails that $H_{t_1} \leq \! H_{t_0}$ and the previous inequality implies $H_{t_1} \! = \! H_{t_0}$, which completes the proof of the lemma. \cqfd

\subsubsection{Excursion of the height process}
\label{excuhaut} 
Let us make a preliminary remark: the results of this section are recalled from Le Gall \& 
D.~\cite{DuLG02} Chapter 1 that only deals with the critical or subcritical cases. However, they easily extend to the supercritical cases thanks to (\ref{supsubis}).

Here we assume that (\ref{recontH}) holds, which implies that
$X$ has unbounded variation sample paths. Then, basic results of fluctuation theory entail that $X\! -\! I$ is a strong Markov process in $[0, \infty)$, that $0$ is regular for 
$(0, \infty)$.  
Moreover, $\! -I$ 
is a local time at $0$ for $X\! -\! I$ (see Bertoin \cite{Be96}, Theorem VII.1). 
We denote by $\bN$ the corresponding excursion measure of $X\! -\! I$ above $0$. It is not difficult to derive from the previous approximations of $H$ that $H_t$ only depends on the excursion of $X\! -\! I$ above $0$ that straddles $t$. Moreover, the following holds true: 
\begin{equation}
\label{excusame}
\ccZ = \{ t\! \in \! [0, \infty) : H_t \! =\! 0 \}= \{ t \! \in \!  [0, \infty)  : X_t \! =\! I_t  \} \,.
\end{equation}
Since $-I$ is a local time for $X\! -\! I$ at $0$, the topological support of the Stieltjes measure $d(-I)$ is $\ccZ$. Namely, 
\begin{equation}  
\label{zero}
\textrm{$\bP$-a.s.~for all $s, t\ino [0, \infty)$ such that $s\! < \! t$,} \quad \Big( (s, t) \cap \ccZ \neq \emptyset \Big) \Longleftrightarrow \Big( I_s \! > \! I_t \Big) \,.
\end{equation}
Denote by $(a_i, b_i)$, $i\! \in \! \cI$, the connected components of the open set 
$\{ t\ino [0, \infty): H_t \! >\! 0\}$ and set $H^{i}_s \! =\!  H_{(a_i +s)\wedge b_i}$, $s\! \in \!  [0, \infty) $. We set $\zeta_i \! = \! b_i -a_i$ that is the lifetime of $H^i$. 
In the supercritical cases, there is one excursion with an infinite lifetime; more precisely, there exists $i_0\ino \cI$ such that $-I_{a_{i_0}}\! = \! \sup_{i\in \cI} (-I_{a_i} )\! = \! -I_\infty$ and $\zeta_{i_0}\! = \! 
\infty$ (recall that in the supercritical case, $-I_\infty$ is exponentially distributed with parameter $\varrho$).  
Then, the point measure 
\begin{equation}\label{PoisdecH}
 \sum_{i\in \cI} \delta_{(-I_{a_i} ,\,  H^{i})}   
\end{equation} 
is distributed as $\un_{\{ x \leq \cE \}} \cN (dxdH)$ where $\cN$ is a Poisson point measure on $\bbR_+ \! \times \! \bC([0, \infty), \bbR)$ with intensity $dx \, \bN(dH)$ and where 
$\cE\! =\!  \inf \{ x \ino [0, \infty) : \cN ([0, x] \! \times \! \{ \zeta \! = \! \infty\} ) \! \neq \! 0 \}$. Note that $\cE$ 
is exponentially distributed with parameter $\bN (\zeta\! = \! \infty)$ that is therefore equal to $\varrho$. 
Here we slightly abuse notation by denoting $\bN (dH)$ the "distribution" of $H(X)$ under 
$\bN (dX)$. In the Brownian case, up to scaling, $\bN$ is It\^o's measure of positive excursion of Brownian motion and the decomposition (\ref{PoisdecH}) corresponds to the Poisson decomposition of a reflected Brownian motion above $0$. 


As a consequence of (\ref{excusame}), $X$ and $H$ under $\bN$ have the same \textit{lifetime} $\zeta$ (that is possibly infinite in the supercritical case) that satisfies the following:
\begin{equation}\label{HNcoding}
 \textrm{$\bN$-a.e.}  \; \,  \forall t\ino [\zeta, \infty),  \; \,  X_0\! = \! H_0\! =\! X_{t}\! =\! H_t  \! =\!  0 
  \; \,  \textrm{and} \;  \,   \forall t \! \in \! (0, \zeta) , \; \,  X_{t}-I_{t}>0 \; \textrm{and} \;  H_t \! >\! 0 \,.
\end{equation} 
By (\ref{glurniglou}), in the supercritical cases, $\bN$-a.e.~on the event $\{ \zeta \! = \! \infty \}$, we get $\lim_{t\rightarrow \infty} H_t\! = \! \infty$.

Recall from (\ref{thetabdef}) the definition of the subordinator $(\gamma_x)_{x\in [0, \infty)}$ whose Laplace exponent is $\psi^{-1}$. Note that 
$\lim_{\lambda \rightarrow \infty} \psi^{-1} (\lambda)/\lambda  \! = \! 0$; consequently, 
$ (\subo_x)$ is a pure jump-process and a.s.~$\subo_x\! = \! \sum_{i\in \cI} \un_{[0, x]} (-I_{a_i}) \zeta_i $. Thus (\ref{PoisdecH}) entails 
\begin{equation}
\label{lifetimeexc}
\forall \lambda\in (0, \infty)\, ,  \quad \bN \big[ 1\! -\! e^{-\lambda \zeta}  \big]= \psi^{-1} (\lambda ), 
\end{equation}
with the convention that $e^{-\infty}\! = \! 0$. We next recall the following. 
\begin{equation}
\label{lawGamma}  
\forall a \in (0, \infty), \quad v(a):= \bN \Big(\! \!  \sup_{\;\;  t\in [0,\zeta ]}\!\! \! H_t>a \Big) \quad \textrm{satisfies} \quad \int_{v(a)}^\infty \! 
\frac{d\lambda}{\psi (\lambda)}= a \; .
\end{equation}
In the critical and subcritical cases, (\ref{lawGamma}) is proved in  Le Gall \& D.~\cite{DuLG02} (Chapter 1, Corollary 1.4.2, p.~41).  
As already mentioned, this result remains true in the supercritical cases: we leave the details to the reader.  
Note that $v \! :\!  (0, \infty) \! \rightarrow \! (\varrho, \infty)$ is a bijective decreasing $C^\infty$ function. Elementary arguments derived from (\ref{glurniglou}) and 
(\ref{lawGamma}) entail the following. 
\begin{equation}
\label{Rayxtinct}
\forall x, a\ino (0, \infty), \quad \bP \big( \! \! \! \!  \! \sup_{\quad t\in [0, \subo_x ]} \! \! \! \! \! \! H_t\leq a \big) = e^{-x v(a)} \; . 
\end{equation}

\subsection{Properties of the coloured processes}
\label{colourpro} 
Let $\alpha\ino \bbR$, $\beta \ino [0, \infty) $, $\kappa \ino (0, \infty) $ 
and $\mathbf{c}\! = \! (c_j)_{j\geq 1} \ino \elldo_3$. For all $j\! \geq \! 1$, let $(N_j (t))_{t\in [0, \infty)}$ be a homogeneous Poisson process with jump rate $\kappa c_j$; let $B$ be a standard Brownian motion with initial value $0$. We assume that the processes $B$ 
and $N_j$, $j\! \geq \! 1$,  are independent. Let $(B^{\mathtt{r}}; N^\prime_j, j\! \geq \! 1)$ be an independent copy of $(B; N_j, j\! \geq \! 1)$. Recall from (\ref{Xblue}) that for all $t\ino [0, \infty)$, we have set 
$$  X^{\mathtt{b}}_t= -\alpha t + \sqrt{\beta} B_t +  \sum_{j\geq 1}  \!\!\! \,^{\perp}\,  c_j \big( 
N_j (t) \! -\! c_j \kappa t \big) \quad \textrm{and} \quad  X^{\mathtt{r}}_t=  -\alpha t + \sqrt{\beta} B^{\mathtt{r}}_t +  \sum_{j\geq 1}  \!\!\! \,^{\perp}\,  c_j \big( N^\prime_j (t) \! -\! c_j \kappa t \big) \; , $$
where $\sum^\perp_{j\geq 1}$ stands for the sum of orthogonal $L^2$-martingales. Then $X^{\mathtt{b}}$ and $X^{\mathtt{r}}$ are two independent spectrally positive L\'evy processes whose Laplace exponent $\psi$ is defined by 
(\ref{repsidefi}). 
We assume (\ref{varinfinie}), namely: either $\beta \! >\! 0$ or $\sigma_2 (\mathbf{c})\! = \! \infty$.
We recall from (\ref{AetYdef}) the definition 
of $(A_t, Y_t)_{t\in [0, \infty)}$: 
 \begin{equation}
\label{AetYdefbis}
\forall t\ino [0, \infty) , \quad A_t = \frac{_{_1}}{^{^2}} \kappa \beta t^2 +  \sum_{j\geq 1}  c_j \big( N_j (t) \! -\! 1 \big)_+ \quad \textrm{and} \quad Y_t = X^{\mathtt{b}}_t \! -\! A_t . 
\end{equation}
Recall from (\ref{thetabdef}) the following definitions:
\begin{equation}
\label{thetabdefff}
 \forall x, t\ino [0, \infty), \quad \subo^{\mathtt{r}}_x = \inf \{ s \ino [0, \infty):    X^{\mathtt{r}}_s \! < \! -x \} \quad \textrm{and} \quad \theta_t^{\mathtt{b}}= t + \subo^{\mathtt{r}}_{A_t} ,
\end{equation}
 with the convention: $\inf \emptyset \! = \! \infty$. Recall that $\subo^{\mathtt{r}}$ is a possibly killed subordinator with Laplace exponent $\psi^{-1}$. 
%
%
%
Recall from (\ref{T*def}) that: 
\begin{equation}
\label{T*defff} 
T^*\! = \! \sup \{ t\ino [0, \infty)\! :  \theta^{\mathtt{b}}_t \! < \infty \}=   \sup \{ t\ino [0, \infty)\! : A_t \! < \! - I^{\mathtt{r}}_\infty \} \; . 
\end{equation}
In critical and subcritical cases, $T^*\! = \! \infty$ and $ \theta^{\mathtt{b}}$ only takes 
finite values. In supercritical cases, a.s.~$T^* \! < \! \infty$ and we check that $ \theta^{\mathtt{b}} (T^*-) \! < \!  \theta^{\mathtt{b}} (T^*) \! = \! \infty$. We next recall the following from (\ref{2lambda}):
\begin{equation}
\label{2lambdaaa}
\forall t \ino [0, \infty), \quad \Lambda^{\mathtt{b}}_t =  \inf \{ s \ino [0, \infty):    \theta^{\mathtt{b}}_s \! > \! t \} \quad \textrm{and} \quad \Lambda^{\mathtt{r}}_t \! = \! t-    \Lambda^{\mathtt{b}}_t . 
\end{equation} 
The processes $\Lambda^{\mathtt{b}}$ and $\Lambda^{\mathtt{r}}$ are continuous and nondecreasing. 
In critical and subcritical cases, we get a.s.~$\lim_{t\rightarrow \infty}  \Lambda^{\mathtt{b}}_t \! = \! \infty$ and $\Lambda^{\mathtt{b}} (\theta^{\mathtt{b}}_t)\! = \! t$ for all $t\ino [0, \infty)$. 
In supercritical cases, $\Lambda^{\mathtt{b}}  (\theta^{\mathtt{b}}_t)\! = \! t$ for all $t\ino [0, T^*)$ and $\Lambda^{\mathtt{b}}$ is constant and equal to $T^*$ on  
$[\theta^{\mathtt{b}} (T^*-), \infty)$.

\subsubsection{Properties of $A$}
\label{proAsec}

\noi
\textbf{Proof of Lemma \ref{AYnontriv}.} 
We assume (\ref{varinfinie}): either $\beta \! >\! 0$ or $\sigma_2 (\mathbf{c})\! = \! \infty$. 
If $\beta \! >\! 0$, then clearly a.s.~$A$ is increasing. 
Suppose that $\sigma_2 (\mathbf{c})\! = \! \infty$. With the notation of (\ref{AetYdefbis}), observe that for all $s, t\ino (0, \infty)$, 
$$ \sum_{j\geq 1} \un_{\{N_j (t) \geq 1 ; N_j (t+s) -N_j (t) \geq 1  \}} \leq \# \big\{ a \ino (t, t+s] : \Delta A_a \! >\! 0 \big\} \; .$$
Note that $\bP (N_j (t) \! \geq \! 1 ; N_j (t+s) \! -\! N_j (t)\!  \geq \! 1)\! = \! (1\! -\! \exp (-\kappa c_j t))  (1\! -\! \exp (-\kappa c_j s))$. 
Since there exists $K \ino (0, \infty)$ depending on $t$ and $s$ such that $\! (1\! -\! \exp (-\kappa c_j t))  (1\! -\! \exp (-\kappa c_j s))
\!  \geq \!  K c_j^2$ for all $j\! \geq \! 1$, Borel's Lemma implies that a.s.~$ \# \big\{ a \ino (t, t+s] : \Delta A_a \! >\! 0 \big\} \! = \! \infty$. This easily implies that $A$ is strictly increasing. To complete the proof of the lemma, observe that under (\ref{varinfinie}), $X^{\mathtt{b}}$ has infinite variation sample paths. 
Since $A$ is increasing, $Y\! = \! X^{\mathtt{b}} \! -\! A$ has infinite variation sample paths.  \cqfd

\medskip

We shall need the following estimates on $A$ in the proof of Theorem \ref{Xdefthm}.  
\begin{lem}
\label{Ainvest} Let $\alpha\ino \bbR$, $\beta \ino [0, \infty) $, $\kappa \ino (0, \infty) $ 
and $\mathbf{c}\! = \! (c_j)_{j\geq 1} \ino \elldo_3$. Assume (\ref{varinfinie}): namely, either $\beta \! >\! 0$ or $\sigma_2 (\mathbf{c})\! = \! \infty$. For all $t\ino [0, \infty)$ we set 
$A^{-1}_t \! =\!  \inf \big\{s\ino [0, \infty ) : A_s \! >\! t  \big\}$,  
that is well-defined. Then, $A^{-1}$ is continuous and there exists $a_0, a_1, a_2 \ino (0, \infty)$ that depend on $\beta, \kappa$ and $\mathbf{c}$, such that  
\begin{equation}
\label{Ainvesti}
\forall t\ino [0, \infty) , \qquad \bE \big[ A^{-1}_t\big] \leq a_1 t+ a_0 \quad \textrm{and} \quad  \bE \big[ (A^{-1}_t)^2\big] \! \leq \! a_2 t^2+a_1 t+a_0\; .
\end{equation}  
\end{lem}
\noi
\textbf{Proof.} By  Lemma \ref{AYnontriv}, $A$ is strictly increasing, which implies that $A^{-1}$ is continuous by standard arguments. 
We first suppose that $c_1 >0$. Then, by (\ref{AetYdefbis}) $A_t \!  \geq \! c_1(N_1(t)\! -\! 1)_+\! \geq \! c_1N_1(t) \! -\! c_1$. This entails that $A$ tends to $\infty$ and therefore that 
$A^{-1}$ is well-defined. 
Moreover, we get $A^{-1}_t \! \leq \! N^{-1}_1 (1+(t/c_1))$, where:  
$N^{-1}_1(t)\! = \! \inf \{ s\ino [0, \infty): N_1(s)>t\}$. 
Note that $N^{-1}_1(t)$ is the sum of $\leq \lfloor t \rfloor+1$ exponentially distributed r.v.~with parameter $\kappa c_1$, which implies that 
$\bE [N^{-1}_1(t) ]\! \leq  \! (t+1)/ (\kappa c_1)$  and $\bE [N^{-1}_1(t)^2 ]\! \leq  \!(t+1) (t+2)/ (\kappa c_1)^2$. 
Thus, 
$$ \bE \big[ A^{-1}_t \big]   \leq   \frac{_1}{^{\kappa c^2_1}} t +   \frac{_2}{^{\kappa c_1}} \quad \textrm{and} \quad \bE \big[ (A^{-1}_t)^2 \big] \leq \frac{_1}{^{\kappa^2 c_1^4}} t^2+ \frac{_5}{^{\kappa^2 c_1^3}} t + \frac{_6}{^{\kappa^2 c_1^2}}. $$
If $\mathbf{c}\!= \! 0$, then (\ref{varinfinie}) entails 
$\beta \! >\! 0$. Thus, $A^{-1}_t \! = \! \sqrt{2t/(\beta\kappa)}$ and it is then possible to choose $a_0, a_1, a_2\ino (0, \infty)$ such that (\ref{Ainvesti}) holds. \cqfd 

\subsubsection{Proof of Theorem \ref{Xdefthm}} 
\label{Xdefthmpf}

Using stochastic calculus arguments, we provide here a proof for the first statement of Theorem \ref{Xdefthm}, namely, $X$ has the same distribution as $X^{\mathtt{b}}$. See Lemma \ref{trajprop1} (i) in Section \ref{proXYsec} for the proof of the second statement. 

Let us first introduce notation. We first say that a martingale $(M_t)_{t\in [0, \infty)}$ 
is of class $(\mathscr{M})$ if 

\smallskip

\begin{compactenum}
\item[\textbf{(a)}] a.s.~$M_0\! = \! 0$, 

\smallskip

\item[\textbf{(b)}] $M$ is c\`adl\`ag,  

\smallskip

\item[\textbf{(c)}] there exists $c\ino [0, \infty)$ such that a.s.~for all $t\ino [0, \infty)$, $0\! \leq \! \Delta M_t  \! \leq \! c$,  

\smallskip

\item[\textbf{(d)}] for all $t\ino [0, \infty)$, $\bE [M_t^2] \! < \! \infty$. 
\end{compactenum}

\smallskip

\noi
Let $M$ 
be a class $(\mathscr{M})$ martingale with respect to
a 
filtration $(\ccF_{\! t})_{t\in [0, \infty)}$. Then, $\langle M \rangle $ stands for 
the predictable quadratic variation process: namely, the unique $(\ccF_t)$-predictable nondecreasing process (provided by the Doob--Meyer decomposition) such that 
$(M_t^2 \!-\! \langle M \rangle_t)_{t\in [0, \infty)}$ is a $(\ccF_{\! t})$-martingale with initial value $0$. 
We shall repeatedly use the following standard optional stopping theorem: 

\smallskip

\begin{compactenum}

\item[$(\textbf{Stp})$]Let $S$ and $T$ be two $(\ccF_{\! t})$-stopping times such that a.s.~$S\! \leq \! T \! < \! \infty$ and $\bE [\langle M\rangle_T ] \! < \! \infty$. Then, 
$\bE [M_T^2] \! = \! \bE [\langle M\rangle_T ]$ and a.s.~$M_S \! 
= \! \bE [M_T | \ccF_S]$. 

\end{compactenum} 
 
\smallskip

\noi
Then, the \textit{characteristic measure} of $M$ is a random measure $\mathcal{V}$ on $[0, \infty)\! \times \! (0, \infty)$ such that:
\begin{compactenum}

\smallskip

\item[$\bullet$] for all $\epp\ino (0, \infty)$, the process $t \! \mapsto \!  \mathcal{V} \big( 
[0, t] \! \times \! [\epp , \infty) \big)$ is $(\ccF_{\! t})$-predictable;

\smallskip

\item[$\bullet$] for all $\epp\ino (0, \infty)$, $t\longmapsto  \sum_{s\in [0, t]} \un_{[\epp, \infty)}( \Delta M_s) - \mathcal{V} \big( 
[0, t] \! \times \! [\epp , \infty) \big)\, $ is a $(\ccF_{\! t})$-martingale. 
\end{compactenum} 

\smallskip

\noi
(See Jacod \& Shiryaev \cite{JaSh02}, Chapter II, Theorem 2.21, p.~80.) 
The \textit{purely discontinuous part} of $M$ is obtained as the $\mathcal{V}$-compensated sum of its jumps: namely, the $L^2$-limit as $\varepsilon$ goes to $0$ of 
the martingales $t\longmapsto \sum_{s\in [0, t]}\!  \Delta M_s\, \un_{[\epp, \infty)}  ( \Delta M_s) \! -\!  \int_{[0, t] \times [\epp , \infty)} \! r\,  \mathcal{V} (dsdr)$. 
The purely discontinuous part of $M$ is denoted by $M^d$ and it is a $(\ccF_{\! t})$-martingale of class $(\mathscr{M})$. 
The \textit{continuous part} of $M$ is the continuous $(\ccF_{\! t})$-martingale $M^c\! =\! M \! -\! M^d$.  Note that $M^c$ is also a $(\ccF_{\! t})$-martingale of class 
$(\mathscr{M})$.
We call $(\langle M^c \rangle, \mathcal{V})$ the \textit{characteristics} of $M$. 
For more detail on the canonical representation of semi-martingales, see Jacod \& Shiryaev \cite{JaSh02} Chapter II, Definition 2.16 p.~76 and \S 2.d, Theorem 2.34, p.~84.

\smallskip

Let $(\ccF^N_t)_{t\in [0, \infty)}$ (resp.~$(\ccF^B_t)_{t\in [0, \infty)}$) be the right-continuous filtration associated with the natural filtration of the process $(N_j (\cdot))_{j\geq 1}$ (resp.~$B$); 
then, we set $\ccF^0_t \! =\!  \sigma (\ccF^N_t , \ccF^B_t)$, $t\ino [0, \infty)$, and 
$$ \forall t\ino [0, \infty), \quad X^{*\mathtt{b}}_t=X^{\mathtt{b}}_t + \alpha t = \sqrt{\beta} B_t +  \sum_{j\geq 1}  \!\!\! \,^{\perp}\,  c_j \big( 
N_j (t) \! -\! c_j \kappa t \big)\; .$$
By standard arguments on L\'evy processes, $X^{*\mathtt{b}}$ is a $(\ccF^0_t)$-martingale. 
We set $a_3 \! = \! \beta + \kappa \sigma_3 (\mathbf{c} )$ and we easily check that 
\begin{equation}
\label{globiglou} \textrm{$t\! \longmapsto \! (X^{*\mathtt{b}}_t)^2 \! -\! a_3 t\; $ is a 
$(\ccF^0_t)$-martingale.}  
\end{equation}
Moreover, we easily check that $X^{*\mathtt{b}}$ is in the class $(\mathscr{M})$ and that 
its (deterministic) characteristics are the following: its characteristic measure is $dt \otimes 
\pi (dr)$, where $\pi (dr)\! =\!  \sum_{j\geq 1} \kappa c_j \delta_{c_j}  $; its continuous part is  $\sqrt{\! \beta} B$, whose predictable quadratic variation process is $t\mapsto \beta t$. 
To prove Theorem \ref{Xdefthm}, we shall use the converse of this result: namely, 
\textit{a martingale whose characteristics are $dt \otimes \pi (dr)$ and $t\mapsto \beta t$ has necessarily the same law as $X^{*\mathtt{b}}$} (for a proof see Jacod \& Shiryaev \cite{JaSh02} Chapter II, \S 4.c, Corollary 4.18, p.~107). To that end, one computes the characteristics of several time-changes of $X^{*\mathtt{b}}$ and of $X^{\mathtt{r}}$.  

  First, recall from Lemma \ref{Ainvest} that $A^{-1}$ is continuous and note that $A^{-1}_t$ is a $(\ccF^0_r)$-stopping time. 
We set 
$$ \forall t\in [0, \infty), \quad M^{_{(1)}}_t \! =\!  X^{*\mathtt{b}}(A^{-1} (t)) \quad \textrm{and} \quad  
\ccF^1_{\! t} \! =\!  \ccF^0(A^{-1}_t)\; .$$
By (\ref{globiglou}), $\langle X^{*\mathtt{b}} \rangle_t \! =\!  a_3 t$ and 
(\ref{Ainvesti}) combined with $(\textbf{Stp})$ imply that $M^{_{(1)}}_{^{\, \!}}$ is a square integrable $(\ccF^1_t)$-martingale and that 
$\bE [(M^{_{(1)}}_{^{t}})^2]\! = \! a_3 \bE [A^{-1}_t]$. Then, set 
$g(r)\! =\! \inf \{ s\ino [0, \infty) : A^{-1}_s\! = \! r \}$, for all $r\ino [0, \infty)$, so that $g(r)=A_{r-}$. Since $A^{-1}$ is continuous, we see that  
 \begin{equation}
\label{corrjumpi}
  g : \big\{ r\ino [0, A^{-1}_t]: \Delta X^{*\mathtt{b}}_r \! >\! 0 \big\}  \longrightarrow \big\{ s\ino [0, t] : \Delta M^{{(1)}}_s >0 \big\}  \; \, \textrm{is one-to-one.}
\end{equation}
Moreover, if $s$ and $r$ are such that  $g(r)\! = \! s$ and $\Delta M^{_{(1)}}_s \! >\! 0$, then $\Delta M^{_{(1)}}_s\! = \! \Delta X^{*\mathtt{b}}_r$. This implies that $M^{_{(1)}}_{^{\, \! }}$ is of class $(\mathscr{M})$. 

For all $\epp \ino (0, \infty)$, we next set 
$$\forall r\ino [0, \infty), \quad 
 J^\epp_r \! = \! \sum_{r^\prime\in [0, r]} \!\!\!  \un_{[\epp, \infty)}( \Delta X^{*\mathtt{b}}_{r^\prime}) \; - r\pi ([\epp, \infty)) $$ 
 that is a $(\ccF^0_r)$-martingale of class $(\mathscr{M})$ such that $\langle J^\epp\rangle_r\! = \!   \pi ([\epp, \infty)) r$; (\ref{Ainvesti}) combined with $(\textbf{Stp})$ entails that  
$J^\epp    \circ   A^{-1}$ is a square integrable $(\ccF^1_t)$-martingale. Moreover, (\ref{corrjumpi}) entails that 
$J^\epp (A^{-1}_t)\! = \! \sum_{s\in [0, t]} \un_{[\epp, \infty)}( \Delta M^{_{(1)}}_s) \! -\!  A^{-1}_t \pi ([\epp, \infty))$. Since $A^{-1}$ continuous, it is  $(\ccF^1_t)$-predictable and 
$dA^{-1}_t \! \otimes \! \pi (dr)$ 
is the characteristic measure of $M^{_{(1)}}_{^{\, \!}}$. It easily entails that the continuous part of $M^{_{(1)}}_{^{\, \!}}$ is $ \sqrt{\beta}B \circ A^{-1}$. We next set $Q_t \! =\!  \beta B^2_t \! -\! \beta t$, 
which is a martingale. We intend to apply $(\textbf{Stp})$ to show that $Q\circ A^{-1}$ is a martingale. Note that by It\^o's formula $\langle Q \rangle_t \! = \! 4\beta^2 \int_0^t B^2_s ds$ and thus, $\bE [ \langle Q \rangle_t ]\! = \! 2\beta^2t^2$. Since $A^{-1}$ is independent of $B$, 
$\bE [ \langle Q \rangle (A^{-1}_t)] \! = \! 2\beta^2 \bE [(A^{-1}_t)^2]$ that is a finite 
quantity by (\ref{Ainvesti}). Then, by $(\textbf{Stp})$, we see that $\langle \sqrt{\beta} B\! \circ \! A^{-1} \rangle \! =  \! \beta A^{-1} $.  
We have proved that 
$\beta A^{-1}$ and $dA^{-1}_t \! \otimes \! \pi (dr)$ are the characteristics of $M^{_{(1)}}_{^{\, \!}}$. It is easy to realize that $M^{_{(1)}}_{^{\, \!}}$ is also a martingale with respect to the natural filtration of $(A^{-1}, M^{_{(1)}}_{^{\, \!}})$ with the same 
characteristics $\beta A^{-1}$ and $dA^{-1}_t \! \otimes \! \pi (dr)$since $A^{-1}$ is continuous. 

\begin{rem}
\label{contispeci}Note that in the previous arguments, the continuity of $A^{-1}$ plays a key role. Since this property does not hold in the discrete cases, the above arguments cannot be adapted to give a proof in such cases. 
\cq 
\end{rem}

We next prove the following lemma. 
\begin{lem}
\label{tpschgind} Let $E$ be a Polish space and let $(Z_t)_{t\in [0,\infty)}$ be a $E$-valued c\`adl\`ag process. Let $(M_r)_{r\in [0, \infty)}$ be a c\`adl\`ag martingale with respect to
the natural filtration of $Z$. Let 
$(\phi_t)_{t\in[0, \infty)}$ be a nondecreasing c\`ad process that is adapted to a filtration $(\ccG_t)_{t\in [0, \infty)}$. We assume that $Z$ and $\ccG_\infty$ are independent and that for all 
$t \ino [0, \infty)$, $\int \! \bP (\phi_t\ino dr) \bE [|M_r|] \! < \! \infty$. We set $\ccF_t \!  = \! \sigma (Z_{\cdot \wedge \phi_t} , \ccG_t)$, for all $t\ino [0, \infty)$. Then, $M\! \circ \! \phi$ is a c\`adl\`ag  
$(\ccF_t)$-martingale.  
\end{lem}
\noi
\textbf{Proof:} Let $t, r_1, \ldots , r_n \ino [0, \infty)$ and let $s\ino [0, t]$. Let $G: E^n \! \rightarrow \! [0, \infty)$ be bounded and measurable and let $Q$ be a nonnegative bounded $\ccG_s$-measurable random variable. We get:
\begin{eqnarray*}
\bE \big[ M_{\phi_t}Q \, G\big( (Z_{r_k\wedge \phi_s})_{1\leq k\leq n} \big) \big]\!\!  &= &\!\!\!\! 
\int \bP (\phi_t \ino d r^\prime ; \phi_s \ino dr ; Q \ino dq)\,  q \, \bE \big[ M_{r^\prime} G\big( (Z_{r_k\wedge r})_{1\leq k\leq n} \big) \big]     \\
& =& \!\!\!\!    \int \bP (\phi_t \ino d r^\prime ; \phi_s \ino dr ; Q \ino dq) \, q \, \bE \big[ M_{r} G\big( (Z_{r_k\wedge r})_{1\leq k\leq n} \big) \big] \\
&= & \!\!\!\!  \bE \big[ M_{\phi_s}Q \, G\big( (Z_{r_k\wedge \phi_s})_{1\leq k\leq n} \big) \big] 
\end{eqnarray*}
and we completes the proof by use of the monotone class theorem. \cqfd 

\smallskip

Recall from the beginning of Section \ref{colourpro} the definitions of the processes $B^{\mathtt{r}}$, $(N^\prime_j (\cdot))_{j\geq 1}$ and $X^{\mathtt{r}}$; recall that $X^{\mathtt{r}}$ is an independent copy of $X^{\mathtt{b}}$. We next need the following result. For all $t\ino [0, \infty)$, we set $I^{\mathtt{r}}_t\! = \! \inf_{s\in [0, t] } X^{\mathtt{r}}_s$. Then, 
\begin{equation}
\label{martini} \forall p, t\ino (0, \infty), \quad \bE [(- I^{\mathtt{r}}_t)^p ] \! < \! \infty \, . 
\end{equation}
\noi
\textit{Indeed,} recall that $\gamma^\mathtt{r}_x\! = \! \inf \{ t \ino [0, \infty): X^{\mathtt{r}}_t \! <\!  -x \}$ and that $x \! \mapsto \! \gamma^\mathtt{r}_x$ is a (possibly killed) subordinator with Laplace exponent $\psi^{-1}$. Then for all $\lambda \ino (0, \infty)$ we get the following 
\begin{eqnarray*} \bE [(- I^{\mathtt{r}}_t)^p ] \!\! & = & \!\!\!\!  \int_0^\infty \!\!\!\! px^{p-1} \bP ( -I^{\mathtt{r}}_t \! >\!  x) \, dx \, \leq \int_0^\infty \!\!\!\! px^{p-1} \bP ( \gamma^\mathtt{r}_x \leq  t) dx \\
& \leq &\!\!\!\!   \int_0^\infty \!\!\!\! px^{p-1}  e^{\lambda t } \bE [e^{-\lambda  \gamma^\mathtt{r}_x }]  \, dx \, = p \Gamma (p) (\psi^{-1} (\lambda))^{-p} e^{\lambda t} \; , 
\end{eqnarray*}
which entails (\ref{martini}). \cq 

\smallskip

We apply Lemma \ref{tpschgind} to $Z\!=\!  (  A^{-1}\!, M^{_{(1)}}_{^{\, \!}})$, to $\phi_t\! = \! -I^{\mathtt{r}}_t$, to $(\ccG_t)$ that is taken as 
the right-continuous filtration associated with the natural filtration of $(B^{\mathtt{r}}; N^\prime_j, j\! \geq \! 1)$, and first to $M\! = \! M^{_{(1)}}_{^{\, \!}}\!\!$ and next to $M\! = \! J^\epp \! \circ \! A^{-1}$.
Recall that $\bE \big[ (M^{_{(1)}}_{^{t}})^2 \big]\! = \! a_3 \bE [A^{-1}_t] \! \leq  \! a_3(a_1t + a_0)$ and 
$\bE \big[ J^\epp (A^{-1}_t)^2\big] \! = \! \pi ([\epp , \infty))   \bE[ A^{-1}_t] \! \leq  \! \pi ([\epp , \infty)) (a_1t + a_0)$, 
by (\ref{Ainvesti}). In both cases ($M\! = \! M^{_{(1)}}_{^{\, \!}}\!\!$ or $M\! = \! J^\epp \! \circ \! A^{-1}$), we get $\int \! \bP (-I^{\mathtt{r}}_t \ino dr) \bE [M^2_r]\! < \! \infty$, by (\ref{martini}). 
Then, we set for all $t \ino [0, \infty)$, 
$$M^{{(2)}}_{t} \! = \! M^{{(1)}}_{-I^{\mathtt{r}}_t} , \quad  J^{\prime \epp}_t \! = \! J^{\epp}_{A^{-1}( -I^{\mathtt{r}}_t)}\quad  \textrm{and} \quad  
\ccF^2_t\! = \! \sigma \big( \ccG_t ,  A^{-1}_{\cdot \wedge (-I^{\mathtt{r}}_t)}, M^{{(1)}}_{{\cdot \wedge (-I^{\mathtt{r}}_t)}} \big)\; .$$ Lemma \ref{tpschgind} asserts that $M^{_{(2)}}_{^{\, \!}}$ and $J^{\prime \epp}$ are $(\ccF^{2}_t)$-martingales and we proved just above that they are square integrable.  
Since they are c\`adl\`ag processes standard arguments entail they are also $(\ccF^{2}_{t+})$-martingales.

We next set $I^\mathtt{r}_\infty\! = \! \lim_{t\rightarrow \infty} I^\mathtt{r}_t$ (that is a.s.~infinite if $\alpha \! \geq \! 0$ and that is a.s.~finite if $\alpha \! < \! 0$).  
For all $r\ino [0, -I^\mathtt{r}_\infty)$, we next set 
$g^\prime (r)\! = \! \inf \{ s\ino [0, \infty) \! : \! -I^{\mathtt{r}}_s \! = \! r \}$ (note that $g^\prime (r)\! = \! \subo^\mathtt{r}_{r-}$). 
Since $I^\mathtt{r}$ is continuous, it is easy to check that a.s.~for all $t\ino [0, \infty)$, $g^\prime$ 
is a one-to-one correspondence between $\{r \ino [0,-I^{\mathtt{r}}_t ] \! : \! \Delta M^{_{(1)}}_{{r}} \!>\! 0\}$ and $\{ s\ino [0, t] \! :\!  \Delta M^{_{(2)}}_{{s}} \!>\! 0\}$. Moreover, if $s$ and $r$ are such that 
$\Delta M^{_{(2)}}_s \! >\! 0$ and $g^\prime(r)\! = \! s$, then $\Delta M^{_{(2)}}_s\! = \! \Delta M^{_{(1)}}_r$. 
This first entails that $M^{_{(2)}}_{^{\, \!}}$ is of class 
$(\mathscr{M})$. It also implies 
that 
\begin{equation}
\label{flapiyopi}
J^{\prime \epp}_t \! = \! \sum_{s\in [0, t]} \un_{[\epp , \infty)} (\Delta M^{_{(2)}}_{^t}) \! -\! A^{-1} (-I^{\mathtt{r}}_t) \pi ([\epp , \infty)) \; . 
\end{equation} 
Since $t \! \mapsto \! A^{-1} (-I^{\mathtt{r}}_t)$ is continuous and $(\ccF^{2}_{t+})$-adapted it is 
$(\ccF^{2}_{t+})$-predictable; therefore, the characteristic measure of $M^{_{(2)}}_{^{\, \!}}$ (with respect to the filtration $(\ccF^{2}_{t+})$) is $ d(A^{-1} \!  \circ \!  (-I^{\mathtt{r}}))(t)  \otimes  \pi (dr) $.

It is easy to deduce that the continuous part of  $M^{_{(2)}}_{^{\, \!}}$ is $\sqrt{\beta} B\!  \circ\!  A^{-1} \! \circ \! (-I^{\mathtt{r}}) $. We then apply Lemma \ref{tpschgind} to $M \! = \! (B \circ A^{-1})^2 -A^{-1}$: to that end note that 
$\bE [|M_t| ] \!\leq \! 2 \bE [A^{-1}_t]\! \leq \! 2(a_1 t+a_0)$, by (\ref{Ainvesti}); by (\ref{martini}), we get 
$\int \! \bP (-I^{\mathtt{r}}_t \ino dr) \bE [ |M_r|]\! < \! \infty$; thus, Lemma \ref{tpschgind} applies and asserts that $M\! \circ \! (-I^{\mathtt{r}}) $ is a 
$(\ccF^2_t)$-martigale; by standard arguments, it is also a martingale with respect to
$(\ccF^2_{t+})$. This entails that $\beta A^{-1} \! \circ  (-I^{\mathtt{r}})$ is the quadratic variation of 
$\sqrt{\beta} B\!  \circ\!  A^{-1} \! \circ \! (-I^{\mathtt{r}}) $ 
that is the continuous part of $M^{_{(2)}}_{^{\, \!}}$. Thus, $\beta  A^{-1} \! \circ \! (-I^{\mathtt{r}}) $ and $ d(A^{-1} \!  \circ \!  (-I^{\mathtt{r}}))(t)  \otimes  \pi (dr) $ are the $(\ccF^2_{t+})$-characteristics of  $M^{_{(2)}}_{^{\, \!}}$.

Recall from (\ref{thetabdefff}) and (\ref{2lambdaaa}) the definitions of $\theta^{\mathtt{b}}$, $\Lambda^{\mathtt{b}}$ and $\Lambda^{\mathtt{r}}$. Then, we next check that a.s.~
\begin{equation}
\label{bloodymary}
\forall t\ino [0, \infty) , \quad  \theta^{\mathtt{r}}_t := \inf \big\{ s\ino [0, \infty ): \Lambda^{\mathtt{r}}_s >t  \big\} = t + A^{-1} (-I^{\mathtt{r}}_t) \; .
\end{equation}
\noi
\textit{Indeed}, since $A$ is strictly increasing, 
$z\! < \! A^{-1} (-I^{\mathtt{r}}_t)$ implies $A_z \! <\!  -I^{\mathtt{r}}_t$
(resp.~$z\! > \! A^{-1} (-I^{\mathtt{r}}_t)$ implies $A_z \! >\!  -I^{\mathtt{r}}_t$),  
which then yields $ \gamma^{\mathtt{r}} (A_z) \! <\! t$ (resp.~$\gamma^{\mathtt{r}} (A_z) \! >\! t$). Since 
$\theta^{\mathtt{b}}_z \! -\! z \! = \! \gamma^{\mathtt{r}} (A_z)$, we get $ \theta^{\mathtt{b}}_z \! <\! t+ z \! < \! t + A^{-1} (-I^{\mathtt{r}}_t)$ (resp.~$ \theta^{\mathtt{b}}_z \! >\! t+ z \! > \! t + A^{-1} (-I^{\mathtt{r}}_t)$). 
Consequently, $z\! = \! \Lambda^{\mathtt{b}}  ( \theta^{\mathtt{b}}_z ) \! \leq \! \Lambda^{\mathtt{b}} ( t+ A^{-1} (-I^{\mathtt{r}}_t)) $ (resp.~$z\! = \! \Lambda^{\mathtt{b}}  ( \theta^{\mathtt{b}}_z ) \! \geq \! \Lambda^{\mathtt{b}} ( t+ A^{-1} (-I^{\mathtt{r}}_t)) $). We then easily get 
$A^{-1} (-I^{\mathtt{r}}_t) \! = \! \Lambda^{\mathtt{b}} ( t+ A^{-1} (-I^{\mathtt{r}}_t))$. 
By (\ref{2lambdaaa}), we get $\Lambda^{\mathtt{r}} ( t+ A^{-1} (-I^{\mathtt{r}}_t)) \! =\! t$, which completes the proof of (\ref{bloodymary}).  \cq 

\smallskip

We next set for all $t\ino [0, \infty)$, 
$$ X^{*\mathtt{r}}_t= X^{\mathtt{r}}_t + \alpha t \quad \textrm{and} \quad M_t \! =\! X^{*\mathtt{r}}_t +  M^{_{(2)}}_{^{t}} \; .$$
Clearly, $M$ is a $(\ccF^{2}_{t+})$-martingale (as the sum of two such martingales) that belongs to  
the class $(\mathscr{M})$. Note that a.s.~$X^{*\mathtt{r}}$ and $M^{_{(2)}}_{^{\! \, }}$ do not jump simultaneously. Thus, by 
(\ref{flapiyopi}) and (\ref{bloodymary}), we get a.s.
$$ J^{\prime \prime \epp}_t =  J^{\prime \epp}_t + \sum_{s\in [0, t]} \un_{[\epp , \infty)} (\Delta X^{*\mathtt{r}}_s) - t\pi ([\epp , \infty))  = \sum_{s\in [0, t]} \un_{[\epp, \infty)}( \Delta M_s) - 
\theta^{\mathtt{r}}_t  \pi ([\epp , \infty)) \; .$$
Moreover, $J^{\prime \prime \epp}$ is clearly a $(\ccF^{2}_{t+})$-martingale 
(as the sum of two such martingales). Next,  
(\ref{bloodymary}) implies that $\theta^{\mathtt{r}}$ is continuous, strictly increasing and $(\ccF^{2}_{t+})$-adapted; therefore 
it is $(\ccF^{2}_{t+})$-predictable.  
This implies that the characteristic measure of $M$ is $d \theta^{\mathtt{r}}_t \otimes \pi (dr)$. 

Consequently, the continuous component of $M$ is $M^c_t\!= \!   \sqrt{\beta} (B^{\mathtt{r}}_t + B(A^{-1} (-I^{\mathtt{r}}_t)))$. The independence of $B$, $A^{-1}\! $ and $B^{\mathtt{r}}$ easily entails that 
$t\!  \mapsto \! B^{\mathtt{r}}_t +B(A^{-1} (-I^{\mathtt{r}}_t)))$ is a $(\ccF^{2}_{t+})$-martingale. Moreover, recall that the predictable quadratic variation of $\sqrt{\beta} B\!  \circ\!  A^{-1} \! \circ \! (-I^{\mathtt{r}}) $ is equal to 
$\beta .A^{-1} \! \circ \! (-I^{\mathtt{r}})$. Thus, it is easy to see that the predictable quadratic variation of $M^c$ is equal to $\beta  \theta^{\mathtt{r}}$. We then have proved that the characteristics of $M$ are 
$\beta  \theta^{\mathtt{r}}$ and $d \theta^{\mathtt{r}}_t \otimes \pi (dr)$. 

  We next recall from (\ref{bloodymary}) that $\Lambda^{\mathtt{r}}$ is the inverse of $\theta^{\mathtt{r}}$ that is also strictly increasing and continuous. We set $X^* \! = \! M\!  \circ \! \Lambda^{\mathtt{r}}$ and we 
see that a.s.~for all $t\ino [0, \infty)$, 
\begin{eqnarray}
\label{manhattan}
X^*_t &= &X^{*\mathtt{r}} (\Lambda^{\mathtt{r}}_t)+ X^{*\mathtt{b}} (A^{-1} (-I^{\mathtt{r}} (\Lambda^{\mathtt{r}}_t))) \nonumber  \\
&=&   X^{*\mathtt{r}} (\Lambda^{\mathtt{r}}_t)+X^{*\mathtt{b}} (\Lambda^{\mathtt{b}}_t) \nonumber \\
&=&   X^{\mathtt{r}} (\Lambda^{\mathtt{r}}_t)+X^{\mathtt{b}} (\Lambda^{\mathtt{b}}_t) + 
\alpha t \; .
\end{eqnarray}
\noi
\textit{Indeed}, the first equality is a direct consequence of the definition; then recall from (\ref{bloodymary}) that $A^{-1} (-I^{\mathtt{r}}_t)\! = \!  \theta^{\mathtt{r}}_t \! -\! t$, thus, $A^{-1} (-I^{\mathtt{r}} (\Lambda^{\mathtt{r}}_t))\! = \! t\! -\! \Lambda^{\mathtt{r}}_t \! =\!   \Lambda^{\mathtt{b}}_t$, which entails the second equality and also (\ref{manhattan}).  \cq 

\smallskip

Observe that for all $t\ino [0, \infty)$, $\Lambda^{\mathtt{r}}_t$ is a $(\ccF^{2}_{r+})$-stopping time such that $\Lambda^{\mathtt{r}}_t\! \leq \! t$. We then set 
$\ccF_t\! = \! \ccF^2(\Lambda^{\mathtt{r}}_t+)$. 
The optional stopping theorem applies to $M$ and $J^{\prime \prime \epp}$ to show that $X$ and $J^{\prime \prime \epp} \! \circ \! \Lambda^{\mathtt{r}}$ are $(\ccF_t)$-square integrable martingales. Since $\Lambda^{\mathtt{r}}$ is strictly increasing and continuous, $X$ is of class $(\mathscr{M})$ and 
$J^{\prime \prime \epp} (\Lambda^{\mathtt{r}}_t)\! = \! \sum_{s\in [0, t]} \un_{[\epp, \infty)}( \Delta X^*_s) - t  \pi ([\epp , \infty))$. This proves that $dt\otimes \pi(dr)$ is the characteristic measure of $X^*$. Consequently, 
$M^{c} \! \circ \! \Lambda^{\mathtt{r}}$ is the continuous part of $X^*$. Since $\Lambda^{\mathtt{r}}$ is a bounded stopping-time, the optional stopping theorem applies to the martingale 
$(M^c)^2 \! -\!  \beta\theta^{\mathtt{r}}$ and it entails that  $\langle M^{c} \! \circ \! \Lambda^{\mathtt{r}} \rangle_t\! = \! \beta t$. Thus, the characteristics of $X^*$ are $t \! \mapsto \! \beta t$ and $dt \! \otimes \! \pi(dr)$. 
By \cite{JaSh02} Corollary 4.18 in Jacod \& Shiryaev \cite{JaSh02}  (Chapter II, \S 4.c, p.~107), it implies that $X^*$ has the same law as $X^{*\mathtt{b}}$, which completes the proof of Theorem \ref{Xdefthm} by  (\ref{manhattan}). \cqfd 

\subsubsection{Properties of $X$ and $Y$} 
\label{proXYsec}
Recall from the beginning of Section \ref{colourpro} the definition of the processes $X^\mathtt{b}$, $A$, $Y$, $X^\mathtt{r}$, $\subo^\mathtt{r}$, $\theta^\mathtt{b}$, $\Lambda^\mathtt{b}$ and 
$\Lambda^\mathtt{r}$. Recall from (\ref{T*defff}) the definition of $T^*$. In supercritical cases, also recall that 
a.s.~$T^* \! < \! \infty$ and $ \theta^{\mathtt{b}} (T^*-) \! < \!  \theta^{\mathtt{b}} (T^*) \! = \! \infty$.  
Recall from (\ref{Xdef}) that 
$X_t \! = \! X^{\mathtt{b}} (\Lambda^{\! \mathtt{b}}_t ) + X^{\mathtt{r}} (\Lambda^{\! \mathtt{r}}_t )$ for all $t\ino [0, \infty)$. 
Let us mention that the proof of the following lemma does not 
use the fact that $X$ is a L\'evy process, and note that, \emph{(i)} completes the proof of Theorem \ref{Xdefthm} of the previous section.
\begin{lem}
\label{trajprop1} Let $\alpha \ino \bbR$, $\beta \ino [0, \infty) $, $\kappa \ino (0, \infty) $ 
and $\mathbf{c}\! = \! (c_j)_{j\geq 1} \ino \elldo_3$ satisfy (\ref{varinfinie}). Namely: either $\beta \! >\! 0$ or $\sigma_2 (\mathbf{c})\! = \! \infty$.  
Then, the following holds true. 
\begin{compactenum}
\smallskip
\item[$(i)$] $\bP$-a.s.~for all $a\ino [0, T^*)$, $X_{\theta^{\mathtt{b}}_a}\! = \! Y_a $. 
\smallskip
\item[$(ii)$] $\bP$-a.s.~for all $a\ino [0, T^*)$, if $\Delta \theta^{\mathtt{b}}_a \! =\! 0$, then $t\! = \! \theta^\mathtt{b}_a$ is the unique $t \ino [0, \infty)$ such that $\Lambda^{\mathtt{b}}_t \! = \! a$. 
\smallskip
\item[$(iii)$]  $\bP$-a.s.~for all $a\ino [0, T^*]$, if $\Delta \theta^{\mathtt{b}}_a \! >\! 0$, then $\Delta X \! (\theta^{\mathtt{b}}_{a-}) \! = \! \Delta A_a $ and 
$\Delta Y_a\! = \! 0$. Moreover, 
$$ \forall t\ino \big( \theta^{\mathtt{b}}_{a-}, \theta^{\mathtt{b}}_{a} \big), \quad X_t  \! \geq \! X_{t-} \!  > \! 
X_{(\theta^{\mathtt{b}}_{a-})-}\!\! =\! Y_a  \quad \textrm{and if $a<T^*\! $, then }\quad  X_{(\theta^{\mathtt{b}}_{a-})-}\!\! = \! X_{\theta^{\mathtt{b}}_a} .$$ 
\end{compactenum}
\end{lem}
\noi
\textbf{Proof.} Suppose that $a\! < \! T^*$. Thus, $\theta^{\mathtt{b}}_a\! < \! \infty$, by definition of $T^*$. Then observe that $\Lambda^{\! \mathtt{r}}({\theta^{\mathtt{b}}_a}) \! = \! \theta^{\mathtt{b}}_a \! -\!  
\Lambda^{\! \mathtt{b}} ({\theta^{\mathtt{b}}_a}) \! = \! \theta^{\mathtt{b}}_a \! -\! a \!  = \! \subo^{\mathtt{r}} ({A_a})$ and thus, 
$X ({\theta^{\mathtt{b}}_a}) \! = \! X^\mathtt{b}_a + X^{\mathtt{r}} ({\subo^{\mathtt{r}} ({A_a}})) \!\!  = \!  X^\mathtt{b}_a\! - \! A_a\! = \! Y_a$, which proves $(i)$. 

Then, $(ii)$ is an immediate consequence of the definition: indeed 
since $\Lambda^{\! \mathtt{b}}$ is the pseudo-inverse of $\theta^\mathtt{b}
$, if $\Lambda^\mathtt{b}_t\! = \! a \! < \! T^* \! $, then $ \theta^\mathtt{b}_{a-} \! \leq \! t \! \leq \!  \theta^{\mathtt{b}}_a $. 

Let us prove $(iii)$: we fix $a\ino [0, T^*]$ and we suppose that $\Delta \theta^{\mathtt{b}}_a \! >\! 0$. 
Observe that for all $z\ino [0, a)$, 
$\theta^{\mathtt{b}}_z \! <\!  \theta^{\mathtt{b}}_{a-}\! < \! \infty$ and $Y_{a-}\! = \! 
\lim_{b\rightarrow a-} Y_z \! = \!   \lim_{b\rightarrow a-} X(\theta^\mathtt{b}_z) \! =\! X(\theta^\mathtt{b}_{a-}-)$ by $(i)$ and since $\theta^\mathtt{b}$ increases strictly.

We first assume that $\Delta A_a >0$. Since $Y$ and $A$ have distinct jump-times, we get $\Delta Y_a\! = \! 0$ and 
$Y_a\! = \! Y_{a-}\! = \! X(\theta^\mathtt{b}_{a-}-)$. 
Since the processes $A$ and $X^\mathtt{r}$ are independent, it is easy to check that a.s.~$\{x\ino [0, \infty) : \Delta \subo^\mathtt{r}_x \! >\! 0\} \cap \{ A_{a-} ; a\ino [0, \infty) \! : \! \Delta A_a \! >\! 0\} \! = \! \emptyset$. Thus, $\theta^{\mathtt{b}}_{a-}\! = \! a+ \subo^\mathtt{r} (A_{a-} )$ and for all $t\ino [\theta^{\mathtt{b}}_{a-}, \theta^{\mathtt{b}}_{a} )$, $\Lambda^\mathtt{b}_t \! = \! a$ and $\Lambda^\mathtt{r}_t\! = \! t\! -\! a\! = \! t\! -\! \theta^{\mathtt{b}}_{a-} +  \subo^\mathtt{r} (A_{a-} )$. Thus, for all $s\ino [0, \Delta \theta^\mathtt{b}_a )$, 
\begin{equation}
\label{gadafla}
X_{s + \theta^{\mathtt{b}}_{a-} } \! = \! X^\mathtt{b}_a + X^\mathtt{r}_{s+  \subo^\mathtt{r} (A_{a-} )} \! = \! Y_a + A_a +  
X^\mathtt{r}_{s+  \subo^\mathtt{r} (A_{a-} )} .
\end{equation} 
Taking $s\! = \! 0$ in the previous equality first entails 
$X({ \theta^{\mathtt{b}}_{a-}} )\! = \! Y_a + \Delta A_a$. Recall that $Y_a\! = \! Y_{a-}\! = \! X(\theta^\mathtt{b}_{a-}-)$. Thus, 
$\Delta X \! (\theta^{\mathtt{b}}_{a-}) \! = \! \Delta A_a $. 
Moreover, for all $ s\ino (0, \Delta \theta^\mathtt{b}_a )$, $A_a+ X^\mathtt{r} ((s+  \subo^\mathtt{r} (A_{a-} ))-)\!\!  >\! 0$. Thus, (\ref{gadafla}) entails that for all $ t\ino \big( \theta^{\mathtt{b}}_{a-}, \theta^{\mathtt{b}}_{a} \big)$, we get $X_t  \! \geq \! X_{t-} \!  >  \! Y_a  $. 
Furthermore, note that if $a\! < \! T^*$, namely if $\theta^\mathtt{b}_a \! < \! \infty$, then 
by $(i)$, we see that 
$X(\theta^\mathtt{b}_a)\! = \! Y_a$. This proves $(iii)$ when $\Delta A_a \! >\! 0$. 

  We next assume that $\Delta \theta^{\mathtt{b}}_a \! >\! 0$ and $\Delta A_a \! = \! 0$ (which occurs when $\beta \! >\! 0$). Consequently, 
 $\theta^{\mathtt{b}}_{a-}\! = \! a+ \subo^\mathtt{r} ((A_{a})- )$. Since 
$(Y,A)$ and $X^\mathtt{r}$ are independent, we a.s.~get $\{x\ino [0, \infty) : \Delta \subo^\mathtt{r}_x \! >\! 0\} \cap \{ A_{a} ; a\ino [0, \infty) \! : \! \Delta Y_a \! >\! 0\} \! = \! \emptyset$. 
Therefore, $\Delta Y_a\! = \! 0$. 
We also check that 
\begin{equation}
\label{ploopsy} \forall s\ino [0, \Delta \theta^\mathtt{b}_a ), \quad 
X_{s + \theta^{\mathtt{b}}_{a-} } \! = \! Y_a + A_a +  
X^\mathtt{r}_{s+  \subo^\mathtt{r} ((A_{a} )-) } \; .
\end{equation} 
Since $\theta^{\mathtt{b}}_{a}\! = \! a+ \subo^\mathtt{r} (A_{a})$, (\ref{ploopsy}) at $s\! = \! 0$ entails 
 $X( \theta^{\mathtt{b}}_{a-} )\! = \! Y_a$. Since $X(\theta^\mathtt{b}_{a-}-)\! = \! Y_{a-}$, we get 
$\Delta X \! (\theta^{\mathtt{b}}_{a-}) \! = \! \Delta Y_a \! = \! \Delta A_a \! = \! 0$. 
Next observe that for all $s\ino (0, \Delta \theta^\mathtt{b}_a )$, $A_a + X^\mathtt{r} ((s+  \subo^\mathtt{r} (A_{a}- ) )-)\!\!  >\! 0$; 
thus by (\ref{ploopsy}), for all $ t\ino \big( \theta^{\mathtt{b}}_{a-}, \theta^{\mathtt{b}}_{a} \big)$, we get $X_t  \! \geq \! X_{t-} \!  >  \! Y_a  $. 
Furthermore, if $a\! < \! T^*$, namely if $\theta^\mathtt{b}_a \! < \! \infty$, then 
$(i)$ entails that $X(\theta^\mathtt{b}_a)\! = \! Y_a$. This proves $(iii)$ when $\Delta A_a \! =\! 0$ and it completes the proof of the lemma. \cqfd 

\begin{lem}
\label{trajprop2}  Let $\alpha\ino \bbR$, $\beta \ino [0, \infty) $, $\kappa \ino (0, \infty) $ 
and $\mathbf{c}\! = \! (c_j)_{j\geq 1} \ino \elldo_3$ satisfy (\ref{varinfinie}). Namely: either $\beta \! >\! 0$ or $\sigma_2 (\mathbf{c})\! = \! \infty$. Then, the following holds true. 
\begin{compactenum}
\item[$(i)$] $\bP$-a.s.~if $(\Delta X^\mathtt{r} )(\Lambda^\mathtt{r}_t) \! >\! 0$, then there exists $a\ino [0, T^*]$ such that 
$\theta^\mathtt{b}_{a-} \! < \! t  \! < \! \theta^\mathtt{b}_{a} $. 

\smallskip

\item[$(ii)$] $\bP$-a.s.~for all $b\ino [0, \infty)$ such that $\Delta X^\mathtt{r}_b\! >\! 0$, there is a unique $t\ino [0, \infty)$ such that $\Lambda^\mathtt{r}_t \! = \! b$.

\smallskip

\item[$(iii)$] For all $t\ino [0, \infty)$, set $Q^\mathtt{b}_t \! = \! X^{\mathtt{b}}_{\Lambda^{\mathtt{b}}_t }$ and  $Q^\mathtt{r}_t \! = \! X^{\mathtt{r}}_{\Lambda^{\mathtt{r}}_t }$. Then, a.s.~for all $t\ino [0, \infty)$, 
$\Delta Q^\mathtt{b}_t \Delta Q^\mathtt{r}_t\! = \! 0$.   
\end{compactenum}
\end{lem}
\noi
\textbf{Proof.} Suppose that $(\Delta X^\mathtt{r} )(\Lambda^\mathtt{r}_t) \! >\! 0$. To simplify notation, we set $b\! = \! \Lambda^\mathtt{r}_t$ and $x\! = \! -\inf_{s\in [0, b]} X^\mathtt{r}_s$. 
Since $X^\mathtt{r}$ is a spectrally positive L\'evy process, $X^{\mathtt{r}}_b\! >\! -x$. Thus, $b \! <\!  \gamma^{\mathtt{r}}_x$; moreover, since no excursion above the infimum of the spectrally positive L\'evy process $X^\mathtt{r}$ starts with a jump we also get  
$\gamma^\mathtt{r}_{x-}\! <\!  b $. Thus, $\gamma^\mathtt{r}_{x-}\! <\!  b \! <\!  \gamma^{\mathtt{r}}_x$. 
We next set $a\! = \! \sup \{ s\ino [0, \infty) \! :\!  A_s \! <\!  x\}$. Then, $A_{a-} \! \leq  \! x \! \leq \!  A_a $ and we first prove the following:  
\begin{equation}
\label{hermanad}
 \theta^\mathtt{b}_{a-}\!\!  -a \! \leq \! \gamma^\mathtt{r}_{x-} \! < \!  b \! < \!   \gamma^\mathtt{r}_{x}\! \leq \!  \theta^\mathtt{b}_{a}\! - a \; .
\end{equation}
Let us first suppose that $\Delta A_a \! >\! 0$, then a.s~$\gamma^{\mathtt{r}} (A_{s}) \! \rightarrow \! \gamma^{\mathtt{r}} (A_{a-})$ as $s\! \rightarrow \! a-$, 
since a.s.~$ \{ y \ino [0, \infty), \Delta \subo^{\mathtt{r}}_y \! >\! 0 \} \cap \{ A_{z-}; z\ino [0, \infty) : \Delta A_{z} \! >\! 0 \} \! = \! \emptyset$ 
because $A$ and $\gamma^{\mathtt{r}}$ are independent. 
Since $A$ strictly increases, $\gamma^{\mathtt{r}} (A_{a-}) \! 
\leq \! \gamma^\mathtt{r} (x-)$. Similarly, a.s~$\gamma^{\mathtt{r}} (A_{s}) \! \rightarrow \! \gamma^{\mathtt{r}} (A_{a})$ as $s\! \rightarrow \! a+$, which implies that $\gamma^\mathtt{r} (x) \! \leq \! \gamma^{\mathtt{r}} (A_{a})$. Note that $\gamma^{\mathtt{r}} (A_{a-})\! = \! \theta^\mathtt{b}_{a-}\! -\! a$; and that  
$\gamma^{\mathtt{r}} (A_{a})\! = \! \theta^\mathtt{b}_{a}\! -\! a$, by definition. This implies (\ref{hermanad}). 
Now suppose that $\Delta A_a\! = \! 0$. Then, $A_{a-}\! = \! A_a\! = \! x$ and $\theta^\mathtt{b}_{a-}\!\!  - a \! = \! \gamma^\mathtt{r} (A_a -)\! = \! \gamma^{\mathtt{r}} (x-)$, which also implies (\ref{hermanad}). 


We next use (\ref{hermanad}) to prove $(i)$: first observe that it implies that 
$\theta^\mathtt{b}_{a-}\! < \! b+ a \! < \! \theta^\mathtt{b}_{a}$. But for all $s\ino (\theta^\mathtt{b}_{a-}, \theta^\mathtt{b}_{a})$, $\Lambda^\mathtt{b}_s\! = \! a$ and thus $\Lambda^\mathtt{r}_s\! = \! s-a$, which shows that on $(\theta^\mathtt{b}_{a-}, \theta^\mathtt{b}_{a})$, $\Lambda^\mathtt{r}$ is strictly increasing: since $b\! = \! \Lambda^\mathtt{r}_{a+b}\! = \! \Lambda^{\mathtt{r}}_t$, we get $t\! = \! a+b$ and finally $\theta^\mathtt{b}_{a-}\! < \! t \! < \! \theta^\mathtt{b}_{a}$, which completes the proof of $(i)$. 
 
 \smallskip
 
Next observe that $(ii)$ is a simple consequence of (\ref{bloodymary}) that shows that $\Lambda^\mathtt{r}$ is continuous and strictly increasing. 
%
 
 \smallskip

Let us prove $(iii)$: suppose that $\Delta Q^\mathtt{r}_t \! >\! 0$. Since $\Lambda^{\! \mathtt{r}}$ is continous, it implies that 
$ \Delta X^\mathtt{r} (\Lambda^\mathtt{r}_t) \! >\! 0$ and by $(i)$ there exists $a\ino (0, \infty)$, such that 
$\theta^\mathtt{b}_{a-} \! < \! t  \! < \! \theta^\mathtt{b}_{a} $; now observe that for all $s\ino (\theta^\mathtt{b}_{a-}, \theta^\mathtt{b}_{a})$, $\Lambda^\mathtt{b}_s\! = \! a$. Thus, $Q^\mathtt{b}$ is constant on this interval and it implies that 
$\Delta Q^\mathtt{b}_t \! =\!  0$. This proves $(iii)$. \cqfd 

\medskip

The following result is quite useful to deal with the supercritical cases. 
\begin{lem}
\label{trutruc} Let $\alpha\in \bbR$, $\beta \ino [0, \infty) $, $\kappa \ino (0, \infty) $ 
and $\mathbf{c}\! = \! (c_j)_{j\geq 1} \ino \elldo_3$ satisfy (\ref{varinfinie}). Namely: either $\beta \! >\! 0$ or $\sigma_2 (\mathbf{c})\! = \! \infty$. Assume that $\alpha \! < \! 0$ and recall that $\varrho$,  
the largest root of $\psi$, is positive.  
For all $a\ino [0, \infty)$ and for all bounded measurable functionals $F \! : \! \bD ([0, \infty), \bbR)^2\! \rightarrow \! \bbR$  
\begin{equation}
\label{retruc}
\bE \big[ F (Y_{\cdot \wedge a} , A_{\cdot \wedge a})\big]= \bE \big[ e^{\varrho A_a} F (Y_{\cdot \wedge a} , A_{\cdot \wedge a}) ; a<T^*  \big] \;. 
\end{equation}
Moreover, let $R \! \subset \! \Omega$ belong to the $\bP$-completion of the natural filtration generated by $(Y, A)$. If $\{ T^* \! >\! a \} \cap R$ is $\bP$-negligible then 
$R$ is $\bP$-negligible. 
%
%
\end{lem}
\noi
\textbf{Proof.} By (\ref{T*def}), a.s.~$\un_{\{ a < T^* \}}\! = \! \un_{\{ A_a < -I_\infty^{\mathtt{r}}\}}$ where 
$ -I_\infty^{\mathtt{r}} \! =\! -\inf_{[0, \infty)} X^{\mathtt{r}}$;  $-I_\infty^{\mathtt{r}}$ is an exponentially distributed r.v.~with parameter $\varrho$ that is independent from $(Y, A)$, which easily implies (\ref{retruc}). 
Next, denote by $\ccG$ the $\bP$-completion of the natural filtration generated by $(Y, A)$ and let $R\ino \ccG$ be such that $\{ T^*\!  >\!  a \} \cap R$ is $\bP$-negligible. Thus,
$\{ T^*\!  > \! a \} \cap R \ino \ccG$ and 
$$ \textrm{$\bP$-a.s.~} \quad  0 = \bE \big[ \un_{R}\un_{\{ T^* > a \}} | \ccG \big] = \bE \big[ \un_{R}\un_{\{ A_a < -\inf X^\mathtt{r} \}} | \ccG \big]= e^{-\varrho A_a} \un_{R} \; ,$$
which implies that $R$ is $\bP$-negligible. \cqfd

\medskip

We next consider the excursions of $X$ and $Y$ above their respective infimum. 
\begin{lem}
\label{XYexcu} Let $\alpha\in \bbR$, $\beta \ino [0, \infty) $, $\kappa \ino (0, \infty) $ 
and $\mathbf{c}\! = \! (c_j)_{j\geq 1} \ino \elldo_3$ satisfy (\ref{varinfinie}). Namely: either $\beta \! >\! 0$ or $\sigma_2 (\mathbf{c})\! = \! \infty$. For all $t\ino [0, \infty)$, we recall the following notation: $I_t \! = \! \inf_{s\in [0, t]} X_s$ and $J_t =  \inf_{s\in [0, t]} Y_s $. Then, the following holds true. 
\begin{compactenum}

\smallskip

\item[$(i)$] A.s.~for all $t\ino [0, \infty)$, $X_t \! \geq \! Y( \Lambda^\mathtt{b}_t)$. Then, a.s.~for all $t_1, t_2 \ino [0, \infty)$ such that $\Lambda^\mathtt{b}_{t_1} \! < \!  \Lambda^\mathtt{b}_{t_2}$, $\inf_{s\in [t_1, t_2]} X_s \! = \!  
\inf_{a\in [\Lambda^\mathtt{b} (t_1),  \Lambda^\mathtt{b} (t_2)]} Y_a$. It implies that a.s.~for all $t\ino [0, \infty)$, $I_t = J (\Lambda^{\! \mathtt{b}}_t)$.

\smallskip

\item[$(ii)$] A.s.~$\big\{ t \ino [0, \infty) : X_t \! >\! I_t \big\} =  \big\{ t \ino [0, \infty) : Y({\Lambda^{\! \mathtt{b}}_t}) \! >\! 
J({\Lambda^{\! \mathtt{b}}_t}) \big\} $.

\smallskip

\item[$(iii)$] A.s.~the set $\ccE\! = \! \big\{ a \ino [0, \infty) : Y_a \! >\! J_a \big\}$ is open. Moreover, 
if $(l,r)$ is a connected component of $\ccE$, then $Y_l\! = \! Y_r\! = \! J_l \! = \! J_r$ and 
for all $a\ino (l, r)$, we get $J_a\! = \! J_l$ and 
$Y_{a-}\! \wedge \! Y_a \! >\! J_l $.

\smallskip

\item[$(iv)$] Set $ \ccZ^\mathtt{b}\! = \! \{ a\ino [0, \infty) \! : \! Y_a\! = \! J_a \}$. Then, $\bP$-a.s.
\begin{equation}
\label{Jtpsloczer}
\textrm{$\forall a, z\ino [0, \infty)$ such that $a\! < \! z$,} \quad \Big( \ccZ^\mathtt{b}\!  \cap (a, z) \neq \emptyset  \Big) \Longleftrightarrow   \Big(  J_z \! < \!  J_a \Big).  
\end{equation} 

%
\end{compactenum}
\end{lem}
\noi
\textbf{Proof.} We first prove $(i)$. To that end, we fix $t\ino (0, \infty)$ and we set $a\! = \! \Lambda^{\! \mathtt{b}}_t$.  
Thus, $\theta^{\mathtt{b}}_{a-} \! \leq \! t \! \leq \! \theta^{\mathtt{b}}_{a}$. If $\Delta \theta^\mathtt{b}_a \! >\! 0$, then Lemma \ref{trajprop1} $(iii)$ implies that 
$X_s \! \geq \! Y_a=X (\theta^\mathtt{b}_{a-}\! -) $, for all $s\ino [  \theta^\mathtt{b}_{a-}, \theta^\mathtt{b}_{a} ]$. Thus, $X_t \! \geq \! Y_a$ (note that $a$ is possibly equal to $T^*$ and in this case, $\theta^{\mathtt{b}}_{a-} \! < \! \infty \! = \! \theta^{\mathtt{b}}_{a}$).  
If $\Delta \theta^\mathtt{b}_a \! =\! 0$, $t\! = \!\theta^\mathtt{b}_a $ and $X_t\! = \! Y_a$ by Lemma \ref{trajprop1} $(i)$. 
Thus, we have proved that 
a.s.~for all $t\ino [0, \infty)$, $X_t \! \geq \! Y (\Lambda^{\! \mathtt{b}}_t)$. 
 
We next fix $t_1, t_2 \ino [0, \infty)$ such that $a_1\! :=\! \Lambda^\mathtt{b}_{t_1} \! < \!  \Lambda^\mathtt{b}_{t_2}\! =  : a_2$, which implies that $a_1 \! < \! T^*\! $. Since $Y_a \! = \! X(\theta^{\mathtt{b}}_a)$ for all $a \ino [0, T^*)$ and since $Y_{T^*}\! = \! X (\theta^{\mathtt{b}}_{T^*\! -}\! -)$ (by Lemma \ref{trajprop1} $(iii)$ with $a\! = \! T^*$), we get 
$\inf_{[a_1, a_2]} Y\! \geq \! \inf_{[\theta^{\mathtt{b}} (a_1), \theta^{\mathtt{b}} (a_2-)]} X$. Since 
$t_1 \! \leq\!  \theta^\mathtt{b}_{a_1} \! < \!  \theta^\mathtt{b}_{a_2-} \! \leq \! t_2$, we get $\inf_{[a_1, a_2]} Y\! \geq \! \inf_{[t_1, t_2]} X$. 
But we have proved  also that 
$X_t \! \geq \! Y (\Lambda^{\! \mathtt{b}}_t)$ for all $t\ino [0, \infty)$: therefore, $ \inf_{[t_1, t_2]} X \! \geq \! \inf_{[a_1, a_2]} Y$. Thus, $ \inf_{[t_1, t_2]} X \!= \! \inf_{[a_1, a_2]} Y\!\! $, which completes the proof of $(i)$. 

%
%

Let us now prove $(ii)$. We first fix $t\ino (0, \infty)$ such that $X_t \! >\! I_t$ and 
we set 
$$g_t \! = \! \sup \big\{s \! <\!  t\! : \! X_s \! = \! I_s \big\} \quad \textrm{and} \quad d_t \! = \! \inf 
\big\{s \! >\!  t\! : \! X_s \! = \! I_s \big\}, $$
 with the convention that $\inf \emptyset\! = \! \infty$ (namely, in supercritical cases, $d_t\! = \! \infty$ iff $I_t \! =\! \inf_{s\in [0, \infty)} X_s$). 
Standard results on the excursions of spectrally positive processes above their infimum entails that 
$\Delta X (g_t)\! = \! 0$ and that $\Delta X (d_t)\! = \! 0$ if $d_t \! < \! \infty$. 
Consequently, for all $s\ino [g_t , d_t)$, $I_s\! = \! I_t\! = \! X(g_t)$ (and we also get $X(d_t)\! = \! X(g_t)$ if $d_t \! < \! \infty$).  

We still assume that $X_t \! >\! I_t$ and we next suppose that $Y(\Lambda^{\! \mathtt{b}}_t)  \! = \! J (\Lambda^{\! \mathtt{b}}_t) $ and 
to simplify, we set $a\! = \! \Lambda^\mathtt{b}_t$ so that $\theta^{\mathtt{b}}_{a-} \! \leq \! t \! \leq \! \theta^{\mathtt{b}}_{a}$. We first prove that $\Delta \theta^\mathtt{b}_a \! >\! 0$. Indeed suppose the contrary: namely, suppose that   
$\Delta \theta^\mathtt{b}_a \! =\! 0$;  then, Lemma \ref{trajprop1} $(i)$ asserts that 
$X_t\! = \! X( \theta^{\mathtt{b}}_{a})\! = \! Y_a= J_a$; however 
$J_a\! = \! I_t$ by $(i)$, which contradicts $X_t \! > \! I_t$. Thus, 
$\Delta \theta^\mathtt{b}_a \! >\! 0$ and Lemma \ref{trajprop1} $(iii)$ asserts that for all $s\ino (\theta^\mathtt{b}_{a-}, \theta^\mathtt{b}_a)$, 
$X_s \! >\! Y_a \! = \! X((\theta^\mathtt{b}_{a-})-)$.
Since $(i)$ entails $I_t \! = \! J_a$ and since we suppose $Y_a\! = \! J_a$, then 
for all $s\ino (\theta^\mathtt{b}_{a-}, \theta^\mathtt{b}_a)$, we get $X_s \! >I_t \! =\! X((\theta^\mathtt{b}_{a-})-) $. 
Thus, $g_t \! = \! \theta^\mathtt{b}_{a-}$ and 
$d_t \! = \! \theta^\mathtt{b}_{a}$ (which includes the possibility that $d_t \! = \! \infty$ in the supercritical cases). Since $\Delta X (g_t)\! = \! 0$, 
Lemma \ref{trajprop1} $(iii)$ entails that $\Delta A_a\! = \! \Delta X(\theta^\mathtt{b}_{a-})\! = \! 0$. 
Thus, we have proved that a.s.~for all $t\ino (0, \infty)$, if $X_t \! >\! I_t$ and if $Y(\Lambda^{\! \mathtt{b}}_t)  \! = \! J (\Lambda^{\! \mathtt{b}}_t)$ then 
$g_t \! = \! \theta^\mathtt{b}_{a-} \! <\! d_t \! = \! \theta^\mathtt{b}_{a}$ (where $a\! = \! \Lambda^\mathtt{b}_t$) and $\Delta A_a \! = \! 0$.  
We next use the following: for all $ \epp \ino (0, \infty)$, 
\begin{equation} 
\label{excpasrou}
\textrm{$\bP$-a.s.~} \quad 
\sum_{a\in [0, \infty)} \un_{\{ \Delta A_a = 0 \, ; \, \Delta \theta^\mathtt{b}_a >\epp \,  ; \, Y_a = J_a \}} \; = \; 0. 
\end{equation}
Before proving (\ref{excpasrou}), let us complete the proof of $(ii)$: first note that (\ref{excpasrou}) and the previous arguments entail 
that a.s.~for all $t\ino (0, \infty)$,  if $X_t \! >\! I_t$, then $Y(\Lambda^{\! \mathtt{b}}_t)  \! > \! J (\Lambda^{\! \mathtt{b}}_t)$. Then by $(i)$,  
if $X_t \! = \! I_t$, then $ J (\Lambda^\mathtt{b}_t) \! \l= \! Y (\Lambda^\mathtt{b}_t)$, which completes the proof of $(ii)$, provided that (\ref{excpasrou}) holds true. 

\smallskip

\noi
\textit{Proof of (\ref{excpasrou})}. Suppose that 
$\Delta \theta^\mathtt{b}_a \! >\! \epp$ (which does not exclude that $\Delta \theta^\mathtt{b}_a\! = \! \infty$ in the supercritical cases) 
and suppose that $\Delta A_a \! = \! 0$. Then $\theta^\mathtt{b}_{a-} \! = \! a +  \subo^\mathtt{r} (A_a-)$ and thus, 
$\Delta \theta^\mathtt{b}_a\! =\!  (\Delta \subo^\mathtt{r}) (A_a)$. 
Recall that for all $x\ino [0,\infty)$, we have set $A^{-1}_x \!= \! \inf \{ a\ino [0, \infty): A_a \! >\! x\}$. By Lemma \ref{AYnontriv}, a.s.~$A$ is strictly increasing and $A^{-1}$ is continuous and we get $A^{-1}(A_a)\! = \! a$; moreover, by definition we get $A(A^{-1}_x -) \! \leq \! x \! \leq \! A( A^{-1}_x)$ for all $x\ino (0, \infty)$. Then, (\ref{excpasrou}) is clearly a consequence of the following 
\begin{equation} 
\label{blurprou}
\textrm{$\bP$-a.s.~} \quad Q(\epp) =\sum_{x\in [0,\infty)} 
\un_{\{ \Delta A (A^{-1}_x) = 0 \, ; \, \Delta \subo^\mathtt{r} (x)>\epp \,  ; \, Y (A^{-1}_x)  = J (A^{-1}_x) \} } \; = \; 0. 
\end{equation}
Let us prove (\ref{blurprou}). Recall from (\ref{lifetimeexc}) that the measure $\bN (\zeta \ino dr)$ on $(0, \infty]$ 
is the L\'evy measure of the possibly defective subordinator $\subo^\mathtt{r}$ and recall that 
$\bN (\zeta\! = \! \infty)\! = \! \varrho$, that is the largest root of $\psi$. More precisely, set $\cM\! = \! \sum \delta_{(x, \Delta \gamma^\mathtt{r} (x))} $, where the sum is over the countable set of times $x$ where $\Delta \gamma^\mathtt{r} (x) \! >\! 0$; then, $\cM $ is distributed as $\un_{\{x \leq  \cE  \}} \cN(dxdr)$ where $\cN$ is a Poisson point process on $(0, \infty)\! \times \! (0, \infty] $ with intensity $dx \, \bN(\zeta \ino dr)$ and where 
$\cE\! =\!  \inf \{ x \ino [0, \infty) : \cN ([0, x] \! \times \! \{ \infty\} ) \! \neq \! 0 \}$ (therefore, $\cE$ 
is exponentially distributed with parameter $\bN (\zeta\! = \! \infty)\! = \! \varrho$). 
Since $(Y, A)$ and $X^\mathtt{r}$ are independent, we get the following. 
\begin{eqnarray*}
\bE \big[  Q(\epp) \big| (Y, A)\big]& = &  \bN( \zeta \ino (\epp, \infty]) 
 \int_0^{\infty} \!\! \!\!\! dx \, e^{-\varrho x}  \un_{\{  \Delta A (A^{-1}_x) = 0  \,  ; \, Y (A^{-1}_x)  = J (A^{-1}_x)  \} } \\
 & =& \bN( \zeta \ino (\epp, \infty])  \int_0^\infty \! \! \!\!\! dA_a e^{-\varrho A_a} \un_{\{ \Delta A_a = 0  \,  ; \, Y_a  = J_a \} }, 
 \end{eqnarray*}
the last equality being a consequence of an easy change of variable and of the fact that $A(A^{-1}_x) \! = \! x$ if $\Delta A (A^{-1}_x) \! = \! 0$.  
Set $c\! = \! \bN( \zeta \ino (\epp, \infty])$ and 
observe that $dA_a = \kappa \beta . a \, da + \sum_{a^\prime  \in [0, \infty)} \Delta A_{a^\prime} \delta_{a^\prime} (da) $. 
Thus, 
\begin{equation}
\label{furnouille}
 \bE \big[  Q(\epp) \big| (Y, A)\big] = c\kappa \beta \!  \int_0^\infty \!  \! \! \!\!\!   ae^{-\varrho A_a}  \un_{\{ Y_a  = J_a \} } da \;. 
\end{equation}
Then, we only need to prove that a.s.~$\mathtt{Leb} (\{ a\ino [0, \infty): Y_a \! = \! J_a \})\! = \! 0$. To that end, recall from Lemma  \ref{trajprop1} $(i)$ that for all $a\ino [0, T^*)$, $X(\theta^\mathtt{b}_a  ) \! =\!  Y_a$ and recall from $(i)$ above that 
$I( \theta^\mathtt{b}_a  ) \! =\!  J_a$. 
The change of variable $t\! = \! \theta^\mathtt{b}_a $, entails that 
$$ \int_0^{T^*}\!\!\!\!  \un_{\{ Y_a= J_a \}}\,  da =  \int_0^{T^*} \!\!\! 
 \un_{\{ X(\theta^\mathtt{b}_a  ) = I(\theta^\mathtt{b}_a  )  \} }\,  da \leq 
 \int_0^\infty  \!\!\!   \un_{\{ X_t= I_t  \} }\,  d\Lambda^\mathtt{b}_t  = 0\;,$$ 
 because $\Lambda^\mathtt{b}$ is $1$-Lipschitz and because a.s.~the Lebesgue measure of $\{ t\ino [0, \infty)\! : \! X_t \! = \! I_t \}$ is null. It then shows that a.s.~on $\{ a\! < \! T^*\}$, the Lebesgue measure of $\{ b\ino [0, a] \! : \! Y_b \! = \! J_b \}$ is null. By Lemma \ref{trutruc}, a.s.~$\mathtt{Leb} (\{ a\ino [0, \infty): Y_a \! = \! J_a \})\! = \! 0$ which implies  (\ref{blurprou}) (and thus (\ref{excpasrou})) by (\ref{furnouille}). It completes the proof of $(ii)$.

Let us prove $(iii)$. By standard results, 
$\ccE^\prime \! :=\! \{ t\ino [0, \infty)\! : \! X_t \! >\! I_t \}$ is open and if $(g,d)$ is a connected component of $\ccE^\prime$, then $X_g\!  = \! I_g\! = \! I_d$ (and $X_g\! = \! X_d$ if $d\! <\! \infty$) 
and for all $t\ino (g,d) $, $X_{t-} \! \wedge X_t \! >\! I_g$. Let $a\ino [0, \infty)$ be such that $ Y_a \! >\! J_a$; we set $l\! = \! \sup \{ b \ino [0, a] \! : \! Y_b\! = \! J_b\}$ and $r\! = \! \inf \{ b \ino (a, \infty) \! : \! Y_b\! = \! J_b\}$. We assume that $r\! < \! T^*$. 
Then, $X(\theta^{\mathtt{b}}_a)\! = \! Y_a>J_a \! = \! I(\theta^{\mathtt{b}}_a)$; let $(g, d)$ be the connected component of $\ccE^\prime$ containing $\theta^{\mathtt{b}}_a$. By $(i)$, $X_{d}\! = \! X_{g}\! = \! Y(\Lambda^{\! \mathtt{b}}_g)\! = \!  Y(\Lambda^{\! \mathtt{b}}_d)\! = \! 
J (\Lambda^{\! \mathtt{b}}_g) \! = \! J(\Lambda^{\! \mathtt{b}}_d)$ and by $(ii)$, $\Lambda^\mathtt{b}_g \! = \! l$ and $\Lambda^\mathtt{b}_d \! = \! r$ and $(l, r)$ is a connected component of  $\ccE\! = \! \big\{ b \ino [0, \infty)\! :\!  Y_b \! >\! J_b \big\}$. It easily shows that $[0, T^*)\cap \ccE$ is an open subset. 
Together with Lemma \ref{trutruc} this concludes the proof of $(iii)$. 
%

Let us prove $(iv)$. First recall from Section \ref{excuhaut} the notation: $\ccZ\! = \! \big\{ t \ino [0, \infty) \! :\!  X_t \! =\! I_t \big\}$ 
and recall that the continuous process $t\mapsto\!  -I_t$ is a local-time for $\ccZ$: 
in particular, recall from (\ref{zero}) that 
$\ccZ \! \cap \! (s,t)\! \neq \! \emptyset $ iff $I_t \! < \! I_s$. By $(ii)$, $\ccZ \! = \! \big\{ t \ino [0, \infty) : Y({\Lambda^{\! \mathtt{b}}_t}) \! =\! J({\Lambda^{\! \mathtt{b}}_t}) \big\} $; it easily implies the following: 
$\ccZ^\mathtt{b}\!  \cap (a, z) \neq \emptyset$ iff $\ccZ \cap (\theta^\mathtt{b}_a , \theta^\mathtt{b}_z) \! \neq \! \emptyset $ which is equivalent to $I (\theta^\mathtt{b}_z)\! = \! J_z \! < \! J_a \! = \! I(\theta^\mathtt{b}_a)$ (by $(i)$), which completes the proof of $(iv)$. \cqfd 

\medskip

We next recall the following result due to Aldous \& Limic \cite{AlLi98} (Proposition 14, p.~20).   
\begin{prop}[Proposition 14 \cite{AlLi98}]  
\label{AldLim1} Let $\alpha\in \bbR$, $\beta \ino [0, \infty) $, $\kappa \ino (0, \infty) $ 
and $\mathbf{c}\! = \! (c_j)_{j\geq 1} \ino \elldo_3$ satisfy (\ref{varinfinie}). Namely: either $\beta \! >\! 0$ or $\sigma_2 (\mathbf{c})\! = \! \infty$. Then, the following holds true. 
\begin{compactenum}

\smallskip

\item[$(i)$] For all $a\ino [0, \infty)$, $\bP (Y_a \! = \! J_a) \! = \! 0$. 

\smallskip

\item[$(ii)$] $\bP$-a.s.~the set $\{ a\ino [0, \infty)\! : \! Y_a\!  =\! J_a \}$ contains no isolated points. 

\smallskip

\item[$(iii)$] Set $M_a\! = \! \max \{ r\! -\! l\, ; \; r\! \geq \! l\! \geq \! a \! : \! \textrm{$(l,r)$ 
is an excursion interval of $Y\! -\! J$ above $0$}\}$. Then, $M_a \! \rightarrow \! 0$ in probability as $a\! \rightarrow \! \infty$.  
\end{compactenum}
\end{prop} 
\noi
\textbf{Proof.} The process $(Y_{s/\kappa})_{s\in [0, \infty)}$ is  the the process $W^{\kappa^\prime, -\tau, \mathbf{c}}$ 
in \cite{AlLi98}, where $\kappa^\prime \! = \! \beta / \kappa$ and $\tau\! =\! \alpha/ \kappa$ (note that the letter 
$\kappa$ plays another role in \cite{AlLi98}). 
Then $(i)$ (resp.~$(ii)$ and $(iii)$) is Proposition 14 \cite{AlLi98} $(b)$ (resp.~$(d)$ and $(c)$). \cqfd 

\bigskip

Thanks to Proposition \ref{AldLim1} $(iii)$, the excursion intervals of $Y\! -\! J$ above $0$ can be listed as follows
\begin{equation}
\label{excYre}
\{ a\ino [0, \infty)\! :  Y_a\! >  \! J_a \}= \bigcup_{k\geq 1} (l_k, r_k) \; .
\end{equation}
where $\zeta_k \! = \! r_k \! -\! l_k $, $k\! \geq \! 1$, is non-decreasing. 
Then, as a consequence of Theorem 2 in Aldous \& Limic \cite{AlLi98}, p.~4, we recall the following. 
\begin{prop}[Theorem 2 \cite{AlLi98}]  
\label{AldLim2} Let $\alpha\in \bbR$, $\beta \ino [0, \infty) $, $\kappa \ino (0, \infty) $ 
and $\mathbf{c}\! = \! (c_j)_{j\geq 1} \ino \elldo_3$ satisfy (\ref{varinfinie}). Namely: either $\beta \! >\! 0$ or $\sigma_2 (\mathbf{c})\! = \! \infty$. 
Then, $(\zeta_k)_{k\geq 1}$, that is the ordered sequence of lengths of the excursions of $Y\! -\! J$ above $0$, is distributed as the 
$(\beta / \kappa, \alpha/ \kappa, \mathbf{c})$-multiplicative coalescent (as defined in \cite{AlLi98}) taken at time $0$. In particular, we get a.s.~$\sum_{k\geq 1} \zeta_k^2 \! < \! \infty$. 
\end{prop}

\subsubsection{Proof of Theorem \ref{cHdefthm}}
\label{cHdefthmpf}
We first prove the following lemma. 
\begin{lem}
\label{rebout} Let $\alpha\ino \bbR$, $\beta \ino [0, \infty) $, $\kappa \ino (0, \infty) $ 
and $\mathbf{c}\! = \! (c_j)_{j\geq 1} \ino \elldo_3$ satisfy (\ref{recontH}). Then, the following holds true. 

\smallskip

\begin{compactenum}
\item[$(i)$] Almost surely, $a\ino [0, T^*) \! \mapsto \! H(\theta^\mathtt{b}_a)$ is continuous. 

\smallskip

\item[$(ii)$] For all $a\ino [0, \infty)
$, there exists a Borel measurable functional  $F_a \! : \! \bD ([0, \infty), \bbR)^2 \! \rightarrow \! 
 [0, \infty)$ such that a.s.~on the event $\{ a\! < \! T^*\}$, $ H(\theta^\mathtt{b}_a)\! = \! F_a (Y_{\! \cdot \, \wedge a} , A_{ \cdot \wedge a} )$. 
\end{compactenum}

\end{lem}
\noi
\textbf{Proof.} We first prove $(i)$. 
Since $H$ is continuous, $H\circ \theta^\mathtt{b}$ is c\`adl\`ag on $[0, T^*)$ and 
the left-limit at $a$ is $H (\theta^\mathtt{b}_{a-})$. If $\Delta  \theta^\mathtt{b}_{a}\! = \! 0$, then 
$H (\theta^\mathtt{b}_{a-}) \! = \! H (\theta^\mathtt{b}_{a})$.  
Suppose that $\Delta  \theta^\mathtt{b}_{a}\! > \! 0$: by Lemma  \ref{trajprop1} $(iii)$, for all 
$t\ino (\theta^\mathtt{b}_{a-}, \theta^\mathtt{b}_{a})$, we get  
$X_t \! >\! X ((\theta^\mathtt{b}_{a-})-) \! =\!  
X (\theta^\mathtt{b}_{a})$; we then apply Lemma \ref{cucueille} to $t_0\! = \! \theta^\mathtt{b}_{a-}$ and $t_1 \! = \! \theta^\mathtt{b}_{a}$: in particular we get $H_{t_0}\! =\! H_{t_1}$, namely: $H (\theta^\mathtt{b}_{a-}) \! = \! H (\theta^\mathtt{b}_{a})$. This proves that a.s.~$H\circ \theta^\mathtt{b}$ is continuous on $[0, T^*)$. 

We next prove $(ii)$ and we first suppose that $\beta \! >\! 0$. In that case recall from (\ref{brocueille}) that $H_t$ is the $2/\beta$ times the Lebesgue measure of the set $\{ \inf_{r\in [s, t]} X_r ; s\ino [0, t] \}$. 
We fix $a\! < \! T^*$ and $s \ino [0,\theta^\mathtt{b}_{a}]$. To simplify, we also set $b=\Lambda^\mathtt{b}_s$. If $s\ino [\theta^\mathtt{b}_{a-},  \theta^\mathtt{b}_{a}]$, then $b\! = \! a$ and 
Lemma \ref{trajprop1} $(iii)$ entails that $\inf_{r\in [s, \theta^\mathtt{b}_a]} X_r = X( \theta^\mathtt{b}_a)\! = \! Y_a$; if $s\ino [0,   \theta^\mathtt{b}_{a-})$, then $b<a$ and Lemma \ref{XYexcu} $(i)$ entails that $\inf_{z\in [b, a]} Y_z\! = \! \inf_{r\in [s, \theta^\mathtt{b}_a]} X_r $. This easily implies that 
$$ \big\{  \!\!\!\!\!  \inf_{\quad r\in [s,  \theta^\mathtt{b}_{a}]} \!\!\!\!\! X_r \, ; \;  s\in [0,  \theta^\mathtt{b}_{a}]  \big\}= 
\big\{ \!\! \inf_{\; \; z\in [b, a]} \!\!\!  Y_z\,  ;\;  b \ino [0, a]\big\} \; .$$
Consequently, a.s.~on $\{ a\! < \! T^*\}$, $H (\theta^\mathtt{b}_{a}) \! = \! 2\beta^{-1} \mathtt{Leb}  (\{  
\inf_{z\in [b, a]}   Y_z\,  ;\;  b \ino [0, a] \})$, which implies $(ii)$ when $\beta\! >\! 0$. 

We next suppose that $\beta \! = \! 0$. Recall that (\ref{recontH}) implies (\ref{varinfinie}), so we get $\sigma_2 (\mathbf{c})\! = \! \infty$. Then, recall from (\ref{infideff}) the following notation: 
$ \ccH^\varepsilon(t) \! = \! \big\{s\ino (0, t]: X_{s-} + \varepsilon \! <\!  \inf_{r\in [s, t]} X_r   \big\}$, for all $\epp, t\ino (0, \infty)$. 
We next fix $a\! < \!  T^*$ such that $\Delta A_a\! >\! 0$. 
Note that we necessarily have $\Delta \theta^\mathtt{b}_a \ino (0, \infty)$. Let $s_1, s_2 \ino  \ccH^\varepsilon 
(\theta^\mathtt{b}_{a-})$ be such that $s_1 \! < \! s_2 \! < \! \theta^\mathtt{b}_{a-}$. For all $i\ino \{ 1, 2 \}$, we set $b_i \! = \! \Lambda^\mathtt{b}_{s_i}$; by definition $ \theta^\mathtt{b}_{b_i-} \! \leq \! s_i \! \leq \!  \theta^\mathtt{b}_{b_i}$. 
If $s_i \! >\!  \theta^\mathtt{b}_{b_i-}$, then Lemma \ref{trajprop1} $(iii)$ implies that $X_{s_i-}\! \geq \! 
X ( \theta^\mathtt{b}_{b_i})$, which contradicts $s_i  \ino  \ccH^\epp (\theta^\mathtt{b}_{a-})$. 
Therefore, $s_i \! =\!  \theta^\mathtt{b}_{b_i-}$, which implies that $b_1 \! < \! b_2 <\! a$. On the other hand,  Lemma \ref{trajprop1} $(i)$ and Lemma \ref{XYexcu} $(i)$ tell us that $Y_{b_{i}-}=X(\theta^{\mathtt{b}}_{b_{i}})$ and 
$\inf_{z\in [b_i, a]} Y_z\! = \! \inf_{r\in [s_i , \theta^{\mathtt{b}}_{a-}] } X_r$. 
Therefore, $b_1, b_{2}$ belong to the set $\ccY^\epp_a\! := \! \{ b \ino [0, a ] \! : \! Y_{b-}+ \epp < \inf_{z\in [b, a]} Y_z \}$. 
Next, suppose that $b_{1}, b_{2}\in \ccY^\epp_a$. Since a.s.~$\Delta Y_{a}=0$ (because $Y$ and $A$ have distinct jump times), we get $b_{i}<a$, $i=1,2$. Then a similar argument based on Lemma \ref{trajprop1} $(i)$ and Lemma \ref{XYexcu} $(i)$ show that $s_{i}= \theta^{\mathtt{b}}_{b_{i}-}<\theta^{\mathtt{b}}_{a-}$ and that \mm{$b_{i}$} is an element of $ \ccH^\varepsilon(\theta^{\mathtt{b}}_{a-})$. 
We have proved that 
$\Lambda^\mathtt{b}$ is one-to-one from $ \ccH^\epp (\theta^\mathtt{b}_{a-})\! \setminus \! \{\theta^{\mathtt{b}}_{a-}\}$ 
onto $\ccY^\epp_a$, which then implies  
the following: for all $\epp \ino (0, \infty)$, 
$$ \textrm{a.s.~for all $a\ino [0, T^*)$ such that $\Delta A_a\! >\! 0$,} \quad \# \ccY^\epp_a\; \le \# \ccH_{\theta^\mathtt{b}_{a-}}^\varepsilon \le  \# \ccY^\epp_a+1. $$

Then, (\ref{cueille}) easily implies that there is $(\epp_k)_{k\in \bbN}$ decreasing to $0$ such that a.s.~for all $t\ino [0, \infty)\cap \bbQ$ and for all $s\ino [0, t]$ such that $X_{s-} \! \leq \! \inf_{r\in [s, t]} X_r$, 
 $H_s\!  =\!  \lim_{k \rightarrow \infty}
(\# \ccH^{\varepsilon_k}_s)/ q(\varepsilon_k)$. Then observe that for all $a\ino [0, T^*)$ such that $\Delta A_a\! >\! 0$, there is $t\ino \bbQ \cap (\theta^\mathtt{b}_{a-}, \theta^\mathtt{b}_{a})$ and by Lemma \ref{trajprop1} $(iii)$, we get $X(\theta^\mathtt{b}_{a-}\! -)\! \leq \! \inf_{r\in [\theta^\mathtt{b}_{a-}, t] } X_r$. Thus, a.s.~for all $a\ino [0, T^*)$ such that $\Delta A_a\! >\! 0$, we get 
\begin{equation} 
\label{glabouboul}
H(\theta^\mathtt{b}_{a})\! = \! H(\theta^\mathtt{b}_{a-})\! = \! \lim_{k\rightarrow \infty} \frac{_1}{^{q(\epp_k)}} \# \ccH_{\theta^\mathtt{b}_{a-}}^{\epp_k}\! = \!  \lim_{k\rightarrow \infty} \frac{_1}{^{q(\epp_k)}}  
\# \ccY^{\epp_k}_a\; .
\end{equation}
Since $\sigma_2 (\mathbf{c}) \! = \! \infty$, the jump-times of $A$ form a dense subset of $[0, \infty)$. Thus it makes sense to set the following for all $a\ino (0, \infty)$: 
$$F_a (Y_{\! \cdot \, \wedge a} , A_{ \cdot \wedge a} )\! = \!\!\!\!\! 
  \limsup_{\quad z\rightarrow a-,\,  \Delta A_{z} >0} \!\! 
 \; \,  \limsup_{k\rightarrow \infty}  \; \frac{_1}{^{q(\epp_k)}}  
\# \ccY^{\epp_k}_{z} \; \textrm{if this quantity is finite and $0$ otherwise.}$$ 
Since $H \circ  \theta^\mathtt{b}$ is continuous on $[0, T^*)$, (\ref{glabouboul}) entails $(ii)$. \cqfd

\medskip

Note that in (sub)critical cases, the previous lemma proves Theorem \ref{cHdefthm}. We next assume that $\alpha \! < \! 0$: then $\varrho \! >\! 0$ and 
$T^* \! < \! \infty$ a.s. For all $a\ino [0, \infty)$, we set 
$$ \cH_{a} = \limsup_{q\rightarrow a-, \, q\in \bbQ} F_q (Y_{  \cdot  \wedge q}, 
A_{ \cdot  \wedge q} ) \; \textrm{if this quantity is finite and $0$ otherwise.}$$
For all $a\ino (0, \infty)$, we set $R_a \! = \! \{ \omega\ino \Omega  \! :  z\ino [0, a]\! \mapsto \! \cH_{z} (\omega) \; \textrm{not continuous} \}$. By Lemma \ref{rebout} $(ii)$, a.s.~on $\{ T^*\!>\! a\}$, for all $q\ino [0, a] \cap \bbQ$, $F_q (Y_{  \cdot  \wedge q}, 
A_{ \cdot  \wedge q} )\! = \! H(\theta^\mathtt{b}_q)$; consequently, by Lemma \ref{rebout} $(i)$, 
a.s.~on $\{ T^*\! >\! a\}$, for all $z \ino [0, a]$, $\cH_{z} \! = \! H( \theta^\mathtt{b}_{z})$ and $R_a \cap \{ T^*\! >\! a\}$ is a $\bP$-negligible set. Then, Lemma \ref{trutruc} entails that $R_a $ is $\bP$-negligible, which implies that 
$\bP$-a.s.~$\cH$ is continuous since $a$ can be chosen arbitrarily large.
%
%
%
Moreover, the previous arguments imply that a.s.~$\cH \! = \! H\circ \theta^\mathtt{b}$ on $[0, T^*)$. 
Therefore, for all $a\ino [0, \infty)$ a.s.~on $\{ T^* \! >\! a\}$, for all $z\ino [0, a]$, $\cH_{z} \! = \! F_{z} (Y_{\cdot \wedge z}, A_{\cdot \wedge z})$. By Lemma \ref{trutruc}, it implies that for all $a\ino [0, \infty)$, a.s.~$\cH_a \! = \! F_a (Y_{\cdot \wedge a}, A_{\cdot \wedge a})$, which completes the proof of Theorem \ref{cHdefthm}. \cqfd 

\medskip

We shall need the following lemma that concerns the excursions of $\cH$ above $0$.  
\begin{lem}
\label{AHeuer} Let $\alpha\in \bbR$, $\beta \ino [0, \infty) $, $\kappa \ino (0, \infty) $ 
and $\mathbf{c}\! = \! (c_j)_{j\geq 1} \ino \elldo_3$ satisfy (\ref{recontH}). Then, the following holds true. 

\smallskip

\begin{compactenum}
\item[$(i)$]  Almost surely for all $t\ino [0, \infty)$, $H_t \! \geq \! \cH( \Lambda^\mathtt{b}_t)$ and a.s.~for all $t_1, t_2 \ino [0, \infty)$ such that $\Lambda^\mathtt{b}_{t_1} \! < \!  \Lambda^\mathtt{b}_{t_2}$, 
$\inf_{s\in [t_1, t_2]} H_s \! = \!  
\inf_{a\in [\Lambda^\mathtt{b} (t_1),  \Lambda^\mathtt{b} (t_2)]} \cH_a$.

\smallskip

\item[$(ii)$] Almost surely 
$\big\{ a \ino [0, \infty) : Y_a \! >\! J_a \big\} \! =\!  \big\{ a \ino [0, \infty) : \cH_a \! >\! 0 \big\}$. 
\end{compactenum}

\end{lem}
\noi
\textbf{Proof.} Let $t\ino [0, \infty)$ and set $a\! = \!  \Lambda^{\! \mathtt{b}}_t$. Then, by definition 
$\theta^{\mathtt{b}}_{a-} \! \leq \! t \! \leq \!  \theta^{\mathtt{b}}_{a}$. If $\Delta \theta^{\mathtt{b}}_a\! = \! 0$, then 
$\theta^{\mathtt{b}}(\Lambda^{\mathtt{b}}_t)\! = \! t$ and 
$H_t \!= \! \cH (\Lambda^{\mathtt{b}}_t)$. Suppose next that $\Delta \theta^{\mathtt{b}}_a\! > \! 0$; by Lemma \ref{trajprop1} $(iii)$, 
$X_s\! >\! X ((\theta^{\mathtt{b}}_{a-})-)\! = \! Y_a\! = \! X (\theta^{\mathtt{b}}_a)$, for all $s\ino  (\theta^{\mathtt{b}}_{a-}, \theta^{\mathtt{b}}_{a})$; thus, we can apply Lemma \ref{cucueille} 
with $t_0\! = \!  \theta^{\mathtt{b}}_{a-}$ and $t_1 \! = \!  \theta^{\mathtt{b}}_{a}$, to get 
$H_t \! \geq \! H( \theta^{\mathtt{b}}_{a-}) \! = H (\theta^{\mathtt{b}}_{a})\! = \! \cH_a \! = \! \cH (\Lambda^{\! \mathtt{b}}_{t})$. To complete the proof of $(i)$ we argue exactly as in Lemma \ref{XYexcu} $(i)$. 
 
%

Let us prove $(ii)$. 
Recall from (\ref{excusame}) and from Lemma \ref{XYexcu} $(ii)$ 
that $\big\{ t \ino [0, \infty) : X_t \! >\! I_t \big\}\! = \! \big\{ t \ino [0, \infty) : H_t \! >\! 0 \big\} \! = \! \big\{ t \ino [0, \infty) : Y({\Lambda^{\! \mathtt{b}}_t}) \! >\! J({\Lambda^{\! \mathtt{b}}_t}) \big\} $.  
Then, observe that on $\{ T^* \! >\! a\}$, $Y_a \! >\! J_a$ iff $X(\theta^\mathtt{b}_a)\! =\!  Y_a \! >\! J_a\! = \!  I(\theta^\mathtt{b}_a)$, which is also equivalent to 
$\cH_a \! = \! H (\theta^\mathtt{b}_a) \! >\! 0$. This proves that for all $a\ino (0, \infty)$ 
a.s.~on $\{ T^* \! >\!  a\}$, $\{z\ino (0, a)\! : \! Y_{z} \! >\! J_{z}\}\! = \! \{z\ino (0, a)\! : \! \cH_{z} \! >\! 0\}$, which proves $(ii)$ in (sub)critical cases; in supercritical cases, Lemma~\ref{trutruc} applies. \cqfd

\subsubsection{Embedding into a L\'evy tree. Proof of Proposition \ref{fractCIRG}} 
\label{embedGinT}
We now explain how continuous multiplicative graphs are embedded in L\'evy trees. 
We first fix $a\ino (0, \infty)$ and we argue on the event $\{ T^*\! >\! a \}$. 
Let $(l,r)$ be an excursion interval of $\cH$ above $0$ such that $r\! < \! a$. 
By Lemma \ref{XYexcu} $(ii)$, there exists an excursion interval $(\mathbf{l}, \mathbf{r})$ of $X$ above its infimum process $I$, or equivalently (see (\ref{excusame})) an excursion interval of $H$ above $0$, such that 
$\Lambda^\mathtt{b}_{\mathbf{l}} \! = \! l$ and $\Lambda^\mathtt{b}_{\mathbf{r}} \! = \! r$, which also entails $\theta^\mathtt{b}_{l}\! = \! \mathbf{l}$ and $\theta^\mathtt{b}_{r}\! = \! \mathbf{r}$. To simplify notation, for all $s\ino [0, \infty)$,we set : 
\begin{equation*}
\label{}
 \mathtt{H}(s)\! = \! \cH_{(l+s)\wedge r}, \quad 
 \mathbf{H}_s\! = \! H_{(\mathbf{l} +s)\wedge \mathbf{r} } ,  \quad \btheta_{\! s} = \theta^\mathtt{b}_{(l+s )\wedge r} - \mathbf{l} , \quad  \zeta\! = \! r\! -\! l \quad \textrm{and} \quad \bzeta \! = \! \mathbf{r}\! -\! \mathbf{l}
\end{equation*}
Recall that 
$H (\theta^{\mathtt{b}}_a)\! = \! \cH_a$ on $(l, r)$, thus 
$\mathbf{H} (\btheta_{\! a} ) \! = \! \mathtt{H} (a)$, $a\ino [0, \zeta]$.  
Recall from (\ref{pseudometric}) the definition of the pseudometric $d_h$ coded by a function $h$. The previous arguments, combined with Lemma \ref{AHeuer} $(i)$ imply the following. 
\begin{equation}
\label{tatapioca}
 \forall a, b\ino [0, \zeta], \quad d_{\mathtt{H}} (a,b )= 
 d_{\mathbf{H}} \big( \btheta_a, \btheta_{b} \big)\; .
 \end{equation}
Recall from (\ref{codef}) in Section \ref{codtreesec} that $(T_h, d_h, \rho_h, m_h)$ stands for the rooted compact measured real tree 
coded by $h$ and recall that $p_h \! : \! [0, \zeta_h) \! \rightarrow \! T_h$ is the canonical projection. 
To simplify notation, we set 
$$ \big(\bT_{\!} , \delta,  \rho, \m^* \big):= \big(T_{\mathbf{H}} , d_{\mathbf{H}}, 
\rho_{\mathbf{H}},  
\bm_{\mathbf{H}} \big). $$
 Then (\ref{tatapioca}) implies the following: set $\cT_{\! } \! = \! p_{\mathbf{H}} (\btheta  ([0, \zeta]))$, then 
\begin{equation}
\label{ttppcc}
 (\cT_{\! } , \delta \,_{\! | \cT_{\! } \times \cT_{\! } }, \, \rho, \, \m^* (\cdot \cap \cT_{\! }) \big) \; \textrm{is isometric to} \; \big(T_{\mathtt{H}} , d_{\mathtt{H}}, \rho_{\mathtt{H}}, m_{\mathtt{H}} \big) .
\end{equation}
Namely, we view the tree coded by $\mathtt{H}$ as a compact subtree (namely, a compact connected subset) of the $\psi$-L\'evy tree coded by the excursion $\mathbf{H}$. 

Next, recall from (\ref{Poisurcon}) and (\ref{pinchset}) the definition of the set of pinching times 
$\Ptt$ on $(0, \infty)^2$; for all $k\! \geq \! 1$, recall from  (\ref{defPik}) the definition of the set of pinching times $\Ptt_k$ that fall in the $k$-th longest excursion of $Y$ above its infimum process $J$; note that there exists a $k$ such that $l_k\! = \! l$ and $r_k \! = \! r$ and then set $\Pi\! = \! \Ptt_k$ that is therefore the set of pinching times of $\Ptt$ that fall in the excursion interval $(l,r)$.   
Denote by $(G,d,\rho, m)$ the compact metric space coded by $\mathtt{H}$ and the pinching setup $(\Pi, 0)$ as defined in (\ref{defgrgraf}) (then, $(G,d,\rho, m)\! =\! (\bG_k, d_k, \rho_k, \bm_k)$, the 
$k$-th largest connected component of the multiplicative graph). 
We then set $\Pi^*\!  =\!  \big\{  \big(p_{\mathbf{H}} (\btheta_{\! s} ) , p_{\mathbf{H}} (\btheta_{t} )\big); (s,t) \ino \Pi \big\}$. Then, thanks to (\ref{tatapioca}), we see that: 
\begin{multline}
\label{tapiopioca}
(G, d, \rho, m)  \; \textrm{is isometric to the $(\Pi^*, 0)$-metric space} \\ 
\textrm{associated to} \;  (\cT_{\!} , \delta_{} \,_{| \cT_{\!} \times \cT_{\! } }, \rho, \m^* (\cdot \cap \cT) \big) \; .
\end{multline} 
%
To summarise, up to the identifications given by (\ref{ttppcc}) and (\ref{tapiopioca}), \textit{the connected component $G$ of the multiplicative continuous random graph corresponding to excursion interval $(l,r)$ is obtained as a finitely pinched metric space associated with the real tree $\cT_{\! }$ coded by $\mathtt{H}$ that is a subtree of the L\'evy tree $\bT_{\! }$ coded by $\mathbf{H}$}. This allows to prove Proposition \ref{fractCIRG} as follows.

\medskip

\noi
\textbf{Proof of Proposition \ref{fractCIRG}.} 
 We introduce the following exponents: 
 $$\gamma \! = \! \sup \{ r\ino [0, \infty): \lim_{\lambda \rightarrow \infty}\psi (\lambda) \lambda^{-r}\! = \! \infty\}\quad \textrm{and} \quad \eta \! = \!  \inf 
\{ r\ino [0, \infty): \lim_{\lambda \rightarrow \infty}\psi (\lambda) \lambda^{-r}\! = \! 0 \}\; .$$ Recall from (\ref{PoisdecH}) that $\bN$ stands for the excursion measure of the $\psi$-height process $H$ above $0$ and denote by $(T_H, d_H, \rho_H, m_H)$ the generic rooted compact measured real tree coded by $H$ under $\bN(dH)$.  
Theorem 5.5 in Le Gall \& D.~\cite{DuLG05}, p.~590, asserts that if $\gamma \! >\! 1$, then 
$\bN(dH)$-a.e.~$\mathrm{dim}_{\mathrm{H}} (T_H)= \eta/ (\eta \! -\! 1)$ and $\mathrm{dim}_{\mathrm{p}} (T_H)= \gamma/ (\gamma \! -\! 1)$ (this statement is a specific case of Theorem 5.5 
in \cite{DuLG05} where $E\! =\! [0, \infty)$). 
Moreover, in the proof of \ts{the} Theorem 5.5 \cite{DuLG05}, two estimates for the local upper- and lower-densities of the mass measure $m_H$ are given at (45) and (46) in \cite{DuLG05}, p.~593: 
for all 
$u \ino (0,  \frac{_\eta}{^{\eta  - 1}})$ and $v \ino (0,  \frac{_\gamma}{^{\gamma  - 1}})$, 
$\bN(dH)$-a.e.~for $m_H$-almost all $\sigma\ino T_H$, $\limsup_{r \rightarrow 0}  r^{-u} m_H (B(\sigma, r)) \! < \! \infty $ and $\liminf_{r \rightarrow 0}  r^{-v} m_H (B(\sigma, r)) \! < \! \infty $ (actually, within the notations of 
\cite{DuLG05}, if $E\! = \! [0, \infty)$, then $d(E)\! = \! 1$ and $\kappa (d\sigma)\! = \! m_H (d\sigma)$). Since $\bT$ is the tree coded by $\mathbf{H}$ that is an excursion of $H$ above $0$, 
the previous estimates and Theorem 5.5 in \cite{DuLG05} show that for all 
$u \ino (0,  \frac{_\eta}{^{\eta  - 1}})$ and for all $v \ino (0,  \frac{_\gamma}{^{\gamma  - 1}})$, we have the following:
\begin{multline}
\label{dindonne}
\textrm{$\bP$-a.s.~on $\{ T^*>a\}$,} \quad  \mathrm{dim}_{\mathrm{H}} (\bT)\! = \! \frac{\eta}{\eta\! -\! 1} \quad \textrm{and} \quad 
 \mathrm{dim}_{\mathrm{p}} (\bT)\! = \! \frac{\gamma}{\gamma \! -\! 1} \; \, \\
\textrm{and for $\m^*$-almost all $\sigma\ino \bT$} \quad  \limsup_{r \rightarrow 0}   
\frac{\m^* (B(\sigma, r))}{r^u} \! < \! \infty \quad  \textrm{and} \quad  \liminf_{r \rightarrow 0}  \frac{\m^* (B(\sigma, r))}{r^v}  \! < \! \infty . 
\end{multline} 
We now apply Lemma \ref{dimpinch} in Appendix Section \ref{Pinfrac} to $E_0\! = \! \bT$, $E\! = \! \cT$ and $E^\prime\! = \! G$ (Lemma \ref{dimpinch} applies since $\m^* (\cT)\! >\! 0$). 
Informally speaking, Lemma \ref{dimpinch} asserts the following: since $G$ is obtained from $\cT$ by identifying only a finite number of points, it does not modify Hausdorff and packing dimensions that are obtained here as local exponents of $\m^*$; moreover since $\m^* (\cT)\! >\! 0$ and since the local exponents are constant on L\'evy trees, the local exponents of $\m^*$ are the same on $\cT$ (and thus on $G$) 
and on $\bT$. Thus, Lemma \ref{dimpinch} and (\ref{dindonne}) entail that  $\bP$-a.s.~on $\{ T^*>a\}$, $\mathrm{dim}_{\mathrm{H}} (G)\! = \! \eta/ (\eta \! -\! 1)$ and $\mathrm{dim}_{\mathrm{p}} (G)\! = \!  \gamma / (\gamma \! -\! 1)$. Lemma \ref{trutruc} easily entails that it holds true for all excursions of $\cH$ above $0$ (not just for those ending before $T^*$). 
Namely, $\bP$-a.s.~for all $k\! \geq \! 1$, $\mathrm{dim}_{\mathrm{H}} (\bG_{\! k})\! = \! \eta/ (\eta \! -\! 1)$ and $\mathrm{dim}_{\mathrm{p}} (\bG_{\! k})\! = \!  \gamma / (\gamma \! -\! 1)$.

Observe that if $\beta>0$, then $\gamma=\eta=2$. Thus, to complete the proof of Proposition \ref{fractCIRG}, it remains to prove that the exponents $\gamma$ and $\eta$ are given by (\ref{expone}) when $\beta \! = \! 0$: set $\pi (dr)\! = \! 
\sum_{j\geq 1} \kappa c_j \delta_{c_j}$, the L\'evy measure of $X$. 
By an immediate calculation, we get  $\psi^\prime(\lambda)-\alpha \! = \! \int_{(0, \infty)} (1-e^{-\lambda r}) r\pi (dr) $. We next introduce 
$$J (x)\! =\!  x^{-1}\!\!\! \int_0^x \!\!\! du \! \int_{(u, \infty)} \!\!\!\!\!\!  \!\!\! \! r\pi (dr) =\int_{(0, \infty)} \!\!\! \! \!\!\!\!\!\!
r(1 \wedge (r/x)) \pi(dr)\! = \! \sum_{j\geq 1 }\! \kappa c_j^2 (1\wedge (c_j/x)) $$
as in Proposition \ref{fractCIRG} $(ii)$. 
As explained in Bertoin's book \cite{Be96} Chapter III, 
general arguments on the Laplace exponents of 
subordinators entail that there exist two universal constants $k_1, k_2 \ino (0, \infty)$ such that 
$$ k_1 J (1/\lambda) \leq \psi^\prime(\lambda)-\alpha\leq  k_2 J (1/\lambda) \; .$$
Since $\psi (\lambda) \! \leq \!\lambda  \psi^\prime(\lambda) \! \leq \! 4\psi (\lambda)$, by convexity, the previous inequality entails: $\gamma \! = \! 1+ \sup \{ r \ino (0, \infty): \lim_{x\to 0+} x^{r}J(x)\! = \! \infty \}$ and $\eta \! = \! 1+ \inf \{ r \ino (0, \infty): \lim_{x\to 0+} x^{r}J(x)\! = \! 0 \}$, which completes the proof of Proposition \ref{fractCIRG}. \cqfd

\appendix
\section{Pinched metric spaces and their fractal dimensions.}
\label{Pinfrac}
Let $(E, d)$ be a metric space. We briefly recall from Section~\ref{pinmetpar} the definition of pinched metrics: for all $i\ino \{ 1, \ldots , p \}$, let $(x_i, y_i)\ino E^2$; set $\Ptt\! = \! ((x_i, y_i))_{1\leq i\leq p}$; 
let $\epp \ino [0, \infty)$. Set  $A_E \! =\! \{ (x,y) ; x,y\ino E\}$ and for all $e\! = \! (x,y)\ino A_E$, we set $\underline{e}\! = \! x$ and $\overline{e}\! = \! y$. A path $\gamma$ joining $x$ to $y$ is a sequence of $e_1, \ldots , e_q \ino A_E$ such that $\underline{e}_1\! = \! x$, $\overline{e}_q\! = \! y$ and $\overline{e}_i\! = \! \underline{e}_{i+1}$, for all $i\ino \{ 1, \ldots , q-1 \}$. Next, we set $A_{\Ptt}\! = \! \{ (x_i, y_i), (y_i, x_i) ; 1\! \leq \! i\! \leq \! p \}$ and we define the length $l_e$ 
of an edge $e$ by setting $l_e\! = \! \epp \! \wedge \! d(\underline{e}, \overline{e})$ if $e \ino A_{\Ptt}$ and $l_e\! = \! d(\underline{e}, \overline{e})$ otherwise.  The length of a path $\gamma\! = \! (e_1, \ldots , e_q)$ is given by 
$l(\gamma)\! = \! \sum_{1\leq i\leq q} l_{e_i}$.  Then, recall from (\ref{pinchdist}), that the $(\Ptt, \epp)$-pinched pseudo-distance between $x$ and $y$ in $E$ is given by $d_{\Ptt , \epp} (x,y)\! = \! \inf \{ l(\gamma); \textrm{$\gamma$ is a path joining $x$ to $y$}\}$. 
We easily check that 
\begin{multline}
\label{vnaglurns}
d_{\Ptt , \epp} (x,y)\! = \! d(x,y) \wedge \min \big\{ \, l (\gamma)\, ; \;  \gamma \!  = \! (e_0, e^\prime_0, \ldots ,e_{r-1}, e^\prime_{r-1}, e_r), \\ \textrm{a path joining $x$ to $y$ such that}\;  e_0^\prime, \ldots e^\prime_{r-1}\ino A_{\Ptt} \; \textrm{and} \; r\! \leq \!  p  
\big\} .  
\end{multline}
Clearly, $d_{\Ptt, \epp}$ is a pseudo-metric and we denote the equivalence relation $d_{\Ptt, \epp}  (x,y)\! = \! 0$
by $x \equiv_{\Ptt, \epp} y$; the quotient space $E/\!\! \equiv_{\Ptt, \epp}$ equipped with $d_{\Ptt, \epp}$ is the 
$(\Ptt, \epp)$-pinched metric space associated with $(E, d)$. Recall that $\varpi_{\Ptt, \epp}\! : \! E \! \rightarrow \!  E/\!\! \equiv_{\Ptt, \epp}$ stands for the canonical projection that is $1$-Lipschitz. Note of course that if 
$\epp \! >\! 0$, then $d_{\Ptt , \epp}$ is a true metric on $E$, which is obviously identified with 
$E/\!\! \equiv_{\Ptt, \epp}$ and $\varpi_{\Ptt, \epp}$ is the identity map on $E$. 

Next, set $S\! = \! \{ x_i , y_i; 1\! \leq \! i \! \leq \! p \}$ and for all $x\ino E$, set 
$d(x, S)\! = \! \min_{y\in S} d(x,y)$. Then, (\ref{vnaglurns}) immediately entails that 
\begin{equation}
\label{populonia}
\forall x, y \ino E, \quad d(x, y) \!  \leq \! d(x, S)+d(y, S) \; \Longrightarrow \; d(x,y)\! = \! d_{\Ptt, \epp} (x,y). 
\end{equation}
Then, for all $r\ino (0, \infty)$, denote by $B_d (x, r)$ (resp.~by $B_{d_{\Ptt, \epp}} (\varpi_{\Ptt , \epp} (x), r)$) the open ball in $(E, d)$ (resp.~in $(E/\!\! \equiv_{\Ptt, \epp}, d_{\Ptt, \epp})$) with center $x$ (resp.~$\varpi_{\Ptt , \epp} (x)$) and radius $r$. 
Then, (\ref{populonia}) entails the following: if $x\ino E\backslash S$ and if $0\! < \! r \! < \!  \frac{1}{4} d(x, S))$, then 
\begin{equation}
\label{Perouse}
\varpi_{\Ptt, \epp} \! :\!  B_d (x, r)\! \rightarrow \!  B_{d_{\Ptt, \epp}} (\varpi_{\Ptt , \epp} (x), r) \; \textrm{is a surjective isometry.}
\end{equation}
Namely, outside the pinching points, the metric is locally unchanged. 

We now prove a result on Hausdorff and packing dimensions that is used in the proof of Proposition \ref{fractCIRG}. To that end, we suppose that there exists $(E_0, d)$, a compact metric space such that 
$E \! \subset E_0$ and such that $E$ is a compact subset of $E_0$. To simplify notation we set 
$(E^\prime, d^\prime, \varpi)\! = \! \big( E/\!\! \equiv_{\Ptt, \epp}, d_{\Ptt, \epp}, \varpi_{\Ptt, \epp}\big)$.  We denote by $\mathrm{dim}_{\mathrm{H}} $ and $\mathrm{dim}_{\mathrm{p}}$ resp.~the Hausdorff and the packing dimensions.

\begin{lem}   
\label{dimpinch} We keep the notations from above. We first 
assume that $\mathrm{dim}_{\mathrm{H}} (E_0)\ino (0, \infty) $ and $\mathrm{dim}_{\mathrm{p}} (E_0)\ino (0, \infty)$. 
Let $a\ino (0, \mathrm{dim}_{\mathrm{H}} (E_0))$ and $b\ino (0, \mathrm{dim}_{\mathrm{p}} (E_0))$; we 
assume that there is a finite measure $m_0$ on the Borel subsets of $E_0$ such that 
$m_0(E)\! >\! 0$ and 
\begin{equation}
\label{beffiage}
 \textrm{for $m_0$-almost all $x\ino E_0$} \; \, \limsup_{r\rightarrow 0} \frac{m_0(B_d(x, r))}{r^a} \! < \! \infty \; \, \textrm{and} \; \,   
\liminf_{r\rightarrow 0} \frac{m_0(B_d(x, r))}{r^b} \! < \! \infty . 
\end{equation}
Then, $a \! \leq \!  \mathrm{dim}_{\mathrm{H}} (E^\prime) \! \leq \! \mathrm{dim}_{\mathrm{H}} (E_0)$ 
and $b \! \leq \!  \mathrm{dim}_{\mathrm{p}} (E^\prime) \! \leq \! \mathrm{dim}_{\mathrm{p}} (E_0)$. 
\end{lem}
\textbf{Proof.} Since $\varpi$ is Lipschitz, 
$\mathrm{dim}_{\mathrm{H}} (E^\prime)\! \leq \! \mathrm{dim}_{\mathrm{H}} (E)\! \leq \! \mathrm{dim}_{\mathrm{H}} 
(E_0)$, with the same inequality for packing dimensions. We set $m\! = \! m_0 (\cdot \cap  E)$ and $m^\prime\! = \! 
m\circ \varpi^{-1}$ that is the pushforward measure of  $m$ via $\varpi$. Since $m(E)\! >\! 0$, (\ref{beffiage}) holds true with $m_0$
replaced by $m$. Observe that (\ref{beffiage}) implies that $m_0$ has no atom. Thus, $m$ has no atom and 
since there is a finite number of pinching points, (\ref{Perouse}) entails that (\ref{beffiage}) holds true for $m^\prime$ which 
entails $\mathrm{dim}_{\mathrm{H}} (E^\prime) \! \geq \! a$ and $\mathrm{dim}_{\mathrm{p}} (E^\prime) \! \geq \! b$ by standard 
comparison results on Hausdorff and packing measures due to Rogers \& Taylor in \cite{RoTa61} (Hausdorff case) 
and Taylor \& Tricot in \cite{TaTr85} (packing case) in Euclidian spaces  that have been extended in Edgar \cite{Ed07} 
(see Theorem 4.15 and Proposition 4.24 for the Hausdorff case and see Theorem 5.9 for the packing case). \cqfd

{\small
\setlength{\bibsep}{.3em}
\bibliographystyle{acm}
\bibliography{Refs}
}

\end{document}